\documentclass{conm-p-l}

\usepackage{amssymb,amsmath,cite}

\copyrightinfo{2008}{International Press}

\usepackage{epsfig}
\usepackage{color,psfrag}
\usepackage{graphics}
\usepackage{graphicx}

\newcommand{\N}{{\mathbb N}}

\newcommand{\C}{{\mathbb C}}
\newcommand{\R}{{\mathbb R}}

\newcommand{\be}{\begin{equation}}
\newcommand{\ee}{\end{equation}}

\newcommand{\irn}{\int_{\R^2}}
\newcommand{\figurescale}{0.41}

\newcommand{\eps}{\hbar}

\renewcommand{\theta}{{\vartheta}}

\numberwithin{equation}{section}
\newtheorem{theorem}{Theorem}[section]
\newtheorem{proposition}[theorem]{Proposition}

\newtheorem{lemma}[theorem]{Lemma}

\theoremstyle{definition}
\newtheorem{remark}[theorem]{Remark}

\numberwithin{equation}{section}

\begin{document}

\title[Ground and excited states for 2D Gross--Pitaevskii systems]{Location and phase segregation of ground
and \\ excited states for 2D Gross--Pitaevskii systems}

\author{Marco Caliari}

\curraddr{Dipartimento di Informatica,
Universit\`a di Verona,
C\`a Vignal 2, Strada Le Grazie~15, I-37134 Verona, Italy.}
\email{marco.caliari@univr.it}

\author{Marco Squassina}

\curraddr{Dipartimento di Informatica,
Universit\`a di Verona,
C\`a Vignal 2, Strada Le Grazie~15, I-37134 Verona, Italy.}
\email{marco.squassina@univr.it}

\thanks{The research of the second author was partially supported by the MIUR national research
project {\it Variational and Topological Methods in the Study of
Nonlinear Phenomena}}


\dedicatory{Communicated by Charles Li, received December 4, 2007 \\
and, in revised form, February 20, 2008.}


\begin{abstract}
We consider a system of Gross--Pitaevskii equations in $\R^2$ modelling a mixture of two Bose--Einstein condensates
with repulsive interaction. We aim to study the qualitative behaviour of ground and excited state solutions.
We allow two different harmonic and off-centered trapping potentials and study 
the spatial patterns of the solutions within the Thomas--Fermi approximation as well as phase segregation phenomena 
within the large-interaction regime.
\end{abstract}
\maketitle

\tableofcontents

\section{Introduction}
The first successful experimental realization of Bose--Einstein condensates for atomic gases \cite{BecR}, 
which goes back to 1995, gave rise to various  
numerical and theoretical investigations on the macroscopic equation ruling these phenomena, that
is the Gross--Pitaevskii equation (GPE)
\begin{equation*}
\eps i\partial_t\psi=-\frac{\eps^2}{2m}\Delta\psi+V(x)\psi+\theta|\psi|^2\psi, \quad \text{$x\in \R^n$},\quad\theta\geq 0.
\end{equation*}
More recently, in 1997, Bose--Einstein condensation for a mixture 
of two different interacting atomic species with the same mass
was firstly realized at JILA \cite{BecSR} exhibiting a partial
overlap between the wave functions.  The vector nature of the order parameter 
gives rise to some intriguing structures and dynamics that are absent in the single component case.
This, again, stimulated various succeeding studies of numerical and theoretical nature.
For both single or binary condensates we refer the reader to \cite{Bigpap,PiTaBook} and to the references therein. 
Recently, some efficient numerical techniques have been developed to compute 
ground state solutions of GPE \cite{Baoal,CalThal}, which can be used to investigate the vector case. 
On this basis, in this paper we deal with the rigorous analysis of the spatial configurations for the 
standing wave solutions (ground and excited states) of the system in $\R^2$
\begin{equation}
\label{GPs}
\begin{cases}
\eps i\partial_t\psi_1=-\frac{\eps^2}{2m_1}\Delta \psi_1
+V_1(x_1,x_2)\psi_1+\theta_{11}\eps^2|\psi_1|^2\psi_1+\theta_{12}\eps^2|\psi_2|^2\psi_1,  \\
\noalign{\vskip6pt}
\eps  i\partial_t\psi_2=-\frac{\eps^2}{2m_2}\Delta \psi_2
+V_2(x_1,x_2)\psi_2+\theta_{21}\eps^2|\psi_1|^2\psi_2+\theta_{22}\eps^2|\psi_2|^2\psi_2,  
\end{cases}
\end{equation}
for the unknown $\psi_i:\R^2\to\C$, $i=1,2$, where $\eps$ denotes the (reduced) Planck constant and the coefficients
$\theta_{ij}\geq 0$ (defocusing case), with $\theta_{12}=\theta_{21}$, are given by the formula \cite{esrygreene} 
$$
\theta_{ij}=2\pi\sigma_{ij} \frac{m_i+m_j}{m_im_j},\quad \sigma_{ij}=\sigma_{ji}, \quad i,j=1,2,
$$
being $\sigma_{ij}$ related to the scattering lengths and $m_i$ the atomic masses of the two species 
composing the mixture. Considering the 2D case is not restrictive as there are
various situations where the full 3D system can be reduced to a 2D system with suitably modified coefficients
(see e.g.\ Section 2.2 of \cite{Babao}).
The coefficients $\theta_{ii}$ and  $\theta_{12}$ play the role of repulsive intra-species and inter-species parameters
respectively. As we will see, when $\theta_{12}$ is sufficiently large, then some interesting overlap and spatial
segregation phenomena between the wave densities occur. Concerning the 
potentials, we let $V_i(x_1,x_2)=\frac{m_i}{2}(\omega_{i1}^2(x_1-x_{i1})^2+\omega_{i2}^2(x_2-x_{i2})^2)$ for $(x_1,x_2)$ in $\R^2$,
where $\omega_{ij}> 0$, $i,j=1,2$. A typical situation is when the $V_i$s have the same 
center, without loss of generality, the origin. On the other hand, there are some physical situations reported in literature, which lead to consider
off-centered potentials. See, for instance, \cite{ribolimod}, where the vertical direction  in the potential
is not aligned with the symmetry axis of the trap.
Similar  equations have also arisen as governing equations for electromagnetic pulse propagation in ``left-handed'' materials 
with Kerr-type nonlinearity \cite{LazTsi},  in the modified Hubbard model in 
the long-wavelength approximation \cite{Mak,Lind}, in quadratic nonlinear materials with
suitable phase matching \cite{Kaloc} and in nonlinear optics,
for instance in the pro\-pagation of pulses in a nonlinear optical fiber of bi-modal type
due to the presence of some birefringence effects generating
two pulses with different polarization directions \cite{men}. For a wide discussion on 
nonlin\-ear Schr\"o\-din\-ger systems we refer the interested reader
to \cite{abpr,abl2,abl3} and to the references therein.
\vskip2pt
Let ${\mathcal H}$ be the Hilbert subspace of 
$H^1(\R^2,\C)\times H^1(\R^2,\C)$ defined by
$$
{\mathcal H}=\left\{(\psi_1,\psi_2)\in H^1(\R^2,\C)\times H^1(\R^2,\C):\,\,\,\irn V_i(x_1,x_2)|\psi_i|^2<\infty,\,\,\,i=1,2\right\},
$$
which is the natural framework for bound state solutions, endowed with the norm
$$
\|(\psi_1,\psi_2)\|_{{\mathcal H}}^2=\sum_{i=1}^2\frac{\eps^2}{2m_i}\irn |\nabla \psi_i|^2+\irn V_i(x_1,x_2)|\psi_i|^2,
$$
and consider the total energy functional $E:{\mathcal H}\to\R$ associated with \eqref{GPs}
\begin{equation}
\label{totenergf}
E(\psi_1(\cdot,t),\psi_2(\cdot,t))=\sum_{i=1}^2 E_i(\psi_i(\cdot,t))
+\theta_{12}\eps^2\irn |\psi_1(\cdot,t)|^2|\psi_2(\cdot,t)|^2, 
\end{equation}
where, for $i=1,2$, we set
\begin{equation*}
E_i(\psi_i(\cdot,t))=\frac{\eps^2}{2m_i}\irn |\nabla 
\psi_i(\cdot,t)|^2+\irn V_i(x_1,x_2)|\psi_i(\cdot,t)|^2+\frac{\theta_{ii}\eps^2}{2}\irn |\psi_i(\cdot,t)|^4.
\end{equation*}
By multiplying the first equation of system \eqref{GPs} by $\partial_t\bar\psi_1$
and the second by $\partial_t\bar\psi_2$, taking the real parts, integrating and adding the resulting identities, it is readily seen that 
$E$ is constant on the solutions, namely $E(\psi_1(\cdot,t),\psi_2(\cdot,t))=E(\psi_1(\cdot,0),\psi_2(\cdot,0))$
for any $t\geq 0$.  Also, as for the case of the single equation,  by multiplying the first equation of \eqref{GPs} 
by $\bar\psi_1$ and the second by $\bar\psi_2$, taking the imaginary parts, integrating and adding the resulting identities, 
it turns out that the total number of particles $N_i$ of the $i$-th species 
is time independent (preservation of the particle number), namely
\begin{equation}
\label{totmasses}
\irn |\psi_i(\cdot,t)|^2=N_i,\qquad t\geq 0,\,\, i=1,2.
\end{equation}
The {\em ground state} (or least energy) solution of \eqref{GPs} is a solution with ansatz 
\begin{equation}
\label{ansatz}
\psi_i(x_1,x_2,t)=e^{-{\rm i}\frac{\mu_it}{\eps}}\phi_i(x_1,x_2),\qquad (x_1,x_2)\in\R^2,\,\, t\geq 0,\,\,\, i=1,2,
\end{equation}
where the pair $(\phi_1,\phi_2)$ is real valued and minimizes functional \eqref{totenergf}
constrained to conditions \eqref{totmasses} (with $\phi_i$ in place of $\psi_i$).  
As a consequence, the $\phi_i$s solve the nonlinear eigenvalue problem in $\R^2$
\begin{equation}
\label{systemGPGen}
\begin{cases}
-\frac{\eps^2}{2m_1}\Delta \phi_1+V_1(x_1,x_2)\phi_1+\theta_{11}\eps^2|\phi_1|^2\phi_1
+\theta_{12}\eps^2|\phi_2|^2\phi_1=\mu_1\phi_1, \\
\noalign{\vskip6pt}
-\frac{\eps^2}{2m_2}\Delta \phi_2+V_2(x_1,x_2)\phi_2+\theta_{21}\eps^2|\phi_1|^2\phi_2+\theta_{22}\eps^2|\phi_2|^2\phi_2=\mu_2\phi_2, \\
\noalign{\vskip6pt}
\,\,\,\displaystyle\irn \phi_1^2=N_1,\,\,\,\irn \phi_2^2=N_2.
\end{cases}
\end{equation}
Testing the first equation of  \eqref{systemGPGen} by $\bar\psi_1$
and the second by $\bar\psi_2$, we have a formula for the eigenvalues $\mu_i$ 
(also known as chemical potentials) versus the eigenvectors $\phi_i$
\begin{equation}
\label{eig-systemGPGen}
N_i\mu_i=E_i(\phi_i)+\frac{\theta_{ii}}{2}\eps^2\irn |\phi_i|^4+\theta_{12}\eps^2\irn |\phi_1|^2|\phi_2|^2,\qquad i=1,2.
\end{equation}
The existence of ground state solutions to \eqref{GPs} in $ {\mathcal H}$ is reduced to the existence of minima for the energy functional
\eqref{totenergf} constrained to  
\begin{equation}
\label{sfera}
{\mathcal S}=\{(\phi_1,\phi_2)\in {\mathcal H}:\,\|\phi_i\|_{L^2}^2=N_i,\,\,\text{$i=1,2$}\}.
\end{equation}
As the $\theta_{ij}$ are positive, 
$$
\frac{\theta_{11}}{2}|\phi_1|^4+\theta_{12}|\phi_1|^2|\phi_2|^2
+\frac{\theta_{22}}{2}|\phi_2|^4\geq 0,
$$
so that the energy functional $E$ is coercive, 
bounded from below and weakly lower semi-continuous over ${\mathcal S}$. Hence, the 
existence of a ground state solution is immediately guaranteed. 
Usually, with reference to the solutions of the form \eqref{ansatz} (standing waves),
there are two possible (physically different) approaches depending on whether one considers the 
chemical potentials $\mu_i$ as fixed (hence searching for solutions to the first two equations in
\eqref{systemGPGen} but with possibly different $L^2$ norms)
or the total masses $\irn |\phi_i|^2$ as fixed (thus solving the nonlinear eigenvalue problem
\eqref{systemGPGen}, which is the case we deal with). Any other solution $(\phi_1,\phi_2)$ of system \eqref{systemGPGen}
of the form \eqref{ansatz} not having minimal energy for $E$ will be called {\em excited state} (or higher energy solution).
\vskip1pt
The main goal of this paper is to prove some geometrical properties (clearly confirmed by some numerical simulations) 
for ground and excited state of \eqref{systemGPGen}, particularly under the 
influence of strong interaction effects (namely $\theta_{12}\to\infty$). See 
e.g.\ Propositions \ref{Gsspatseg} and \ref{finalPropES}.
\vskip1pt
In Section~\ref{TFapp} we derive
the location of ground state solutions via the Thomas--Fermi approximation and classify the relative
configuration of $\phi_i$ with respect to $\phi_j$. In Section~\ref{SIr} we study the phase separation process
(spatial segregation) in the large competition regime by means of suitable limiting energy levels which provide 
$\kappa$-independent upper bounds for the energy of solutions. 
In Section~\ref{Numset}, for the sake of completeness, we briefly describe the functional framework of 
the numerical scheme used to compute the solutions.

In the following the Hilbert space $L^2(\R^2,\C)$ is endowed with the standard scalar product  
$(f,g)_2=\irn f\bar g$, $f,g\in L^2(\R^2)$ and the induced norm is denoted by $\|\cdot\|_{L^2}$.

\medskip
\section{Location and Thomas--Fermi approximation}
\label{TFapp}

If the distance between the centers of the trapping potentials $V_i$ 
is sufficiently small compared with the radii of the supports of the ground state solutions $\phi_i$,
then the condensates share a region where they coexist (with a partial or full overlap, that is
one condensate is partially or entirely included in the other). In the opposite case the supports
of the wave functions are disjoint.
Hence, we can encounter three different patterns
for the spatial wave functions $\phi_i$, which we are going to 
discuss, namely: {\em no} overlap, {\em partial} overlap and {\em full} overlap.
It should be noted that {\em support} just means here the planar region where the mass of the ground state solution is mainly 
concentrated, being (exponentially) vanishing on the outside.
In the Thomas--Fermi regime, an approximation of the ground state solutions of system~\eqref{systemGPGen},
which is very good for sufficiently large values of the coupling constants, can be obtained by simply
dropping the diffusion terms $-\Delta\phi_i$, namely the kinetic contributions, thus assuming the wave functions
to be slowly varying (cf.\ \cite{fermi,thomas,lieb,liebsimon}). In turn, \eqref{systemGPGen} reduces to 
the algebraic system (here we let $\eps=1$)
\begin{equation}
\label{systemGPGenTF}
\begin{cases}
2\theta_{11}|\phi_1|^2+2\theta_{12}|\phi_2|^2=2\mu_1-(x_1-x_{11})^2-(x_2-x_{12})^2, \\
\noalign{\vskip5pt}
2\theta_{21}|\phi_1|^2+2\theta_{22}|\phi_2|^2=2\mu_2-(x_1-x_{21})^2-(x_2-x_{22})^2, 
\end{cases}
\end{equation}
where the $\mu_i$s should be computed through the normalization conditions \eqref{totmasses} 
(if, for instance, $\theta_{12}=0$,  it holds $\mu_i\propto\sqrt{\theta_{ii}}$ for $i=1,2$).
In general, as the left-hand sides are positive, this system is satisfied only on a (possibly empty) 
subset ${\mathcal O}\subset\R^2$ (${\mathcal O}={\mathcal O}_1\cap
{\mathcal O}_2$ in the notations introduced below), namely the 
overlap region. It is natural to introduce the circumferences defined by
$(x_1-x_{i1})^2+(x_2-x_{i2})^2=r_i^2$, with $r_i(\mu_i)=\sqrt{2\mu_i}$, $i=1,2$.
The intersection of the corresponding disks ${\mathcal D}_i$ gives the region where 
system \eqref{systemGPGenTF} makes sense. Outside the region ${\mathcal O}$, the wave functions take
the usual form of the solutions of the GPE ($\theta_{12}=0$) in the Thomas--Fermi regime (see the expressions below of $\phi_i$
on ${\mathcal D}_i\setminus{\mathcal O}$). More precisely, we have the following non-smooth
approximations of the ground state (see also the work by Riboli and Modugno \cite{ribolimod})
\begin{equation*}
\phi_1 =
\begin{cases}
\sqrt{\frac{\theta_{22}(2\mu_1-(x_1-x_{11})^2-(x_2-x_{12})^2)
- \theta_{12}(2\mu_2-(x_1-x_{21})^2-(x_2-x_{22})^2)}{2(\theta_{11}\theta_{22}-\theta_{12}^2)}}, & \text{in ${\mathcal O}$}, \\
\noalign{\vskip3pt}
\sqrt{\frac{2\mu_1-(x_1-x_{11})^2-(x_2-x_{12})^2}{2\theta_{11}}}, & \text{in ${\mathcal D}_1\setminus{\mathcal O}$}, \\
\noalign{\vskip7pt}
0, & \text{in $\R^2\setminus{\mathcal D}_1$},
\end{cases}
\end{equation*}
\begin{equation*}
\phi_2=
\begin{cases}
\sqrt{\frac{\theta_{11}(2\mu_2-(x_1-x_{21})^2-(x_2-x_{22})^2)
- \theta_{12}(2\mu_1-(x_1-x_{11})^2-(x_2-x_{12})^2)}{2(\theta_{11}\theta_{22}-\theta_{12}^2)}}, & \text{in ${\mathcal O}$}, \\
\noalign{\vskip3pt}
\sqrt{\frac{2\mu_2-(x_1-x_{21})^2-(x_2-x_{22})^2}{2\theta_{22}}}, & \text{in ${\mathcal D}_2\setminus{\mathcal O}$}, \\
\noalign{\vskip7pt}
0, & \text{in $\R^2\setminus {\mathcal D}_2$},
\end{cases}
\end{equation*}
that is, equivalently,
\begin{equation}
\label{TFrappppree}
\phi_i=
\begin{cases}
\sqrt{\alpha_i(R_i^2-(x_1-y_{i1})^2-(x_2-y_{i2})^2)}, & \text{in ${\mathcal O}$}, \\
\noalign{\vskip3pt}
\sqrt{\frac{r_i^2-(x_1-x_{i1})^2-(x_2-x_{i2})^2}{2\theta_{ii}}}, & \text{in ${\mathcal D}_i\setminus{\mathcal O}$}, \\
\noalign{\vskip6pt}
0, & \text{in $\R^2\setminus {\mathcal D}_i$},
\end{cases}
\end{equation}
where, according to the notations introduced below, ${\mathcal O}={\mathcal O}_1\cap {\mathcal O}_2$,
${\mathcal D}={\mathcal D}_1\cap {\mathcal D}_2$,
\begin{align*}
{\mathcal O}_i &=\big\{(x_1,x_2)\in{\mathcal D}:\,\,(x_1-y_{i1})^2+(x_2-y_{i2})^2\leq (\geq) R_i^2\big\},\\
{\mathcal D}_i&=\big\{(x_1,x_2)\in\R^2:\,\,(x_1-x_{i1})^2+(x_2-x_{i2})^2\leq r_i^2\big\},
\end{align*}
with $\leq$ (resp.\ $\geq$) in the definition of ${\mathcal O}_i$  for $\alpha_i>0$ (resp.\ $\alpha_i<0$),
if  we set
\begin{align*}
y_{11}&=\frac{\omega_{22}x_{11}-\omega_{12} x_{21}}{\omega_{22}-\omega_{12}}=
x_{11}+\frac{\omega_{12}}{\alpha_1}\Delta_1x,  
\quad
y_{12}=\frac{\omega_{22}x_{12}-\omega_{12} x_{22}}{\omega_{22}-\omega_{12}} =
x_{12}+\frac{\omega_{12}}{\alpha_1}\Delta_2x, \\
\noalign{\vskip5pt}
y_{21}&=\frac{\omega_{11}x_{21}-\omega_{12} x_{11}}{\omega_{11}-\omega_{12}}=
x_{21}-\frac{\omega_{12}}{\alpha_2}\Delta_1x, 
\quad
y_{22}=\frac{\omega_{11}x_{22}-\omega_{12} x_{12}}{\omega_{11}-\omega_{12}}=
x_{22}-\frac{\omega_{12}}{\alpha_2}\Delta_2x, 
\end{align*}
where $\Delta_jx=x_{1j}-x_{2j}$, namely, for $i,j=1,2$ with $i\neq j$,
$$
y_{ij}=x_{ij}-(-1)^i\frac{\omega_{12}}{\alpha_i}\Delta_jx, 
$$
with $\omega_{ij}=\frac{\theta_{ij}}{2{\rm det}\,\Theta}$, $i,j=1,2$,
$\alpha_1=\omega_{22}-\omega_{12}$, $\alpha_2=\omega_{11}-\omega_{12}$, and
\begin{align*}
R_1^2&=\frac{2\omega_{22}\mu_1-2\omega_{12}\mu_2+\omega_{12}x_{21}^2+
\omega_{12} x_{22}^2-\omega_{22} x_{11}^2-\omega_{22} x_{12}^2}{\omega_{22}-\omega_{12}} +y_{11}^2+y_{12}^2 \\
\noalign{\vskip6pt}
R_2^2&=\frac{2\omega_{11}\mu_2-2\omega_{12}\mu_1+\omega_{12}x_{11}^2+
\omega_{12} x_{12}^2-\omega_{11} x_{21}^2-\omega_{11} x_{22}^2}{\omega_{11}-\omega_{12}}  +y_{22}^2+y_{21}^2.
\end{align*}
Setting ${\bf x_i}=(x_{i1},x_{i2})$ and ${\bf y_i}=(y_{i1},y_{i2})$ for $i=1,2$, this reads as
$$
R_i=\sqrt{r_i^2+\frac{2\omega_{12}}{\alpha_i}(\mu_i-\mu_j)
+\frac{\omega_{12}}{\alpha_i}(|{\bf x_j}|^2-|{\bf y_i}|^2)
-\frac{\omega_{jj}}{\alpha_i}(|{\bf x_i}|^2-|{\bf y_i}|^2)}.
$$
Hence, we have four circumferences that rule the geometry of the ground states
\begin{equation*}
\Sigma_i^r:\, (x_1-x_{i1})^2+(x_2-x_{i2})^2=r_i^2,  \qquad
\Sigma_i^R:\, (x_1-y_{i1})^2+(x_2-y_{i2})^2= R_i^2.
\end{equation*}
If $\theta_{12}=0$ (namely no interaction), we deduce that ${\bf x_i}={\bf y_i}$, $R_i=r_i$, and
${\mathcal O}_i={\mathcal D}_i$ ($\theta_{12}=0$ implies $\alpha_i>0$), so that the ground 
state solutions turn into the usual Thomas--Fermi representation for the single GPE
\begin{equation*}
\phi_i=
\begin{cases}
\sqrt{\frac{r_i^2-(x_1-x_{i1})^2-(x_2-x_{i2})^2}{2\theta_{ii}}}, & \text{in ${\mathcal D}_i$}, \\
\noalign{\vskip4pt}
0, & \text{in $\R^2\setminus {\mathcal D}_i$}.
\end{cases}
\end{equation*}
If $\theta_{12}\approx 0$, then
$\omega_{11}\approx\frac{1} {2\theta_{22}}$,
$\omega_{22}\approx\frac{1} {2\theta_{11}}$,
$\omega_{12}\approx 0$,
$\alpha_i\approx\frac{1} {2\theta_{ii}}$ and
$y_{ij}\approx x_{ij}$, $R_i\approx r_i$ for $i,j=1,2$,
so that $\Sigma_i^R\approx \Sigma_i^r$ for $i=1,2$.
If ${\bf x_1}={\bf x_2}$, the $\Sigma_i^R$s have centers 
${\bf y_i}={\bf x_i}$ but different radii
$R_i^2=r_i^2+\frac{2\omega_{12}}{\alpha_i}(\mu_i-\mu_j)$, for any $i\neq j$.

\subsection{Nonoverlap case}
\begin{figure}[h!!!]
\begin{center}
\begin{psfrags}%
\psfragscanon%
%
\psfrag{s10}[][]{\color[rgb]{0,0,0}\setlength{\tabcolsep}{0pt}\begin{tabular}{c} \end{tabular}}%
\psfrag{s11}[][]{\color[rgb]{0,0,0}\setlength{\tabcolsep}{0pt}\begin{tabular}{c} \end{tabular}}%
\psfrag{s12}[l][l]{\color[rgb]{0,0,0}R2}%
\psfrag{s13}[l][l]{\color[rgb]{0,0,0}${\scriptscriptstyle r_1}$}%
\psfrag{s14}[l][l]{\color[rgb]{0,0,0}${\scriptscriptstyle r_2}$}%
\psfrag{s15}[l][l]{\color[rgb]{0,0,0}${\scriptscriptstyle R_1}$}%
\psfrag{s16}[l][l]{\color[rgb]{0,0,0}${\scriptscriptstyle R_2}$}%
%
\psfrag{x01}[t][t]{-20}%
\psfrag{x02}[t][t]{-10}%
\psfrag{x03}[t][t]{0}%
\psfrag{x04}[t][t]{10}%
\psfrag{x05}[t][t]{20}%
\psfrag{x06}[t][t]{30}%
%
\psfrag{v01}[r][r]{-15}%
\psfrag{v02}[r][r]{-10}%
\psfrag{v03}[r][r]{-5}%
\psfrag{v04}[r][r]{0}%
\psfrag{v05}[r][r]{5}%
\psfrag{v06}[r][r]{10}%
\psfrag{v07}[r][r]{15}%
%
\includegraphics[scale=\figurescale]{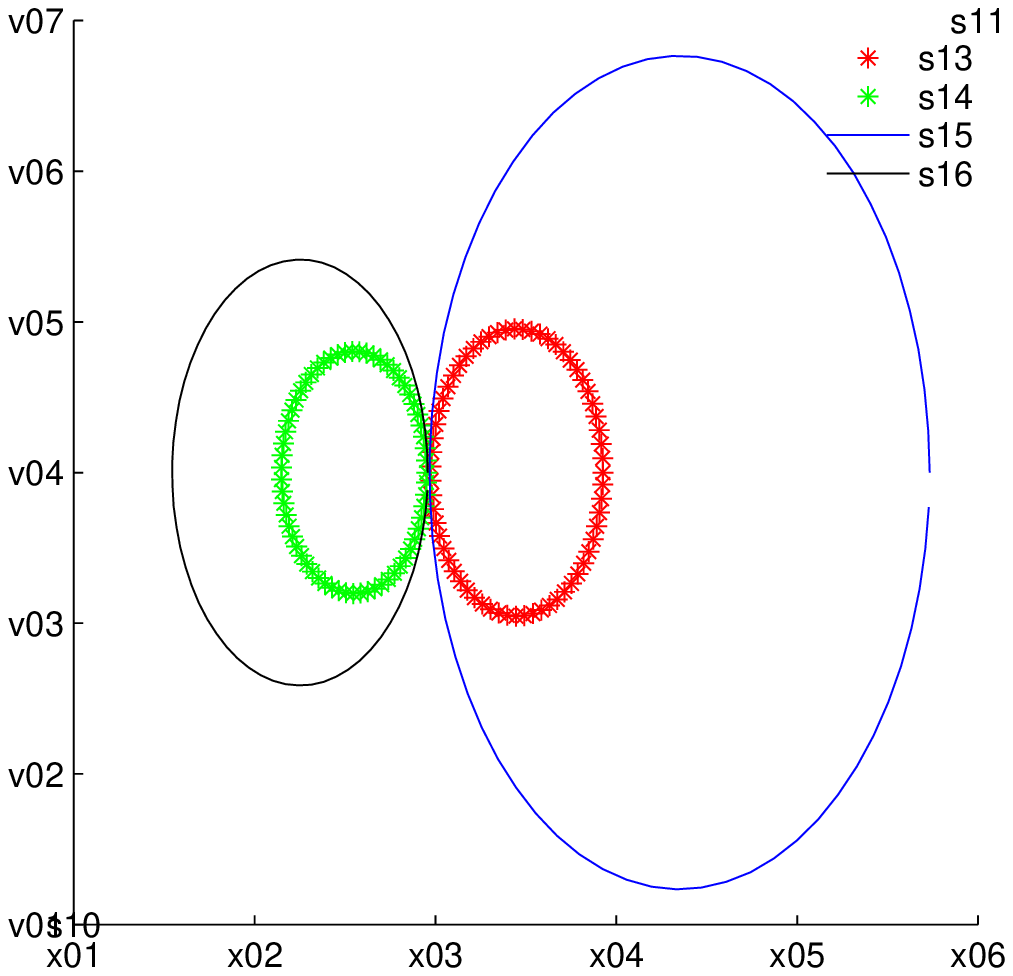}%
\end{psfrags}%
      \hspace{1cm}
\begin{psfrags}%
\psfragscanon%
%
\psfrag{s10}[][]{\color[rgb]{0,0,0}\setlength{\tabcolsep}{0pt}\begin{tabular}{c} \end{tabular}}%
\psfrag{s11}[][]{\color[rgb]{0,0,0}\setlength{\tabcolsep}{0pt}\begin{tabular}{c} \end{tabular}}%
\psfrag{s12}[l][l]{\color[rgb]{0,0,0}R2}%
\psfrag{s13}[l][l]{\color[rgb]{0,0,0}${\scriptscriptstyle r_1}$}%
\psfrag{s14}[l][l]{\color[rgb]{0,0,0}${\scriptscriptstyle r_2}$}%
\psfrag{s15}[l][l]{\color[rgb]{0,0,0}${\scriptscriptstyle R_1}$}%
\psfrag{s16}[l][l]{\color[rgb]{0,0,0}${\scriptscriptstyle R_2}$}%
%
\psfrag{x01}[t][t]{-20}%
\psfrag{x02}[t][t]{-10}%
\psfrag{x03}[t][t]{0}%
\psfrag{x04}[t][t]{10}%
\psfrag{x05}[t][t]{20}%
\psfrag{x06}[t][t]{30}%
\psfrag{x07}[t][t]{40}%
%
\psfrag{v01}[r][r]{-20}%
\psfrag{v02}[r][r]{-15}%
\psfrag{v03}[r][r]{-10}%
\psfrag{v04}[r][r]{-5}%
\psfrag{v05}[r][r]{0}%
\psfrag{v06}[r][r]{5}%
\psfrag{v07}[r][r]{10}%
\psfrag{v08}[r][r]{15}%
\psfrag{v09}[r][r]{20}%
%
\includegraphics[scale=\figurescale]{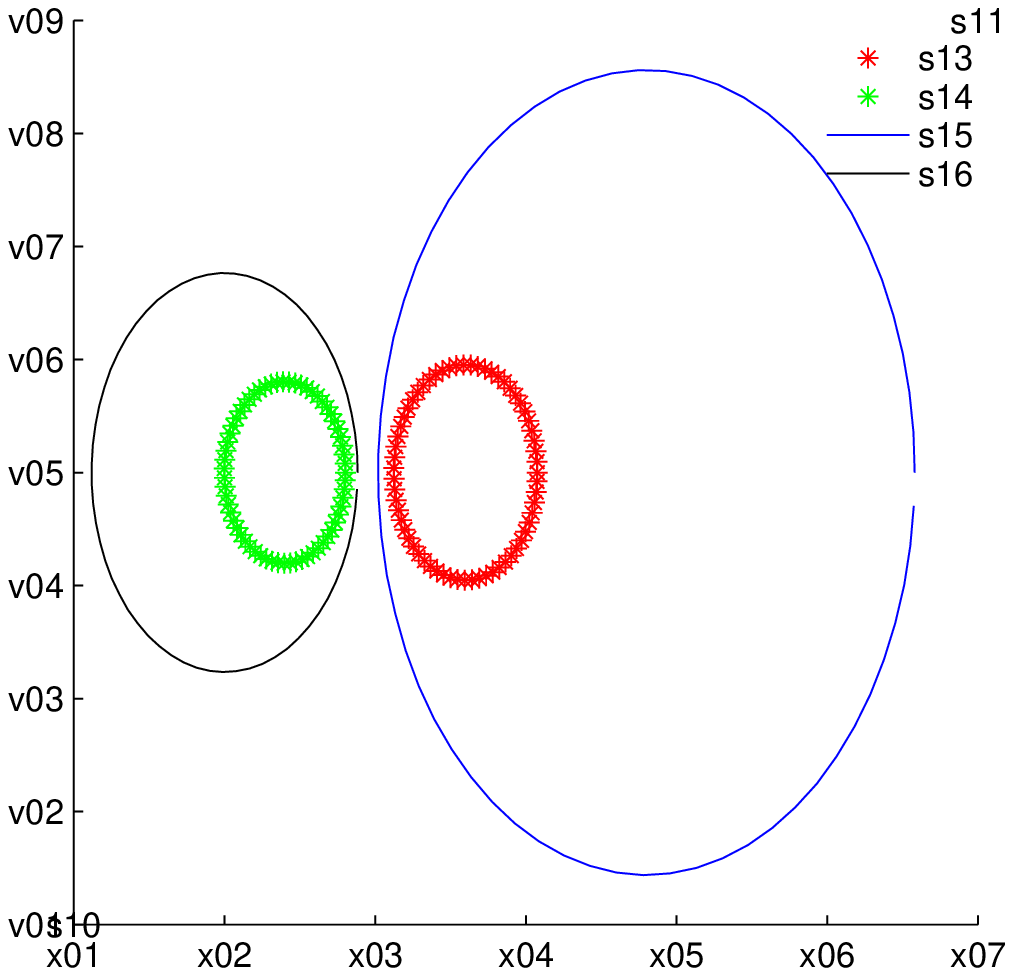}%
\end{psfrags}%
\caption{Supports of  $\phi_1$ and $\phi_2$ (starred) and disks (unstarred) related to the
Thomas--Fermi approximation for $\phi_1$ and $\phi_2$. 
We have taken  $N_1=N_2=m_1=m_2=1$, $\theta_{11}=400$, $\theta_{22}=200$, 
$\theta_{12}=100$, $x_{11}=4.5$, $x_{21}=-4.5$ (left figure, tangential supports);
and $x_{11}=6$, $x_{21}=-6$ (right figure, disjoint supports), with $x_{ij}=0$ for any other $i,j$.}
\label{supporti-position}
\end{center}
\end{figure}

In the case occurring when the constant $\theta_{12}$ is zero (absence of interaction in the mixture),  
the system uncouples into a pair of GPEs (for the single GPE various accurate and
efficient numerical techniques have been recently compared in \cite{CalThal}).
Numerical experiments show that the ground state solution
$\phi_i$ always locates its mass around the minimum point ${\bf x_i}=(x_{i1},x_{i2})$ of $V_i$, for $i=1,2$.
It looks apparent that boosting up the parameter $\theta_{ii}$ in front of the cubic nonlinearity
in the equation of $\phi_i$ has the effect of squeezing down the profile of $\phi_i$
making it flatter and larger. Going back to the case $\theta_{12}\neq 0$, we say that we have 
{\em no overlap} between $\phi_1$ and $\phi_2$, if the centers
${\bf x_1}$ and ${\bf x_2}$ of $\Sigma_i^r$ satisfy the geometric condition
\begin{equation}
\label{non-overl}
(\Delta_1x)^2+(\Delta_2x)^2>
|\sqrt{2\mu_1}+\sqrt{2\mu_2}|^2,
\end{equation}
namely if ${\bf x_1}$ is sufficiently far from ${\bf x_2}$ with respect to the amplitudes $r_i$
of the supports of $\phi_i$. In this situation the ground state solutions
look like those of the decoupled case. In fact, the coupling
terms $C_{ij}=\theta_{12}|\phi_j|^2\phi_i$ with $i\neq j$ are almost everywhere zero as
the supports are disjoint, due to \eqref{non-overl}. Hence the system is actually a small deformation of a pair
of uncoupled GPEs. See Figure \ref{supporti-position}.

\subsection{Partial overlap case}
\begin{figure}[h!!!]
\begin{center}
\begin{psfrags}%
\psfragscanon%
%
\psfrag{s10}[][]{\color[rgb]{0,0,0}\setlength{\tabcolsep}{0pt}\begin{tabular}{c} \end{tabular}}%
\psfrag{s11}[][]{\color[rgb]{0,0,0}\setlength{\tabcolsep}{0pt}\begin{tabular}{c} \end{tabular}}%
\psfrag{s12}[l][l]{\color[rgb]{0,0,0}R2}%
\psfrag{s13}[l][l]{\color[rgb]{0,0,0}${\scriptscriptstyle r_1}$}%
\psfrag{s14}[l][l]{\color[rgb]{0,0,0}${\scriptscriptstyle r_2}$}%
\psfrag{s15}[l][l]{\color[rgb]{0,0,0}${\scriptscriptstyle R_1}$}%
\psfrag{s16}[l][l]{\color[rgb]{0,0,0}${\scriptscriptstyle R_2}$}%
%
\psfrag{x01}[t][t]{-10}%
\psfrag{x02}[t][t]{-5}%
\psfrag{x03}[t][t]{0}%
\psfrag{x04}[t][t]{5}%
\psfrag{x05}[t][t]{10}%
\psfrag{x06}[t][t]{15}%
%
\psfrag{v01}[r][r]{-8}%
\psfrag{v02}[r][r]{-6}%
\psfrag{v03}[r][r]{-4}%
\psfrag{v04}[r][r]{-2}%
\psfrag{v05}[r][r]{0}%
\psfrag{v06}[r][r]{2}%
\psfrag{v07}[r][r]{4}%
\psfrag{v08}[r][r]{6}%
\psfrag{v09}[r][r]{8}%
%
\includegraphics[scale=\figurescale]{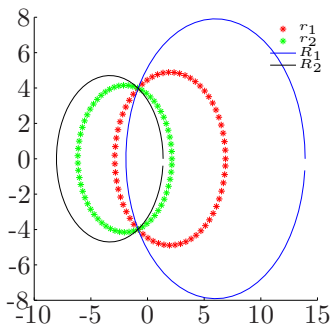}%
\end{psfrags}%
\caption{Supports of  $\phi_1$ and $\phi_2$ (starred) and disks (unstarred) related to the
Thomas--Fermi approximation for $\phi_1$ and $\phi_2$ (partial overlap case). 
We have $N_1=N_2=m_1=m_2=1$, $\theta_{11}=400$, $\theta_{22}=200$, 
$\theta_{12}=100$, $x_{11}=2$, $x_{21}=-2$ and $x_{ij}=0$ for other $i,j$.}
\label{supporti-position2}
\end{center}
\end{figure}

We have {\em partial overlap} between $\phi_1$ and $\phi_2$, if 
$$
|\sqrt{2\mu_1}-\sqrt{2\mu_2}|^2<(\Delta_1x)^2+(\Delta_2x)^2<
|\sqrt{2\mu_1}+\sqrt{2\mu_2}|^2,
$$
namely the disks of boundaries $\Sigma_1^r$ and $\Sigma_2^r$ overlap without being
one completely embedded in the other. We see in the contour plot of Figure \ref{GS2partialOverlp-bis} the partial overlap
in the case $\theta_{11}\gg\theta_{22}$. Apparently, boosting $\theta_{11}$ with respect
to $\theta_{22}$ makes the overlap region more localized. 
Keeping in mind the behaviour of the uncoupled case, in order to give this fact
a very simple empirical explanation, it suffices to argue on the coupling terms $C_{ij}$.
In the region (depending upon the relative magnitude of the $\theta_{ii}$s) 
where both $\phi_i$ are nonzero the contribution of $C_{ij}$ pushes down the profile around the 
origin (the center of trapping for $\phi_2$), provided that $\theta_{12}$ is significantly
large. The support of $\phi_1$ still remains
contractible, but the radial symmetry property of $\phi_1$ is broken (due to strong interaction).
See e.g.\ the situations reported in Figures \ref{supporti-position2} and \ref{GS2partialOverlp-bis}.
As Figure~\ref{TFlargeintfig} shows, while the Thomas--Fermi approximation disks $\Sigma_i^R$ are overlapped to the
support disks $\Sigma_i^r$ when the coupling constant $\theta_{12}$ is much smaller than the $\theta_{ii}$s,
in the case where $\theta_{12}\gg\theta_{ii}$, i.e.\ in the large interaction regime, the four circumferences
intersect in two points and $\Sigma_i^r$ and $\Sigma_i^R$ may have quite different sizes.

\subsection{Full overlap case}
\begin{figure}[h!!!]
\begin{center}
\begin{psfrags}%
\psfragscanon%
%
\psfrag{s10}[][]{\color[rgb]{0,0,0}\setlength{\tabcolsep}{0pt}\begin{tabular}{c} \end{tabular}}%
\psfrag{s11}[][]{\color[rgb]{0,0,0}\setlength{\tabcolsep}{0pt}\begin{tabular}{c} \end{tabular}}%
\psfrag{s12}[l][l]{\color[rgb]{0,0,0}R2}%
\psfrag{s13}[l][l]{\color[rgb]{0,0,0}${\scriptscriptstyle r_1}$}%
\psfrag{s14}[l][l]{\color[rgb]{0,0,0}${\scriptscriptstyle r_2}$}%
\psfrag{s15}[l][l]{\color[rgb]{0,0,0}${\scriptscriptstyle R_1}$}%
\psfrag{s16}[l][l]{\color[rgb]{0,0,0}${\scriptscriptstyle R_2}$}%
%
\psfrag{x01}[t][t]{-10}%
\psfrag{x02}[t][t]{-5}%
\psfrag{x03}[t][t]{0}%
\psfrag{x04}[t][t]{5}%
\psfrag{x05}[t][t]{10}%
%
\psfrag{v01}[r][r]{-8}%
\psfrag{v02}[r][r]{-6}%
\psfrag{v03}[r][r]{-4}%
\psfrag{v04}[r][r]{-2}%
\psfrag{v05}[r][r]{0}%
\psfrag{v06}[r][r]{2}%
\psfrag{v07}[r][r]{4}%
\psfrag{v08}[r][r]{6}%
\psfrag{v09}[r][r]{8}%
%
\includegraphics[scale=\figurescale]{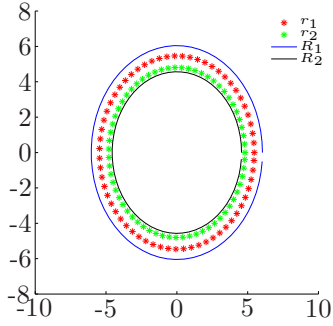}%
\end{psfrags}%
\caption{Supports of  $\phi_1$ and $\phi_2$ (inside of starred disks) and disks related to
the Thomas--Fermi approximation for $\phi_1$ and $\phi_2$ (overlap inside the smaller disk). 
We have taken $N_1=N_2=m_1=m_2=1$, $\theta_{11}=400$, $\theta_{22}=200$, 
$\theta_{12}=100$ and $x_{ij}=0$ for other $i,j=1,2$.}
\label{supporti-position3}
\end{center}
\end{figure}

\begin{figure}[h!!!]
\begin{center}
\begin{psfrags}%
\psfragscanon%
%
\psfrag{s10}[][]{\color[rgb]{0,0,0}\setlength{\tabcolsep}{0pt}\begin{tabular}{c} \end{tabular}}%
\psfrag{s11}[][]{\color[rgb]{0,0,0}\setlength{\tabcolsep}{0pt}\begin{tabular}{c} \end{tabular}}%
\psfrag{s12}[l][l]{\color[rgb]{0,0,0}R2}%
\psfrag{s13}[l][l]{\color[rgb]{0,0,0}${\scriptscriptstyle r_1}$}%
\psfrag{s14}[l][l]{\color[rgb]{0,0,0}${\scriptscriptstyle r_2}$}%
\psfrag{s15}[l][l]{\color[rgb]{0,0,0}${\scriptscriptstyle R_1}$}%
\psfrag{s16}[l][l]{\color[rgb]{0,0,0}${\scriptscriptstyle R_2}$}%
%
\psfrag{x01}[t][t]{-10}%
\psfrag{x02}[t][t]{-5}%
\psfrag{x03}[t][t]{0}%
\psfrag{x04}[t][t]{5}%
\psfrag{x05}[t][t]{10}%
%
\psfrag{v01}[r][r]{-5}%
\psfrag{v02}[r][r]{-4}%
\psfrag{v03}[r][r]{-3}%
\psfrag{v04}[r][r]{-2}%
\psfrag{v05}[r][r]{-1}%
\psfrag{v06}[r][r]{0}%
\psfrag{v07}[r][r]{1}%
\psfrag{v08}[r][r]{2}%
\psfrag{v09}[r][r]{3}%
\psfrag{v10}[r][r]{4}%
\psfrag{v11}[r][r]{5}%
%
\includegraphics[scale=\figurescale]{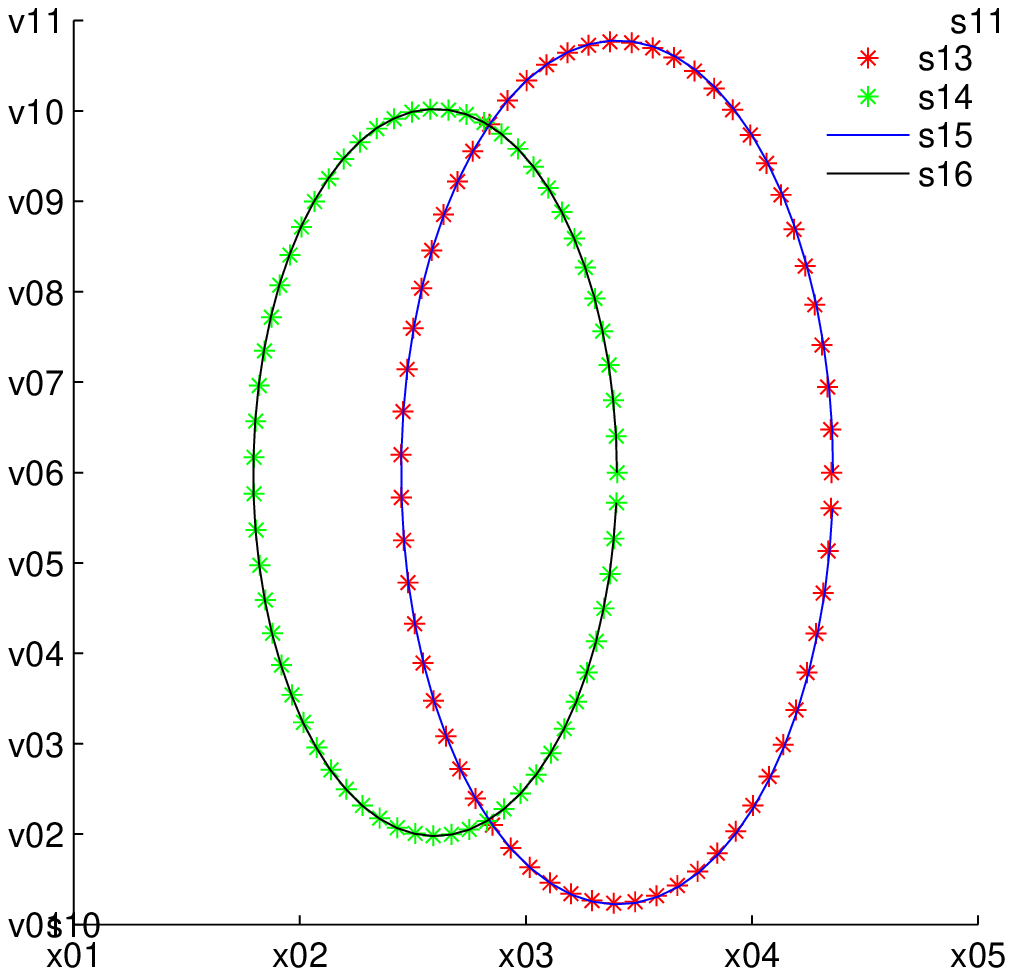}%
\end{psfrags}%
      \hspace{1cm}
\begin{psfrags}%
\psfragscanon%
%
\psfrag{s10}[][]{\color[rgb]{0,0,0}\setlength{\tabcolsep}{0pt}\begin{tabular}{c} \end{tabular}}%
\psfrag{s11}[][]{\color[rgb]{0,0,0}\setlength{\tabcolsep}{0pt}\begin{tabular}{c} \end{tabular}}%
\psfrag{s12}[l][l]{\color[rgb]{0,0,0}R2}%
\psfrag{s13}[l][l]{\color[rgb]{0,0,0}${\scriptscriptstyle r_1}$}%
\psfrag{s14}[l][l]{\color[rgb]{0,0,0}${\scriptscriptstyle r_2}$}%
\psfrag{s15}[l][l]{\color[rgb]{0,0,0}${\scriptscriptstyle R_1}$}%
\psfrag{s16}[l][l]{\color[rgb]{0,0,0}${\scriptscriptstyle R_2}$}%
%
\psfrag{x01}[t][t]{-10}%
\psfrag{x02}[t][t]{0}%
\psfrag{x03}[t][t]{10}%
\psfrag{x04}[t][t]{20}%
\psfrag{x05}[t][t]{30}%
%
\psfrag{v01}[r][r]{-20}%
\psfrag{v02}[r][r]{-15}%
\psfrag{v03}[r][r]{-10}%
\psfrag{v04}[r][r]{-5}%
\psfrag{v05}[r][r]{0}%
\psfrag{v06}[r][r]{5}%
\psfrag{v07}[r][r]{10}%
\psfrag{v08}[r][r]{15}%
\psfrag{v09}[r][r]{20}%
%
\includegraphics[scale=\figurescale]{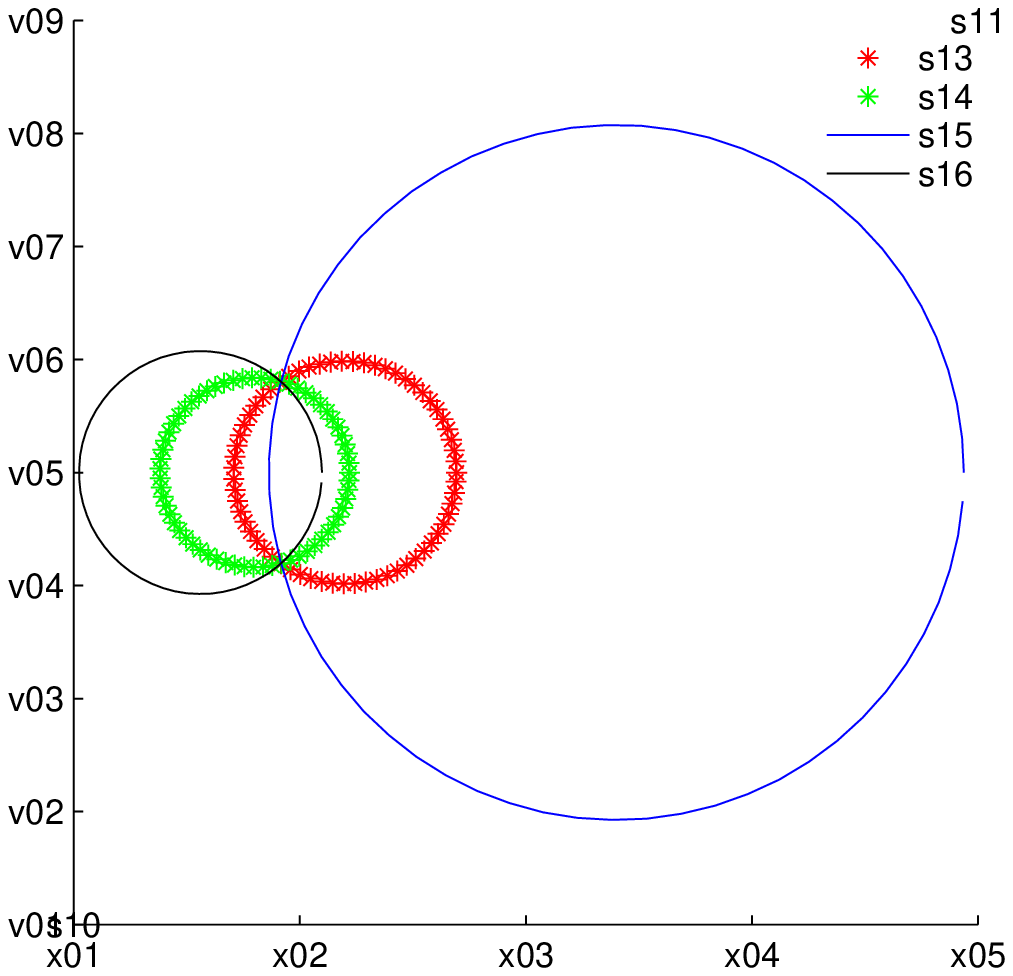}%
\end{psfrags}%
\caption{Supports of $\phi_1$ and $\phi_2$ (starred disks) and disks related to the Thomas--Fermi approximation for $\phi_1$ and
$\phi_2$ (which holds within the intersection of the unstarred disks).  We have taken $N_1=N_2=m_1=m_2=1$, $\theta_{11}=400$, $\theta_{22}=200$, 
$\theta_{12}=1$ (left figure) and $\theta_{12}=150$ (right figure). In both figures we have $x_{11}=2$, $x_{21}=-2$ 
and $x_{ij}=0$ for other $i,j$.}
\label{TFlargeintfig}
\end{center}
\end{figure}

We have {\em full overlap} between $\phi_1$ and $\phi_2$, if 
$$
(\Delta_1x)^2+(\Delta_2x)^2<|\sqrt{2\mu_1}-\sqrt{2\mu_2}|^2,
$$
so that the disks of boundaries $\Sigma_1^r$ and $\Sigma_2^r$ are
included one in the other, see Figure~\ref{supporti-position3}. 
In a highly interacting regime, this configuration leads to
a {\em non-contractible} support for one of the two wave functions, which
looses the symmetry properties of the trap. If the
potentials are both centered at the origin, $\theta_{11}\gg \theta_{22}$, and the coupling $\theta_{12}$
is sufficiently large, as
$\phi_2$ spikes around the origin, $\phi_1$ feels the influence
of the coupling $\theta_{12}|\phi_j|^2\phi_i$, lowing down the profile (around the origin) and
giving rise to a local minimum. This behaviour will be rigorously justified in the forthcoming section
via energy estimates. In the Thomas--Fermi regime the location
of the overlap regions depends upon the values of $\theta_{ij}$ and 
of the centers ${\bf x_i}$ according to formula \eqref{TFrappppree}. 
In order to visualize some situations arising with respect to the position of the trap centers, 
see Figures \ref{fig7} and \ref{fig8} where we kept 
the $\theta_{ij}$s fixed and varied the position of ${\bf x_1}$ ($x_{11}=\pm 2$ in Figure \ref{fig7}
and $x_{12}=\pm 2$ in  Figure \ref{fig8}), while ${\bf x_2}={\bf 0}$.

\medskip
\section{Strong interaction and phase separation}
\label{SIr}

In the next sections we deal with the justification of the phase separation phenomena 
occurring when the repulsive interaction between the condensates gets very strong. We consider 
both ground and excited state solutions. See also \cite{CTVpc,crooks1,crooks2,nambamim,mikawa,dancerzang,squassinaAA2008}
for various studies of spatial segregation phenomena in systems with large interactions.

\subsection{Ground state solutions}
Assume that the intra-species parameters $\theta_{ii}$s are chosen within a bounded range of values and, on the contrary, that
the inter-species interaction rate $\theta_{12}$ becomes very large, say $\theta_{12}=\kappa$, 
where we let the parameter $\kappa\geq 0$ go to infinity. For notational simplicity, we set $\eps=1$.
Let ${\mathcal H}\subset H^1(\R^2)\times H^1(\R^2)$ be the realization of the Hilbert subspace
given in the introduction and consider the energy functional $E_\kappa:{\mathcal H}\to\R$ 
re-written as
\begin{equation}
\label{totenergf-Bis}
E_\kappa(\phi_1,\phi_2)=E_\infty(\phi_1,\phi_2)+\kappa\irn |\phi_1|^2|\phi_2|^2, 
\end{equation}
where
$$
E_\infty(\phi_1,\phi_2)=\sum_{i=1}^2 E_i(\phi_i).
$$
Recalling \eqref{sfera}, the energy level of the ground state solutions is 
$$
c_\kappa=\inf_{(\phi_1,\phi_2)\in{\mathcal S}}E_\kappa(\phi_1,\phi_2).
$$
We also define the value for the limiting segregated least energy value $c_\infty$,
$$
c_\infty=\inf_{(\phi_1,\phi_2)\in{\mathcal S}_\infty}E_\infty(\phi_1,\phi_2),
$$
where we have set 
$$
{\mathcal S}_\infty=\big\{(\phi_1,\phi_2)\in {\mathcal S}:\,\,\phi_1\phi_2=0\,\,\text{a.e.\ in $\R^2$}\big\}.
$$

\begin{figure}[h!!!]
\begin{center}
\begin{psfrags}%
\psfragscanon%
%
\psfrag{s01}[b][b]{\color[rgb]{0,0,0}\setlength{\tabcolsep}{0pt}\begin{tabular}{c}$\phi_1$\end{tabular}}%
\psfrag{s02}[t][t]{\color[rgb]{0,0,0}\setlength{\tabcolsep}{0pt}\begin{tabular}{c}$x$\end{tabular}}%
\psfrag{s03}[b][b]{\color[rgb]{0,0,0}\setlength{\tabcolsep}{0pt}\begin{tabular}{c}$y$\end{tabular}}%
%
\psfrag{x01}[t][t]{-10}%
\psfrag{x02}[t][t]{-5}%
\psfrag{x03}[t][t]{0}%
\psfrag{x04}[t][t]{5}%
\psfrag{x05}[t][t]{10}%
%
\psfrag{v01}[r][r]{-10}%
\psfrag{v02}[r][r]{-8}%
\psfrag{v03}[r][r]{-6}%
\psfrag{v04}[r][r]{-4}%
\psfrag{v05}[r][r]{-2}%
\psfrag{v06}[r][r]{0}%
\psfrag{v07}[r][r]{2}%
\psfrag{v08}[r][r]{4}%
\psfrag{v09}[r][r]{6}%
\psfrag{v10}[r][r]{8}%
\psfrag{v11}[r][r]{10}%
%
\psfrag{z01}[r][r]{0}%
\psfrag{z02}[r][r]{0.05}%
\psfrag{z03}[r][r]{0.1}%
\psfrag{z04}[r][r]{0.15}%
\psfrag{z05}[r][r]{0.2}%
\psfrag{z06}[r][r]{0.25}%
\psfrag{z07}[r][r]{0.3}%
\psfrag{z08}[r][r]{0.35}%
%
\includegraphics[scale=\figurescale]{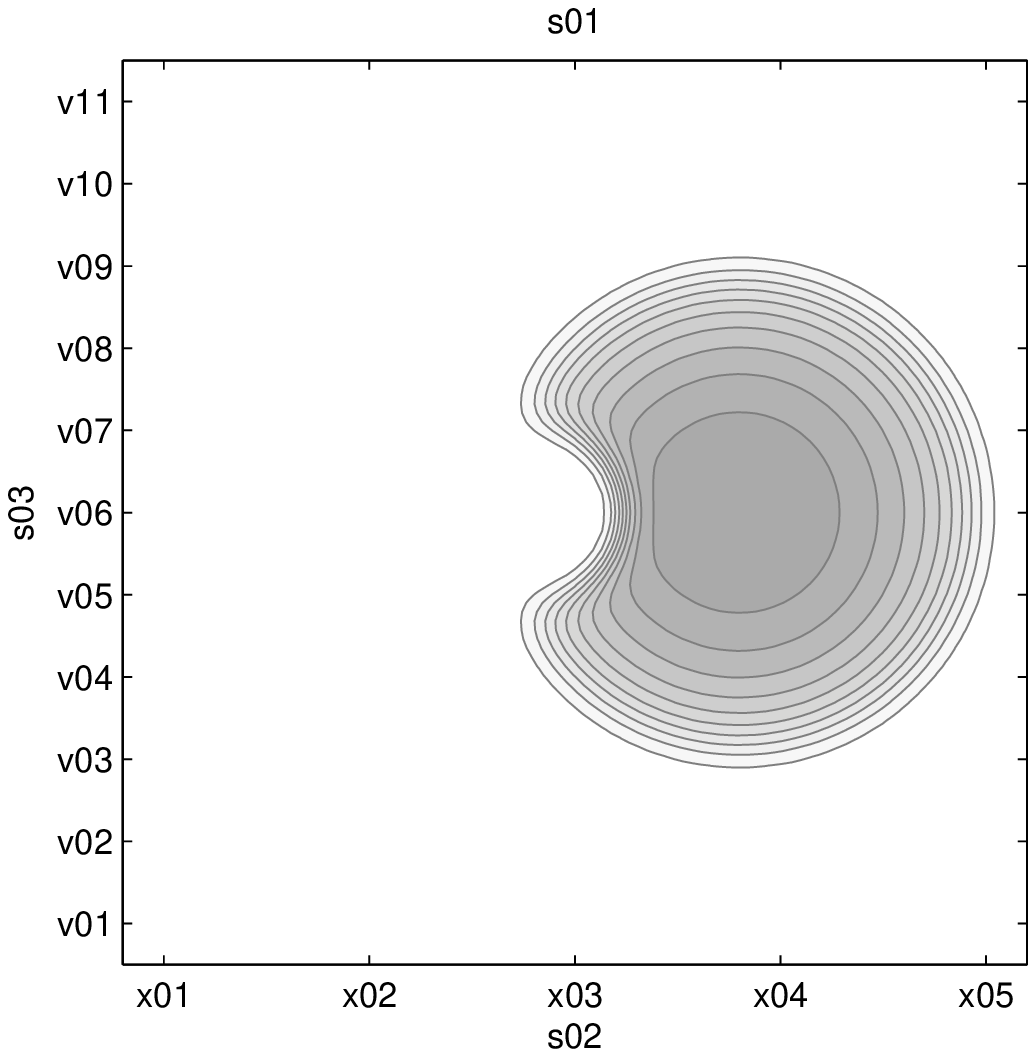}%
\end{psfrags}%
      \hspace{1cm}
\begin{psfrags}%
\psfragscanon%
%
\psfrag{s01}[b][b]{\color[rgb]{0,0,0}\setlength{\tabcolsep}{0pt}\begin{tabular}{c}$\phi_2$\end{tabular}}%
\psfrag{s02}[t][t]{\color[rgb]{0,0,0}\setlength{\tabcolsep}{0pt}\begin{tabular}{c}$x$\end{tabular}}%
\psfrag{s03}[b][b]{\color[rgb]{0,0,0}\setlength{\tabcolsep}{0pt}\begin{tabular}{c}$y$\end{tabular}}%
%
\psfrag{x01}[t][t]{-10}%
\psfrag{x02}[t][t]{-5}%
\psfrag{x03}[t][t]{0}%
\psfrag{x04}[t][t]{5}%
\psfrag{x05}[t][t]{10}%
%
\psfrag{v01}[r][r]{-10}%
\psfrag{v02}[r][r]{-8}%
\psfrag{v03}[r][r]{-6}%
\psfrag{v04}[r][r]{-4}%
\psfrag{v05}[r][r]{-2}%
\psfrag{v06}[r][r]{0}%
\psfrag{v07}[r][r]{2}%
\psfrag{v08}[r][r]{4}%
\psfrag{v09}[r][r]{6}%
\psfrag{v10}[r][r]{8}%
\psfrag{v11}[r][r]{10}%
%
\psfrag{z01}[r][r]{0}%
\psfrag{z02}[r][r]{0.05}%
\psfrag{z03}[r][r]{0.1}%
\psfrag{z04}[r][r]{0.15}%
\psfrag{z05}[r][r]{0.2}%
\psfrag{z06}[r][r]{0.25}%
\psfrag{z07}[r][r]{0.3}%
\psfrag{z08}[r][r]{0.35}%
%
\includegraphics[scale=\figurescale]{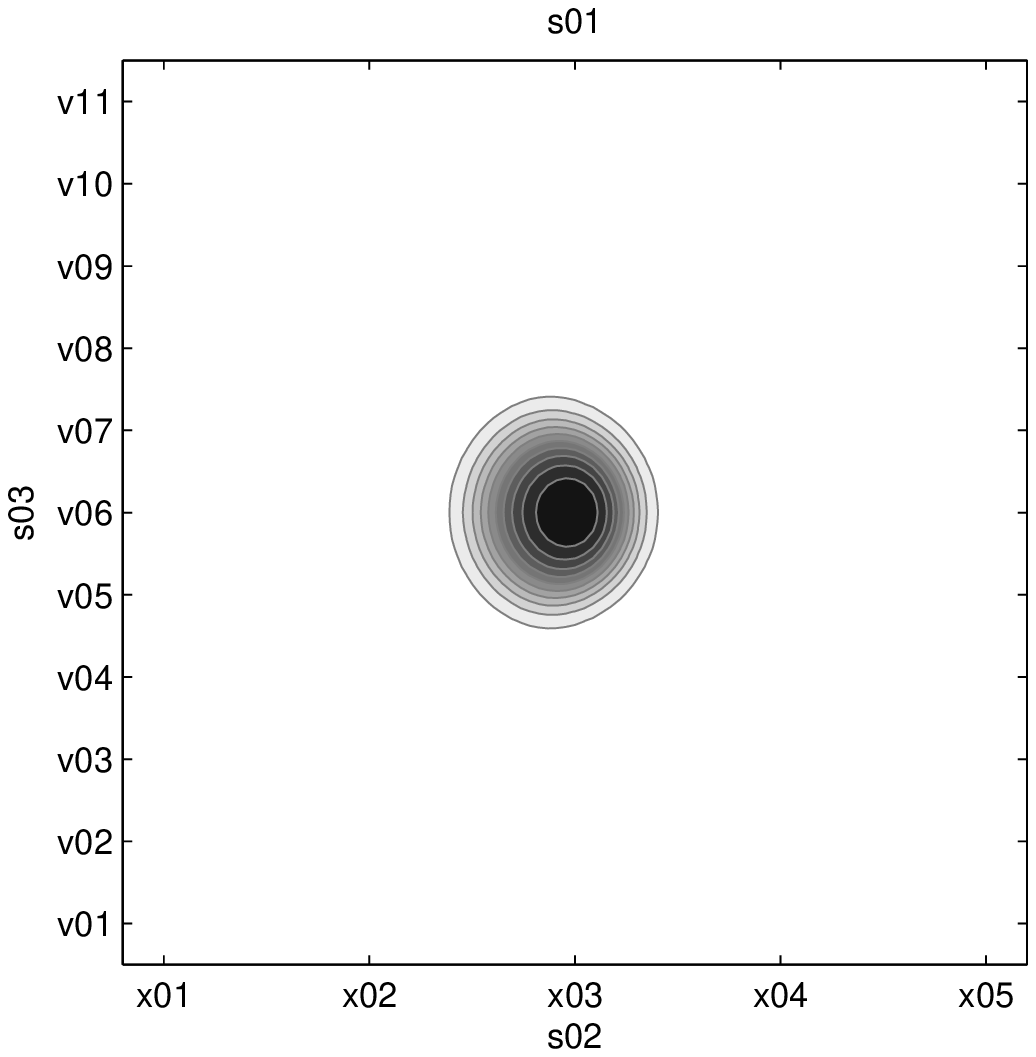}%
\end{psfrags}%
\caption{2D contour plots, in the square $[-11,11]^2$, of the ground state solution for 
$N_1=N_2=m_1=m_2=1$, $\theta_{11}=850$, $\theta_{22}=18$, $\theta_{12}=210$,
where the potentials have centers $x_{11}=4$ and $x_{ij}=0$ for any other $i,j$. The phase separation
is evident around the origin (partial overlap case).}
\label{GS2partialOverlp-bis}
\end{center}
\end{figure}

\begin{figure}[h!!!]
\begin{center}
      \begin{psfrags}%
\psfragscanon%
%
\psfrag{s01}[b][b]{\color[rgb]{0,0,0}\setlength{\tabcolsep}{0pt}\begin{tabular}{c}$\phi_1$\end{tabular}}%
\psfrag{s02}[t][t]{\color[rgb]{0,0,0}\setlength{\tabcolsep}{0pt}\begin{tabular}{c}$x$\end{tabular}}%
\psfrag{s03}[b][b]{\color[rgb]{0,0,0}\setlength{\tabcolsep}{0pt}\begin{tabular}{c}$y$\end{tabular}}%
%
\psfrag{x01}[t][t]{-10}%
\psfrag{x02}[t][t]{-5}%
\psfrag{x03}[t][t]{0}%
\psfrag{x04}[t][t]{5}%
\psfrag{x05}[t][t]{10}%
%
\psfrag{v01}[r][r]{-10}%
\psfrag{v02}[r][r]{-8}%
\psfrag{v03}[r][r]{-6}%
\psfrag{v04}[r][r]{-4}%
\psfrag{v05}[r][r]{-2}%
\psfrag{v06}[r][r]{0}%
\psfrag{v07}[r][r]{2}%
\psfrag{v08}[r][r]{4}%
\psfrag{v09}[r][r]{6}%
\psfrag{v10}[r][r]{8}%
\psfrag{v11}[r][r]{10}%
%
\psfrag{z01}[r][r]{0}%
\psfrag{z02}[r][r]{0.05}%
\psfrag{z03}[r][r]{0.1}%
\psfrag{z04}[r][r]{0.15}%
\psfrag{z05}[r][r]{0.2}%
\psfrag{z06}[r][r]{0.25}%
\psfrag{z07}[r][r]{0.3}%
\psfrag{z08}[r][r]{0.35}%
\psfrag{z09}[r][r]{0.4}%
%
\includegraphics[scale=\figurescale]{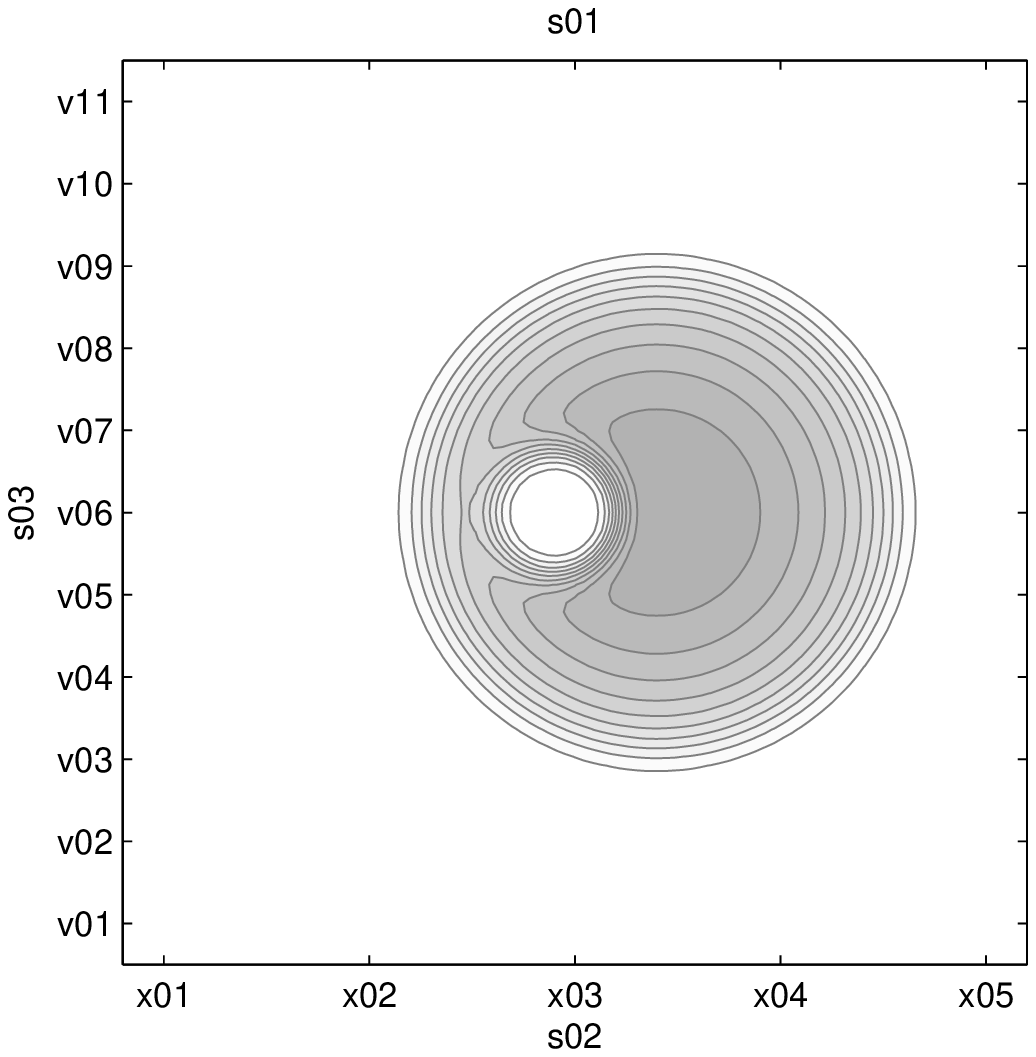}%
\end{psfrags}%
      \hspace{1cm}
     \begin{psfrags}%
\psfragscanon%
%
\psfrag{s01}[b][b]{\color[rgb]{0,0,0}\setlength{\tabcolsep}{0pt}\begin{tabular}{c}$\phi_1$\end{tabular}}%
\psfrag{s02}[t][t]{\color[rgb]{0,0,0}\setlength{\tabcolsep}{0pt}\begin{tabular}{c}$x$\end{tabular}}%
\psfrag{s03}[b][b]{\color[rgb]{0,0,0}\setlength{\tabcolsep}{0pt}\begin{tabular}{c}$y$\end{tabular}}%
%
\psfrag{x01}[t][t]{-10}%
\psfrag{x02}[t][t]{-5}%
\psfrag{x03}[t][t]{0}%
\psfrag{x04}[t][t]{5}%
\psfrag{x05}[t][t]{10}%
%
\psfrag{v01}[r][r]{-10}%
\psfrag{v02}[r][r]{-8}%
\psfrag{v03}[r][r]{-6}%
\psfrag{v04}[r][r]{-4}%
\psfrag{v05}[r][r]{-2}%
\psfrag{v06}[r][r]{0}%
\psfrag{v07}[r][r]{2}%
\psfrag{v08}[r][r]{4}%
\psfrag{v09}[r][r]{6}%
\psfrag{v10}[r][r]{8}%
\psfrag{v11}[r][r]{10}%
%
\psfrag{z01}[r][r]{0}%
\psfrag{z02}[r][r]{0.05}%
\psfrag{z03}[r][r]{0.1}%
\psfrag{z04}[r][r]{0.15}%
\psfrag{z05}[r][r]{0.2}%
\psfrag{z06}[r][r]{0.25}%
\psfrag{z07}[r][r]{0.3}%
\psfrag{z08}[r][r]{0.35}%
\psfrag{z09}[r][r]{0.4}%
%
\includegraphics[scale=\figurescale]{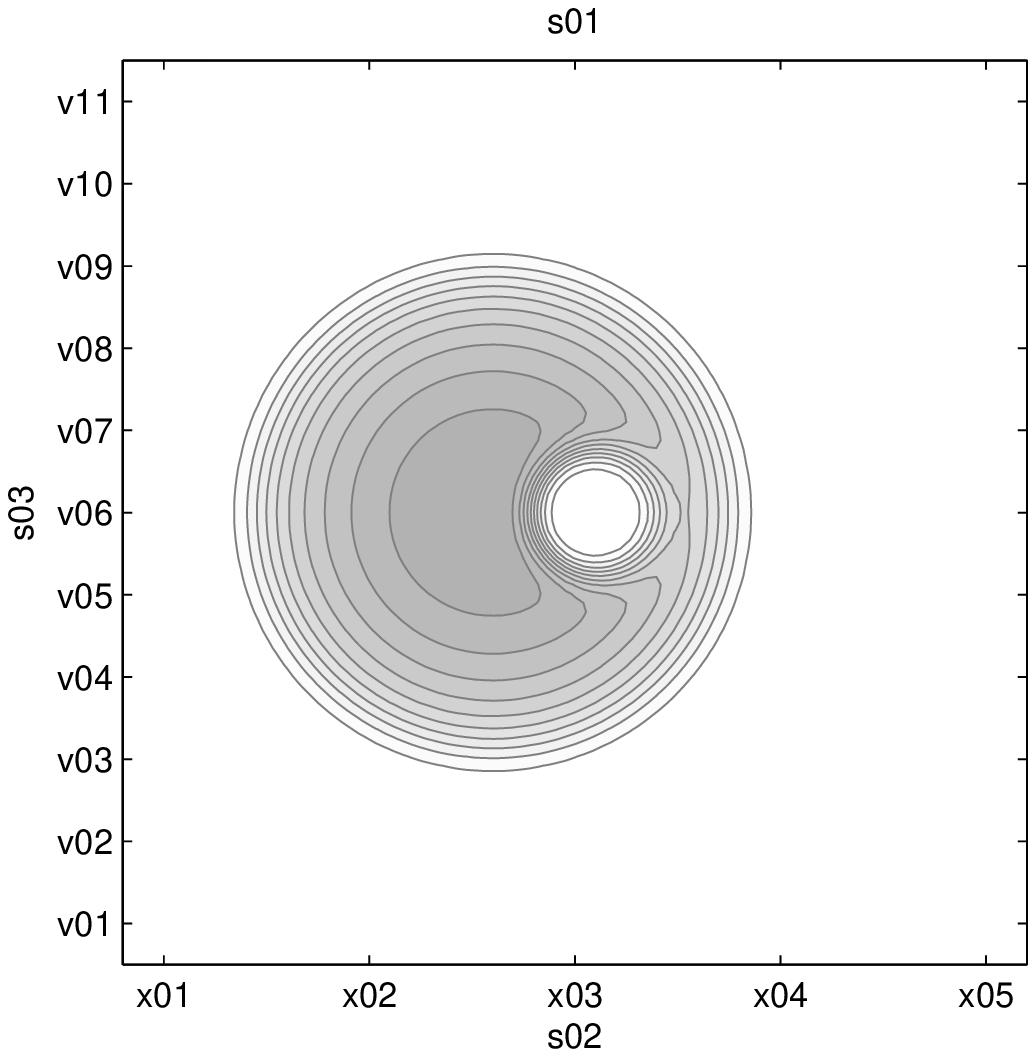}%
\end{psfrags}%
\caption{2D contour plot, in the square $[-11,11]^2$, of the first component of ground state solution for 
$N_1=N_2=m_1=m_2=1$, $\theta_{11}=850$, $\theta_{22}=18$, $\theta_{12}=210$ in the cases
where the potentials have centers $x_{11}=2$ and $x_{ij}=0$ (left) and 
$x_{11}=-2$ and $x_{ij}=0$ (right). The symmetry breaking
is evident (full overlap case).}\label{fig7}
\end{center}
\end{figure}

\begin{figure}[h!!!]
\begin{center}
   \begin{psfrags}%
\psfragscanon%
%
\psfrag{s01}[b][b]{\color[rgb]{0,0,0}\setlength{\tabcolsep}{0pt}\begin{tabular}{c}$\phi_1$\end{tabular}}%
\psfrag{s02}[t][t]{\color[rgb]{0,0,0}\setlength{\tabcolsep}{0pt}\begin{tabular}{c}$x$\end{tabular}}%
\psfrag{s03}[b][b]{\color[rgb]{0,0,0}\setlength{\tabcolsep}{0pt}\begin{tabular}{c}$y$\end{tabular}}%
%
\psfrag{x01}[t][t]{-10}%
\psfrag{x02}[t][t]{-5}%
\psfrag{x03}[t][t]{0}%
\psfrag{x04}[t][t]{5}%
\psfrag{x05}[t][t]{10}%
%
\psfrag{v01}[r][r]{-10}%
\psfrag{v02}[r][r]{-8}%
\psfrag{v03}[r][r]{-6}%
\psfrag{v04}[r][r]{-4}%
\psfrag{v05}[r][r]{-2}%
\psfrag{v06}[r][r]{0}%
\psfrag{v07}[r][r]{2}%
\psfrag{v08}[r][r]{4}%
\psfrag{v09}[r][r]{6}%
\psfrag{v10}[r][r]{8}%
\psfrag{v11}[r][r]{10}%
%
\psfrag{z01}[r][r]{0}%
\psfrag{z02}[r][r]{0.05}%
\psfrag{z03}[r][r]{0.1}%
\psfrag{z04}[r][r]{0.15}%
\psfrag{z05}[r][r]{0.2}%
\psfrag{z06}[r][r]{0.25}%
\psfrag{z07}[r][r]{0.3}%
\psfrag{z08}[r][r]{0.35}%
\psfrag{z09}[r][r]{0.4}%
%
\includegraphics[scale=\figurescale]{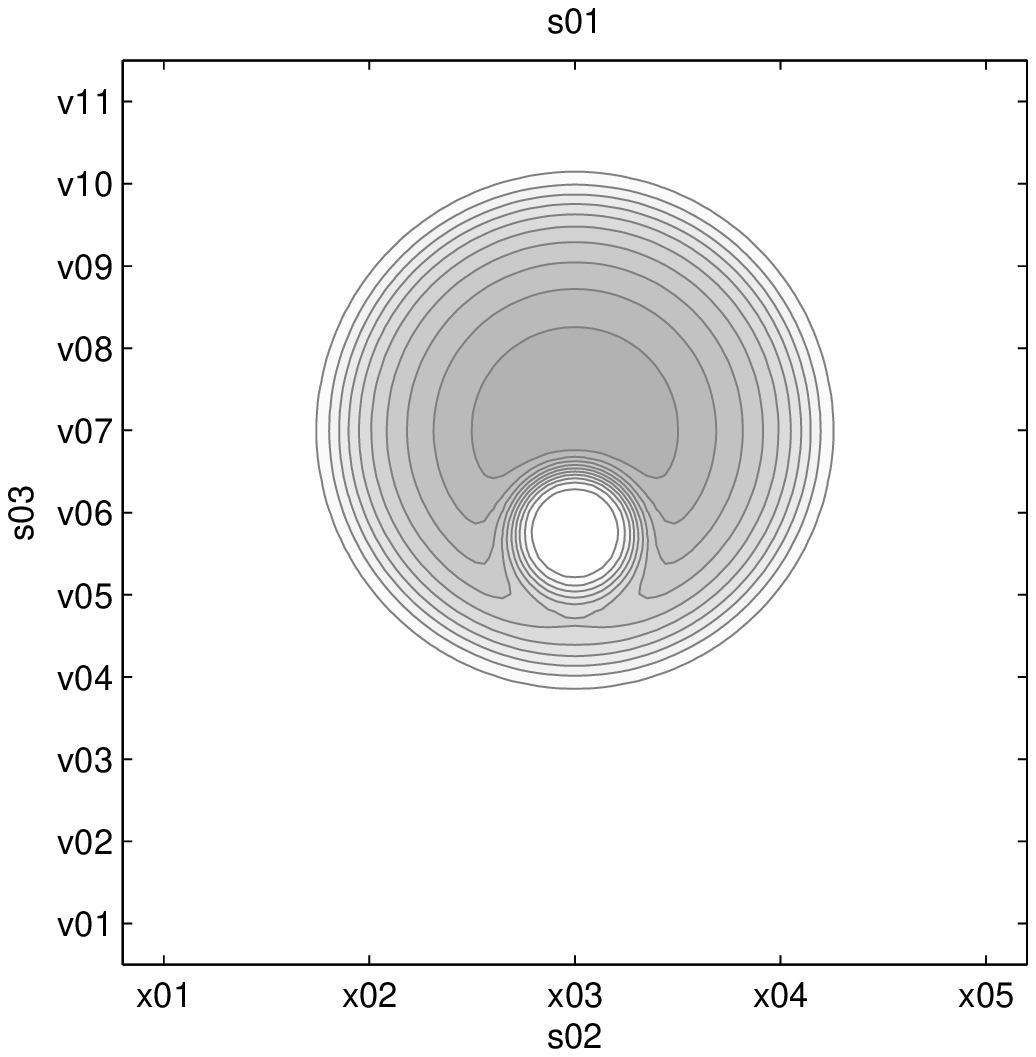}%
\end{psfrags}%
      \hspace{1cm}
\begin{psfrags}%
\psfragscanon%
%
\psfrag{s01}[b][b]{\color[rgb]{0,0,0}\setlength{\tabcolsep}{0pt}\begin{tabular}{c}$\phi_1$\end{tabular}}%
\psfrag{s02}[t][t]{\color[rgb]{0,0,0}\setlength{\tabcolsep}{0pt}\begin{tabular}{c}$x$\end{tabular}}%
\psfrag{s03}[b][b]{\color[rgb]{0,0,0}\setlength{\tabcolsep}{0pt}\begin{tabular}{c}$y$\end{tabular}}%
%
\psfrag{x01}[t][t]{-10}%
\psfrag{x02}[t][t]{-5}%
\psfrag{x03}[t][t]{0}%
\psfrag{x04}[t][t]{5}%
\psfrag{x05}[t][t]{10}%
%
\psfrag{v01}[r][r]{-10}%
\psfrag{v02}[r][r]{-8}%
\psfrag{v03}[r][r]{-6}%
\psfrag{v04}[r][r]{-4}%
\psfrag{v05}[r][r]{-2}%
\psfrag{v06}[r][r]{0}%
\psfrag{v07}[r][r]{2}%
\psfrag{v08}[r][r]{4}%
\psfrag{v09}[r][r]{6}%
\psfrag{v10}[r][r]{8}%
\psfrag{v11}[r][r]{10}%
%
\psfrag{z01}[r][r]{0}%
\psfrag{z02}[r][r]{0.05}%
\psfrag{z03}[r][r]{0.1}%
\psfrag{z04}[r][r]{0.15}%
\psfrag{z05}[r][r]{0.2}%
\psfrag{z06}[r][r]{0.25}%
\psfrag{z07}[r][r]{0.3}%
\psfrag{z08}[r][r]{0.35}%
\psfrag{z09}[r][r]{0.4}%
%
\includegraphics[scale=\figurescale]{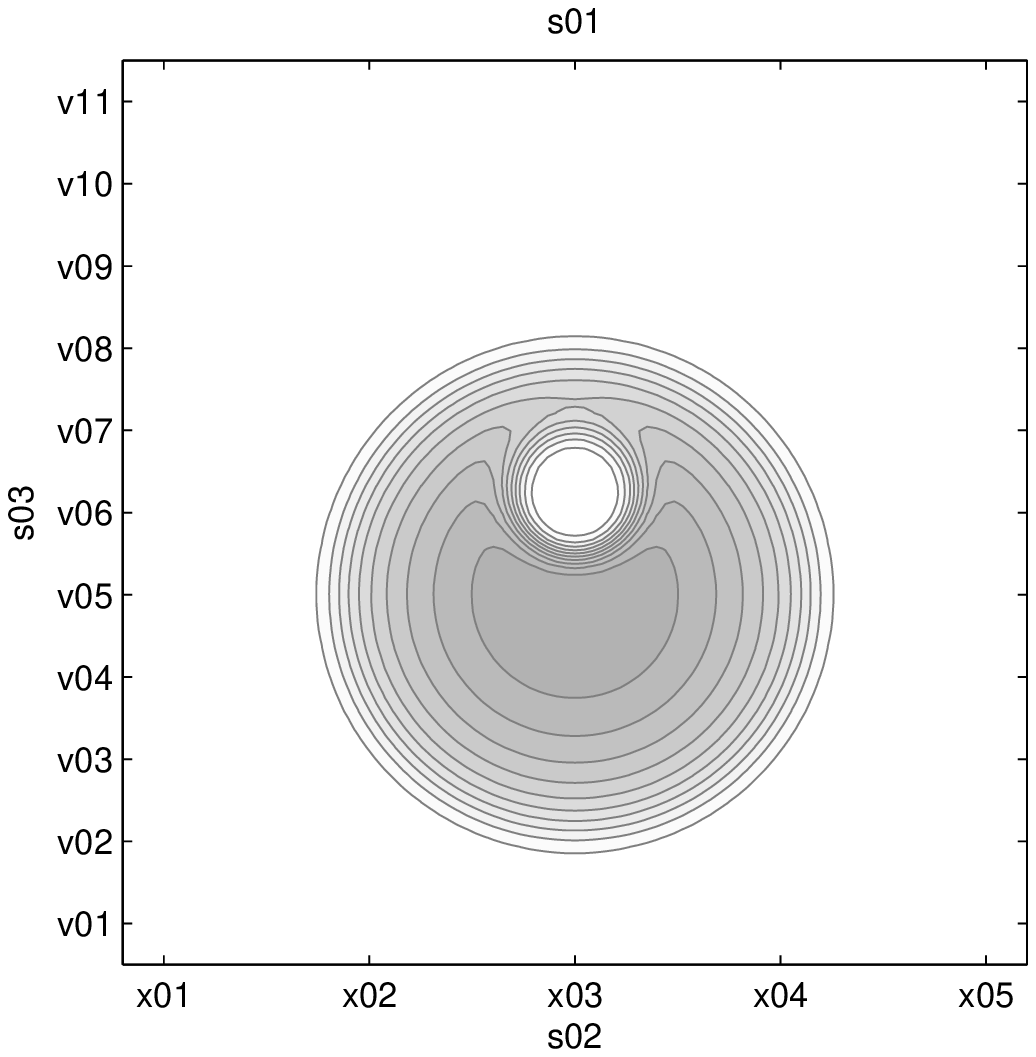}%
\end{psfrags}%
%
\caption{2D contour plot, in the square $[-11,11]^2$, of the first component of ground state solution for 
$N_1=N_2=m_1=m_2=1$, $\theta_{11}=850$, $\theta_{22}=18$, $\theta_{12}=210$ in the cases
where the potentials have centers $x_{12}=2$ and $x_{ij}=0$ (left) and 
$x_{12}=-2$ and $x_{ij}=0$ (right). The symmetry breaking
is evident (full overlap case).}\label{fig8}
\end{center}
\end{figure}

\begin{proposition}
\label{Gsspatseg}
The sequence $(\phi_1^\kappa,\phi_2^\kappa)\subset {\mathcal S}$ of ground state solutions of \eqref{systemGPGen} 
converges in ${\mathcal H}$ to a function $(\phi_1^\infty,\phi_2^\infty)\in {\mathcal S}_\infty$ 
at the energy level $c_\infty$. Furthermore, 
\begin{equation}
\label{varineqqq}
-\frac{1}{2m_i}\Delta \phi_i^\infty+V_i(x_1,x_2)\phi_i^\infty+\theta_{ii}|\phi_i^\infty|^2\phi_i^\infty
\leq \mu_i^\infty\phi_i^\infty,
\end{equation}
where 
$$
N_i\mu_i^\infty=E_i(\phi_i^\infty)+\frac{\theta_{ii}}{2}\irn |\phi_i^\infty|^4,
$$
for $i=1,2$.
\end{proposition}
\begin{proof}
The infimum that defines the value $c_\infty$ is taken over a smaller set with
respect to the one defining $c_\kappa$. Moreover, if the functions $\phi_1$ and $\phi_2$ have disjoint supports, 
$E_\kappa(\phi_1,\phi_2)=E_\infty(\phi_1,\phi_2)$, for any $\kappa>0$. In particular, of course, this 
implies that $c_\kappa\leq c_\infty$, for all $\kappa>0$. Therefore, for the ground state 
solutions $(\phi_1^\kappa,\phi_2^\kappa)\in {\mathcal H}$, $\phi_i^\kappa\not\equiv 0$ for $i=1,2$,
we have $E_\kappa(\phi_1^\kappa,\phi_2^\kappa)=c_\kappa$ and
\begin{equation}
\label{BddkappaInt}
\kappa\irn |\phi_1^\kappa|^2|\phi_2^\kappa|^2\leq E_\infty(\phi_1^\kappa,\phi_2^\kappa)
+\kappa\irn |\phi_1^\kappa|^2|\phi_2^\kappa|^2=c_\kappa\leq c_\infty,
\end{equation}
for every $\kappa>0$. As a consequence,
\begin{equation}
\label{weaksegreg}
\lim_{\kappa\to\infty}\irn |\phi_1^\kappa|^2|\phi_2^\kappa|^2=0.
\end{equation}
Also, for all $\kappa>0$, we have
\begin{equation*}
\|(\phi_1^\kappa,\phi_2^\kappa)\|_{\mathcal H}^2\leq  E_\infty(\phi_1^\kappa,\phi_2^\kappa) 
+\kappa \irn |\phi_1^\kappa|^2|\phi_2^\kappa|^2\leq c_\infty,
\end{equation*}
so that the sequences $(\phi_1^\kappa,\phi_2^\kappa)$ is uniformly bounded in ${\mathcal H}$.
In particular, up to a subsequence, there exist $(\phi_1^\infty,\phi_2^\infty)$ in ${\mathcal H}$ 
such that $(\phi_1^\kappa,\phi_2^\kappa)\rightharpoonup (\phi_1^\infty,\phi_2^\infty)$ in ${\mathcal H}$ 
as $\kappa\to\infty$ and $\phi_i^\kappa(x_1,x_2)\to\phi_i^\infty(x_1,x_2)$ a.e. in $\R^2$. Hence,
by combining Fatou's Lemma with formula \eqref{weaksegreg}, we get 
$$
\irn |\phi_1^\infty|^2|\phi_2^\infty|^2=0,
$$ 
so that
\begin{equation}
\label{segrlimite}
\phi_1^\infty\phi_2^\infty=0,\quad\text{a.e.\ in $\R^2$}.
\end{equation}
Since by definition of ground state solution $\irn|\phi_i^\kappa|^2=N_i$, for any $\kappa>0$, and ${\mathcal H}$
in compactly embedded into $L^2(\R^2)\times L^2(\R^2)$  (see inequality \eqref{embedkeyineq}), 
up to passing to a further subsequence, we conclude that
\begin{equation}
\label{normlimite}
\irn|\phi_i^\infty|^2=N_i,
\end{equation}
for $i=1,2$.
In particular $(\phi_1^\infty,\phi_2^\infty)\in {\mathcal S}_\infty$, by virtue of
\eqref{segrlimite} and \eqref{normlimite}. Observe also that, by virtue of \eqref{eig-systemGPGen},
\eqref{BddkappaInt} and $E_i(\phi_i^\kappa)\leq c_\infty$,
\begin{equation*}
\sup_{\kappa\geq 1}\mu_i^\kappa=\frac{1}{N_i}\sup_{\kappa\geq 1}\left\{E_i(\phi_i^\kappa)+\frac{\theta_{ii}}{2}
\irn |\phi_i^\kappa|^4+\kappa\irn |\phi_1^\kappa|^2|\phi_2^\kappa|^2\right\}<\infty,
\end{equation*}
being $\mu_i^\kappa$ the eigenvalues corresponding to  $\phi_i^\kappa$.
Then, up to a subsequence, $\mu_i^\kappa\to \mu_i^\infty$ as $\kappa\to\infty$. By testing the equations of 
\eqref{systemGPGen} by an arbitrary positive function $\eta$ with compact support, we get
\begin{equation*}
\frac{1}{2m_i}\irn\nabla \phi_i^\kappa\cdot\nabla\eta+\irn V_i(x_1,x_2)\phi_i^\kappa\eta+\theta_{ii}\irn|\phi_i^\kappa|^2
\phi_i^\kappa\eta \leq \mu_i^\kappa\irn \phi_i^\kappa\eta, 
\end{equation*}
for all $\kappa>0$ and any $\eta\in C^\infty_c(\R^2)$ with $\eta\geq 0$.
Hence, letting $\kappa\to\infty$, it turns out that $\phi_i^\infty$
satisfies the variational inequality \eqref{varineqqq}.
Notice that, since $(\phi_1^\infty,\phi_2^\infty)\in {\mathcal S}_\infty$, by the definition of $c_\infty$, we deduce
\begin{align*}
 & \sum_{i=1}^2\frac{1}{2m_i}\irn |\nabla \phi_i^\infty|^2+ \sum_{i=1}^2\irn V_i|\phi_i^\infty|^2+ \sum_{i=1}^2\frac{\theta_{ii}}{2}  
 \irn |\phi_i^\infty|^4  +\lim_{\kappa\to\infty}\kappa\irn |\phi_1^\kappa|^2|\phi_2^\kappa|^2 \\
 \noalign{\vskip2pt}
& \leq  \sum_{i=1}^2\frac{1}{2m_i}\liminf_{\kappa\to\infty}\irn |\nabla \phi_i^\kappa|^2+  
\sum_{i=1}^2\liminf_{\kappa\to\infty}\irn V_i|\phi_i^\kappa|^2+ 
 \sum_{i=1}^2 \frac{\theta_{ii}}{2}\liminf_{\kappa\to\infty}\irn |\phi_i^\kappa|^4  \\
 \noalign{\vskip2pt}
 & +\lim_{\kappa\to\infty}\kappa\irn |\phi_1^\kappa|^2|\phi_2^\kappa|^2 
  \leq \liminf_{\kappa\to\infty} E_\kappa(\phi_1^\kappa,\phi_2^\kappa)=\liminf_{\kappa\to\infty} c_\kappa
  \leq c_\infty\leq E_\infty(\phi_1^\infty,\phi_2^\infty) \\
 \noalign{\vskip2pt}
 & =  \sum_{i=1}^2\frac{1}{2m_i}\irn |\nabla \phi_i^\infty|^2+ 
 \sum_{i=1}^2\irn V_i|\phi_i^\infty|^2+ \sum_{i=1}^2\frac{\theta_{ii}}{2}\irn |\phi_i^\infty|^4,
\end{align*}
which yields $\kappa\irn |\phi_1^\kappa|^2|\phi_2^\kappa|^2\to 0$ as $\kappa\to\infty$,
which is a much stronger conclusion compared with \eqref{weaksegreg}. Consequently, 
the convergence of $\phi_i^\kappa$ to $\phi_i^\infty$ in ${\mathcal H}$
is strong, otherwise, assuming by contradiction that for some $i=1,2$
$$
\irn |\nabla \phi_i^\infty|^2<\lim_{\kappa\to\infty}\irn |\nabla \phi_i^\kappa|^2
\quad\text{or}\quad
\irn V_i |\phi_i^\infty|^2<\lim_{\kappa\to\infty}\irn V_i|\phi_i^\kappa|^2,
$$
the previous inequalities we would become strict, yielding immediately
a contradiction. Finally, as a further consequence, $c_\kappa\to c_\infty$ as $\kappa\to\infty$ and
the value $c_\infty$ is indeed assumed and
$$
c_\infty=\sum_{i=1}^2\frac{1}{2m_i}\irn |\nabla \phi_i^\infty|^2+ 
\sum_{i=1}^2\irn V_i |\phi_i^\infty|^2+ \sum_{i=1}^2\frac{\theta_{ii}}{2}\irn |\phi_i^\infty|^4=E_\infty(\phi_1^\infty,\phi_2^\infty).
$$ 
Finally, the strong convergence and \eqref{eig-systemGPGen} yield 
$$
N_i\mu_i^\infty=E_i(\phi_i^\infty)+\frac{\theta_{ii}}{2}\irn |\phi_i^\infty|^4
$$ 
for any $i=1,2$, which concludes the proof.
\end{proof}
As one can see in the numerical simulations, as the interaction
coefficient gets large, the phase separation becomes rather evident.
See Figures \ref{GS2partialOverlp-bis} and \ref{fig7}-\ref{fig8} (just $\phi_1$ component) where
different choices of the centers of the $V_i$s have been considered.

\subsection{The anisotropic case}
\begin{figure}[h!!!]
\begin{center}
\begin{psfrags}%
\psfragscanon%
%
\psfrag{s01}[b][b]{\color[rgb]{0,0,0}\setlength{\tabcolsep}{0pt}\begin{tabular}{c}$\phi_1$\end{tabular}}%
\psfrag{s02}[t][t]{\color[rgb]{0,0,0}\setlength{\tabcolsep}{0pt}\begin{tabular}{c}$x$\end{tabular}}%
\psfrag{s03}[b][b]{\color[rgb]{0,0,0}\setlength{\tabcolsep}{0pt}\begin{tabular}{c}$y$\end{tabular}}%
%
\psfrag{x01}[t][t]{-5}%
\psfrag{x02}[t][t]{0}%
\psfrag{x03}[t][t]{5}%
%
\psfrag{v01}[r][r]{-5}%
\psfrag{v02}[r][r]{-4}%
\psfrag{v03}[r][r]{-3}%
\psfrag{v04}[r][r]{-2}%
\psfrag{v05}[r][r]{-1}%
\psfrag{v06}[r][r]{0}%
\psfrag{v07}[r][r]{1}%
\psfrag{v08}[r][r]{2}%
\psfrag{v09}[r][r]{3}%
\psfrag{v10}[r][r]{4}%
\psfrag{v11}[r][r]{5}%
%
\psfrag{z01}[r][r]{0}%
\psfrag{z02}[r][r]{0.1}%
\psfrag{z03}[r][r]{0.2}%
\psfrag{z04}[r][r]{0.3}%
\psfrag{z05}[r][r]{0.4}%
\psfrag{z06}[r][r]{0.5}%
\psfrag{z07}[r][r]{0.6}%
\psfrag{z08}[r][r]{0.7}%
%
\includegraphics[scale=\figurescale]{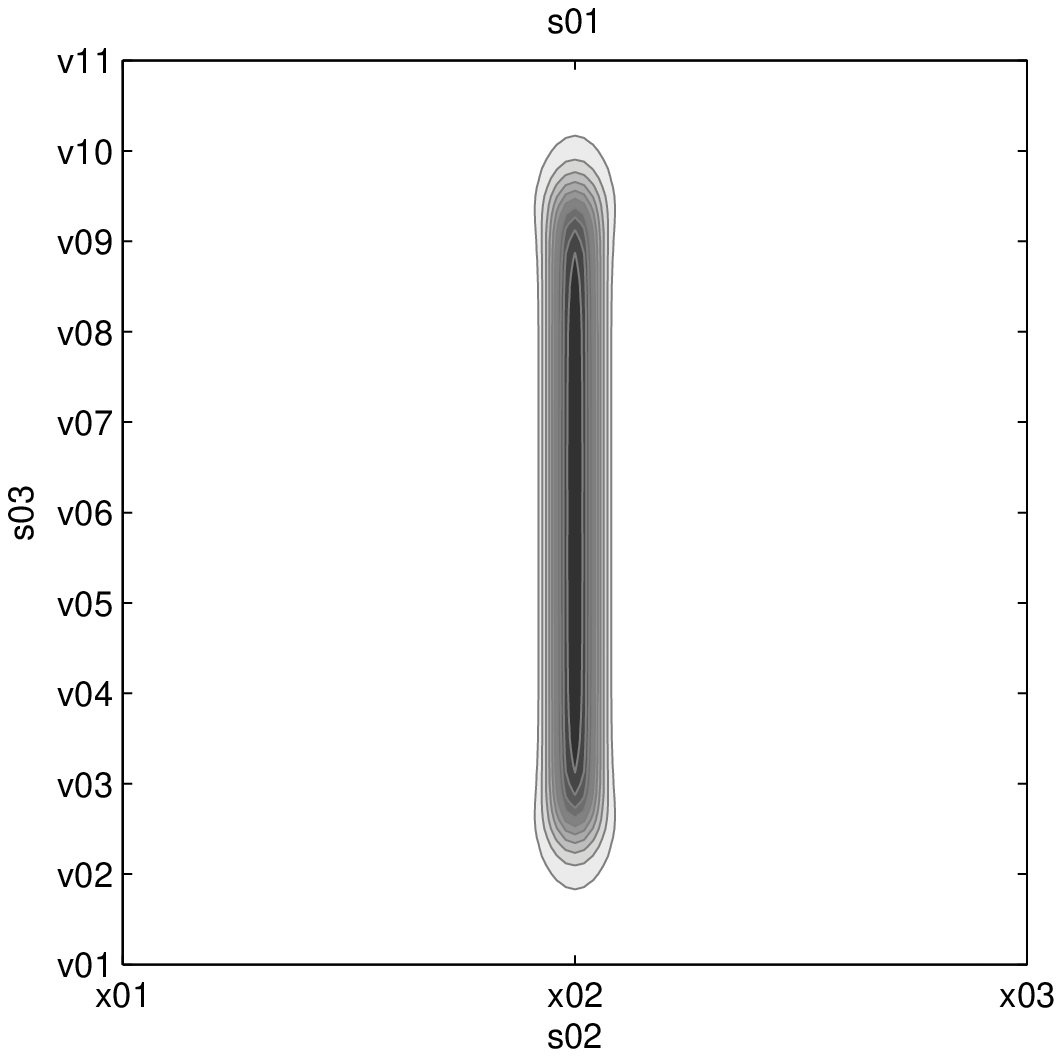}%
\end{psfrags}%
     \hspace{1cm}
\begin{psfrags}%
\psfragscanon%
%
\psfrag{s01}[b][b]{\color[rgb]{0,0,0}\setlength{\tabcolsep}{0pt}\begin{tabular}{c}$\phi_2$\end{tabular}}%
\psfrag{s02}[t][t]{\color[rgb]{0,0,0}\setlength{\tabcolsep}{0pt}\begin{tabular}{c}$x$\end{tabular}}%
\psfrag{s03}[b][b]{\color[rgb]{0,0,0}\setlength{\tabcolsep}{0pt}\begin{tabular}{c}$y$\end{tabular}}%
%
\psfrag{x01}[t][t]{-5}%
\psfrag{x02}[t][t]{0}%
\psfrag{x03}[t][t]{5}%
%
\psfrag{v01}[r][r]{-5}%
\psfrag{v02}[r][r]{-4}%
\psfrag{v03}[r][r]{-3}%
\psfrag{v04}[r][r]{-2}%
\psfrag{v05}[r][r]{-1}%
\psfrag{v06}[r][r]{0}%
\psfrag{v07}[r][r]{1}%
\psfrag{v08}[r][r]{2}%
\psfrag{v09}[r][r]{3}%
\psfrag{v10}[r][r]{4}%
\psfrag{v11}[r][r]{5}%
%
\psfrag{z01}[r][r]{0}%
\psfrag{z02}[r][r]{0.1}%
\psfrag{z03}[r][r]{0.2}%
\psfrag{z04}[r][r]{0.3}%
\psfrag{z05}[r][r]{0.4}%
\psfrag{z06}[r][r]{0.5}%
\psfrag{z07}[r][r]{0.6}%
\psfrag{z08}[r][r]{0.7}%
%
\includegraphics[scale=\figurescale]{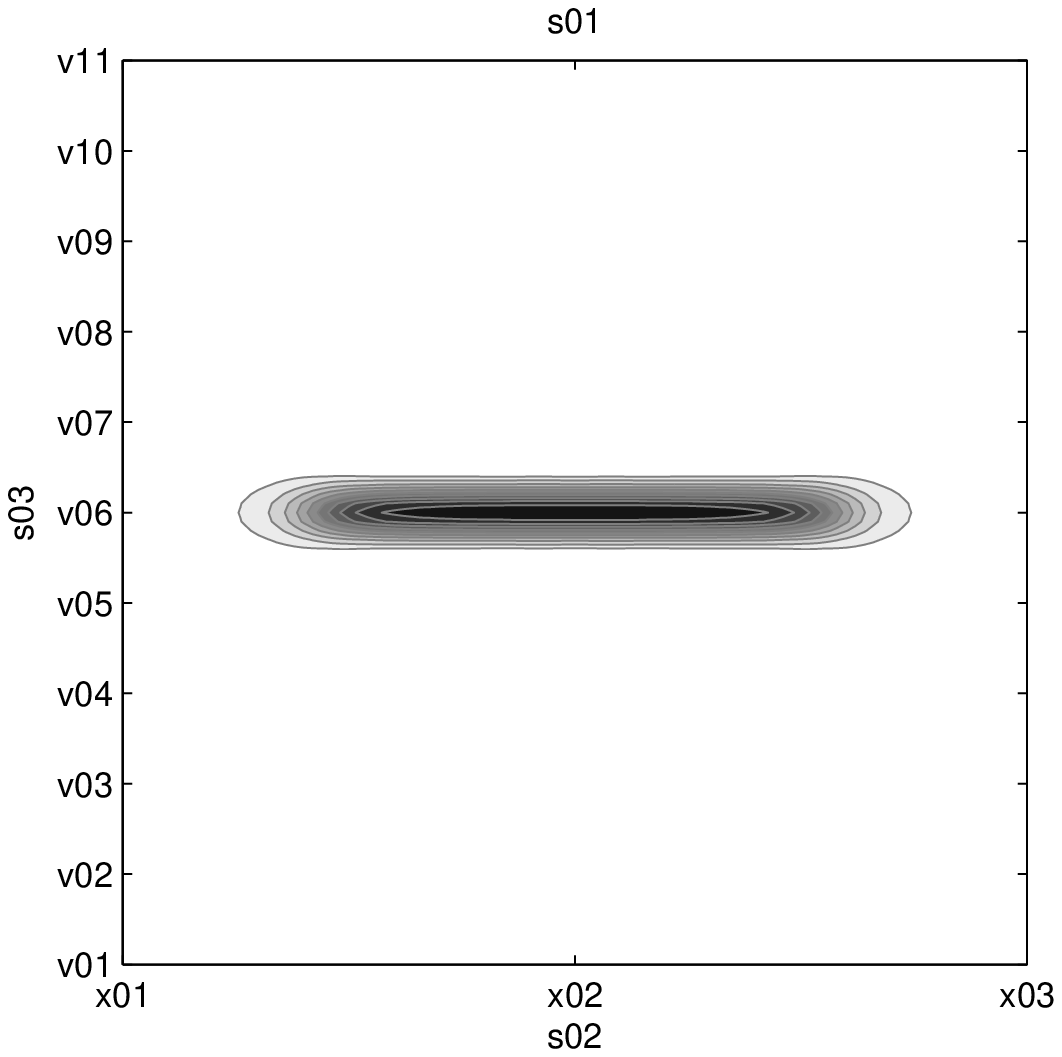}%
\end{psfrags}%
\caption{2D contour plots, in the square $[-5,5]^2$, of the ground state solution for $N_1=N_2=m_1=m_2=1$,
$\theta_{11}=400$, $\theta_{22}=150$, $\theta_{12}=1$, $\omega_{11}=100$, $\omega_{22}=100$,
$\omega_{12}=\omega_{21}=1$ and potentials centered at the origin.}\label{GS2sigaro}
\end{center}
\end{figure}

\begin{figure}[h!!!]
\begin{center}
\begin{psfrags}%
\psfragscanon%
%
\psfrag{s01}[b][b]{\color[rgb]{0,0,0}\setlength{\tabcolsep}{0pt}\begin{tabular}{c}$\phi_1$\end{tabular}}%
\psfrag{s02}[t][t]{\color[rgb]{0,0,0}\setlength{\tabcolsep}{0pt}\begin{tabular}{c}$x$\end{tabular}}%
\psfrag{s03}[b][b]{\color[rgb]{0,0,0}\setlength{\tabcolsep}{0pt}\begin{tabular}{c}$y$\end{tabular}}%
%
\psfrag{x01}[t][t]{-5}%
\psfrag{x02}[t][t]{0}%
\psfrag{x03}[t][t]{5}%
%
\psfrag{v01}[r][r]{-5}%
\psfrag{v02}[r][r]{-4}%
\psfrag{v03}[r][r]{-3}%
\psfrag{v04}[r][r]{-2}%
\psfrag{v05}[r][r]{-1}%
\psfrag{v06}[r][r]{0}%
\psfrag{v07}[r][r]{1}%
\psfrag{v08}[r][r]{2}%
\psfrag{v09}[r][r]{3}%
\psfrag{v10}[r][r]{4}%
\psfrag{v11}[r][r]{5}%
%
\psfrag{z01}[r][r]{0}%
\psfrag{z02}[r][r]{0.1}%
\psfrag{z03}[r][r]{0.2}%
\psfrag{z04}[r][r]{0.3}%
\psfrag{z05}[r][r]{0.4}%
\psfrag{z06}[r][r]{0.5}%
\psfrag{z07}[r][r]{0.6}%
\psfrag{z08}[r][r]{0.7}%
%
\includegraphics[scale=\figurescale]{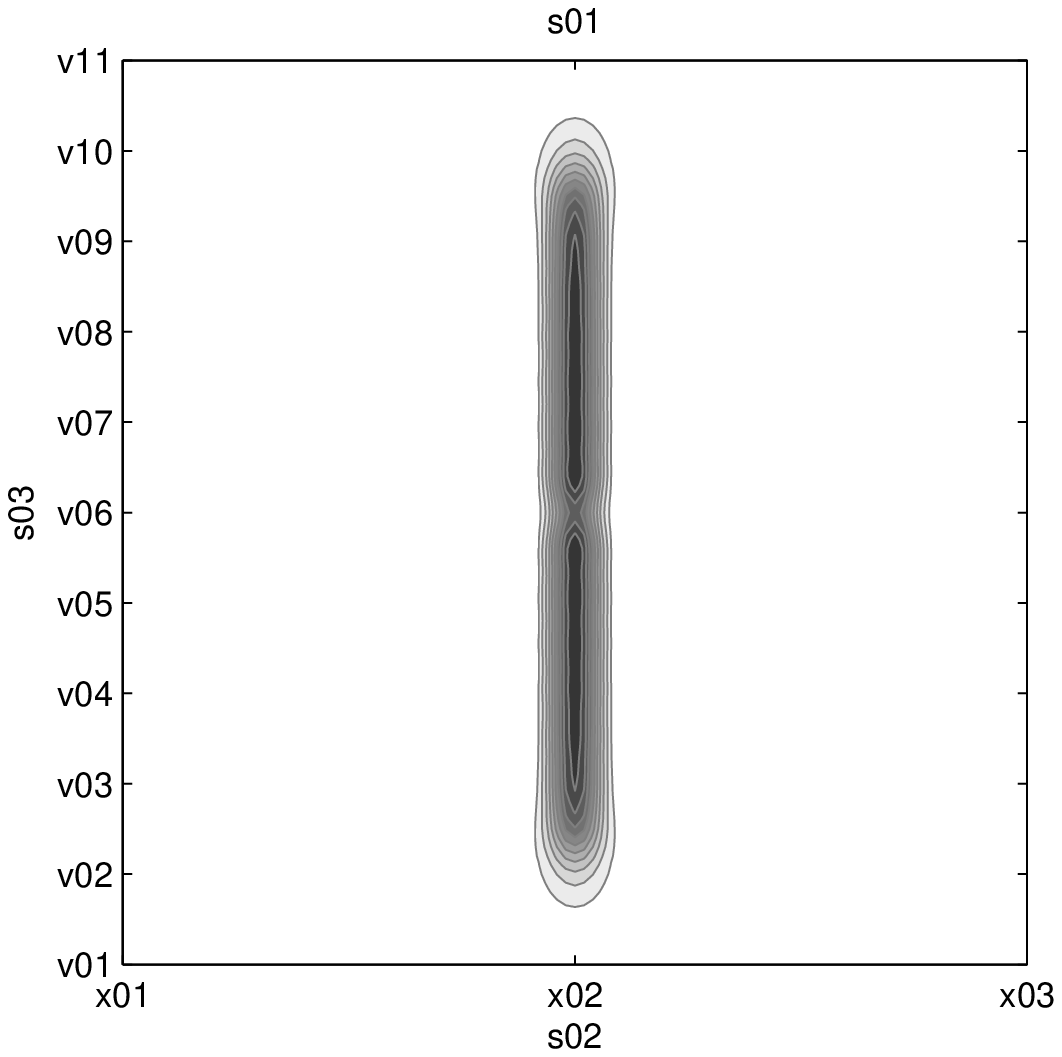}%
\end{psfrags}%
      \hspace{1cm}
\begin{psfrags}%
\psfragscanon%
%
\psfrag{s01}[b][b]{\color[rgb]{0,0,0}\setlength{\tabcolsep}{0pt}\begin{tabular}{c}$\phi_2$\end{tabular}}%
\psfrag{s02}[t][t]{\color[rgb]{0,0,0}\setlength{\tabcolsep}{0pt}\begin{tabular}{c}$x$\end{tabular}}%
\psfrag{s03}[b][b]{\color[rgb]{0,0,0}\setlength{\tabcolsep}{0pt}\begin{tabular}{c}$y$\end{tabular}}%
%
\psfrag{x01}[t][t]{-5}%
\psfrag{x02}[t][t]{0}%
\psfrag{x03}[t][t]{5}%
%
\psfrag{v01}[r][r]{-5}%
\psfrag{v02}[r][r]{-4}%
\psfrag{v03}[r][r]{-3}%
\psfrag{v04}[r][r]{-2}%
\psfrag{v05}[r][r]{-1}%
\psfrag{v06}[r][r]{0}%
\psfrag{v07}[r][r]{1}%
\psfrag{v08}[r][r]{2}%
\psfrag{v09}[r][r]{3}%
\psfrag{v10}[r][r]{4}%
\psfrag{v11}[r][r]{5}%
%
\psfrag{z01}[r][r]{0}%
\psfrag{z02}[r][r]{0.1}%
\psfrag{z03}[r][r]{0.2}%
\psfrag{z04}[r][r]{0.3}%
\psfrag{z05}[r][r]{0.4}%
\psfrag{z06}[r][r]{0.5}%
\psfrag{z07}[r][r]{0.6}%
\psfrag{z08}[r][r]{0.7}%
%
\includegraphics[scale=\figurescale]{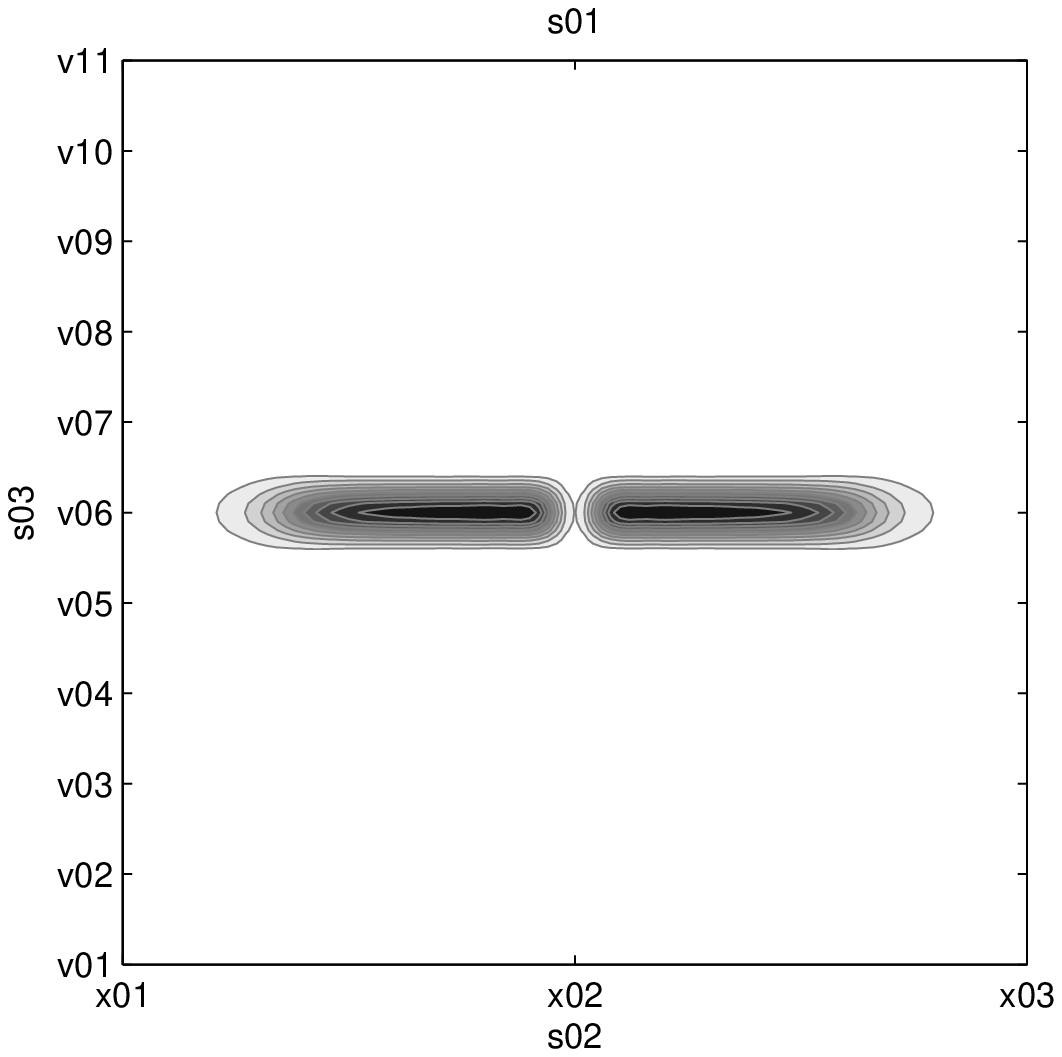}%
\end{psfrags}%
\caption{2D contour plots, in the square $[-5,5]^2$, of the ground state solution for $N_1=N_2=m_1=m_2=1$,
$\theta_{11}=400$, $\theta_{22}=150$, $\theta_{12}=1200$, $\omega_{11}=100$, $\omega_{22}=100$,
$\omega_{12}=\omega_{21}=1$ and potentials centered at the origin. The 
segregation around the origin is evident.}\label{GS2sigaro-bis}
\end{center}
\end{figure}

Depending on the relative magnitude of parameters $\omega_{ij}$, there 
are some directions along which the ground state solutions tends to 
concentrate. For instance, for $\omega_{11}$ (resp.\ $\omega_{12}$) much larger
than $\omega_{12}$ (resp.\ $\omega_{11}$), the component $\phi_1$ has a cigar-like shape along
the $y$-axis (resp.\ $x$-axis). Similar behaviour for $\phi_2$ along the $y$-axis 
(resp.\ $x$-axis) for $\omega_{21}$ (resp.\ $\omega_{22}$) much larger
than $\omega_{22}$ (resp.\ $\omega_{21}$). In Figure \ref{GS2sigaro} we consider the
small interaction case, namely  $\theta_{12}\ll \theta_{ii}$, when $\omega_{ii}=100$ and
$\omega_{ij}=1$ for $i\neq j$. As it is evident from Figure \ref{GS2sigaro-bis}, increasing the inter-specific coupling constant
($\theta_{12}=1200$) the wave functions $\phi_1$ and $\phi_2$ spatially segregate around the origin.

\medskip
\subsection{Excited state solutions}

\begin{figure}[h!!!]
\begin{center}
\begin{psfrags}%
\psfragscanon%
%
\psfrag{s01}[b][b]{\color[rgb]{0,0,0}\setlength{\tabcolsep}{0pt}\begin{tabular}{c}$\phi_1$\end{tabular}}%
\psfrag{s02}[t][t]{\color[rgb]{0,0,0}\setlength{\tabcolsep}{0pt}\begin{tabular}{c}$x$\end{tabular}}%
\psfrag{s03}[b][b]{\color[rgb]{0,0,0}\setlength{\tabcolsep}{0pt}\begin{tabular}{c}$y$\end{tabular}}%
%
\psfrag{x01}[t][t]{-5}%
\psfrag{x02}[t][t]{0}%
\psfrag{x03}[t][t]{5}%
%
\psfrag{v01}[r][r]{-5}%
\psfrag{v02}[r][r]{-4}%
\psfrag{v03}[r][r]{-3}%
\psfrag{v04}[r][r]{-2}%
\psfrag{v05}[r][r]{-1}%
\psfrag{v06}[r][r]{0}%
\psfrag{v07}[r][r]{1}%
\psfrag{v08}[r][r]{2}%
\psfrag{v09}[r][r]{3}%
\psfrag{v10}[r][r]{4}%
\psfrag{v11}[r][r]{5}%
%
\psfrag{z01}[r][r]{-0.2}%
\psfrag{z02}[r][r]{-0.1}%
\psfrag{z03}[r][r]{0}%
\psfrag{z04}[r][r]{0.1}%
\psfrag{z05}[r][r]{0.2}%
\psfrag{z06}[r][r]{0.3}%
\psfrag{z07}[r][r]{0.4}%
%
\includegraphics[scale=\figurescale]{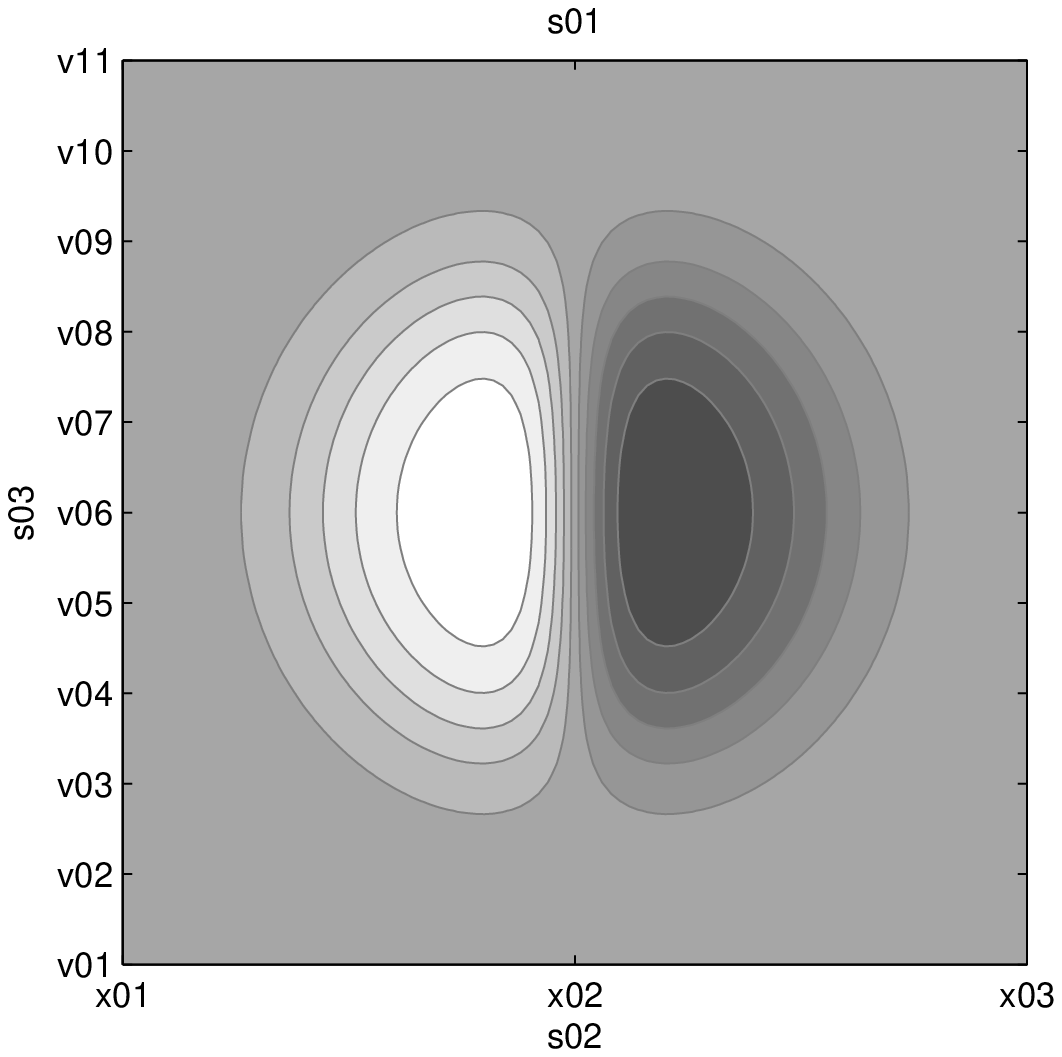}%
\end{psfrags}%
      \hspace{1cm}
\begin{psfrags}%
\psfragscanon%
%
\psfrag{s01}[b][b]{\color[rgb]{0,0,0}\setlength{\tabcolsep}{0pt}\begin{tabular}{c}$\phi_2$\end{tabular}}%
\psfrag{s02}[t][t]{\color[rgb]{0,0,0}\setlength{\tabcolsep}{0pt}\begin{tabular}{c}$x$\end{tabular}}%
\psfrag{s03}[b][b]{\color[rgb]{0,0,0}\setlength{\tabcolsep}{0pt}\begin{tabular}{c}$y$\end{tabular}}%
%
\psfrag{x01}[t][t]{-5}%
\psfrag{x02}[t][t]{0}%
\psfrag{x03}[t][t]{5}%
%
\psfrag{v01}[r][r]{-5}%
\psfrag{v02}[r][r]{-4}%
\psfrag{v03}[r][r]{-3}%
\psfrag{v04}[r][r]{-2}%
\psfrag{v05}[r][r]{-1}%
\psfrag{v06}[r][r]{0}%
\psfrag{v07}[r][r]{1}%
\psfrag{v08}[r][r]{2}%
\psfrag{v09}[r][r]{3}%
\psfrag{v10}[r][r]{4}%
\psfrag{v11}[r][r]{5}%
%
\psfrag{z01}[r][r]{-0.2}%
\psfrag{z02}[r][r]{-0.1}%
\psfrag{z03}[r][r]{0}%
\psfrag{z04}[r][r]{0.1}%
\psfrag{z05}[r][r]{0.2}%
\psfrag{z06}[r][r]{0.3}%
\psfrag{z07}[r][r]{0.4}%
%
\includegraphics[scale=\figurescale]{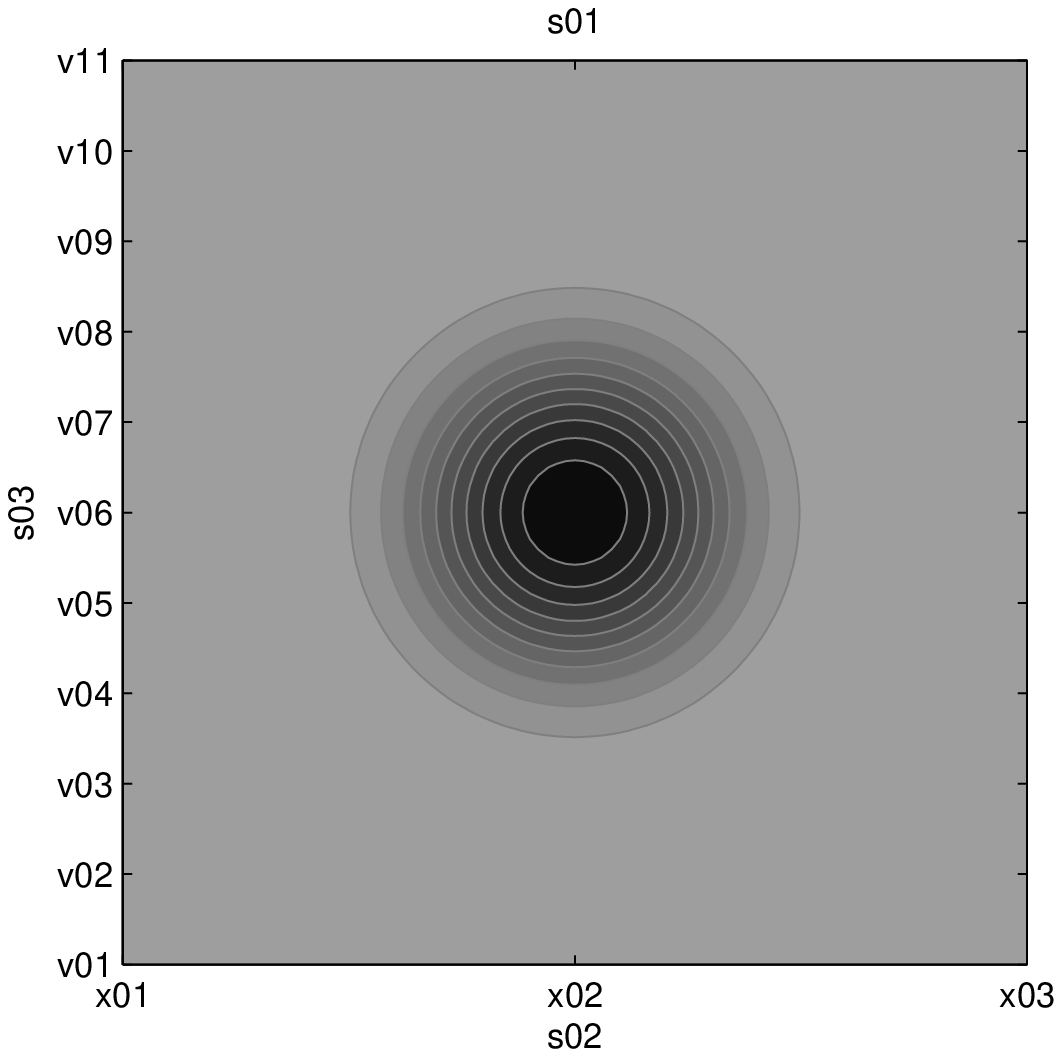}\\[1cm]
\end{psfrags}%
      \begin{psfrags}%
\psfragscanon%
%
\psfrag{s01}[b][b]{\color[rgb]{0,0,0}\setlength{\tabcolsep}{0pt}\begin{tabular}{c}$\phi_1$\end{tabular}}%
\psfrag{s02}[t][t]{\color[rgb]{0,0,0}\setlength{\tabcolsep}{0pt}\begin{tabular}{c}$x$\end{tabular}}%
\psfrag{s03}[b][b]{\color[rgb]{0,0,0}\setlength{\tabcolsep}{0pt}\begin{tabular}{c}$y$\end{tabular}}%
%
\psfrag{x01}[t][t]{-5}%
\psfrag{x02}[t][t]{0}%
\psfrag{x03}[t][t]{5}%
%
\psfrag{v01}[r][r]{-5}%
\psfrag{v02}[r][r]{-4}%
\psfrag{v03}[r][r]{-3}%
\psfrag{v04}[r][r]{-2}%
\psfrag{v05}[r][r]{-1}%
\psfrag{v06}[r][r]{0}%
\psfrag{v07}[r][r]{1}%
\psfrag{v08}[r][r]{2}%
\psfrag{v09}[r][r]{3}%
\psfrag{v10}[r][r]{4}%
\psfrag{v11}[r][r]{5}%
%
\psfrag{z01}[r][r]{-0.2}%
\psfrag{z02}[r][r]{-0.1}%
\psfrag{z03}[r][r]{0}%
\psfrag{z04}[r][r]{0.1}%
\psfrag{z05}[r][r]{0.2}%
\psfrag{z06}[r][r]{0.3}%
\psfrag{z07}[r][r]{0.4}%
\psfrag{z08}[r][r]{0.5}%
%
\includegraphics[scale=\figurescale]{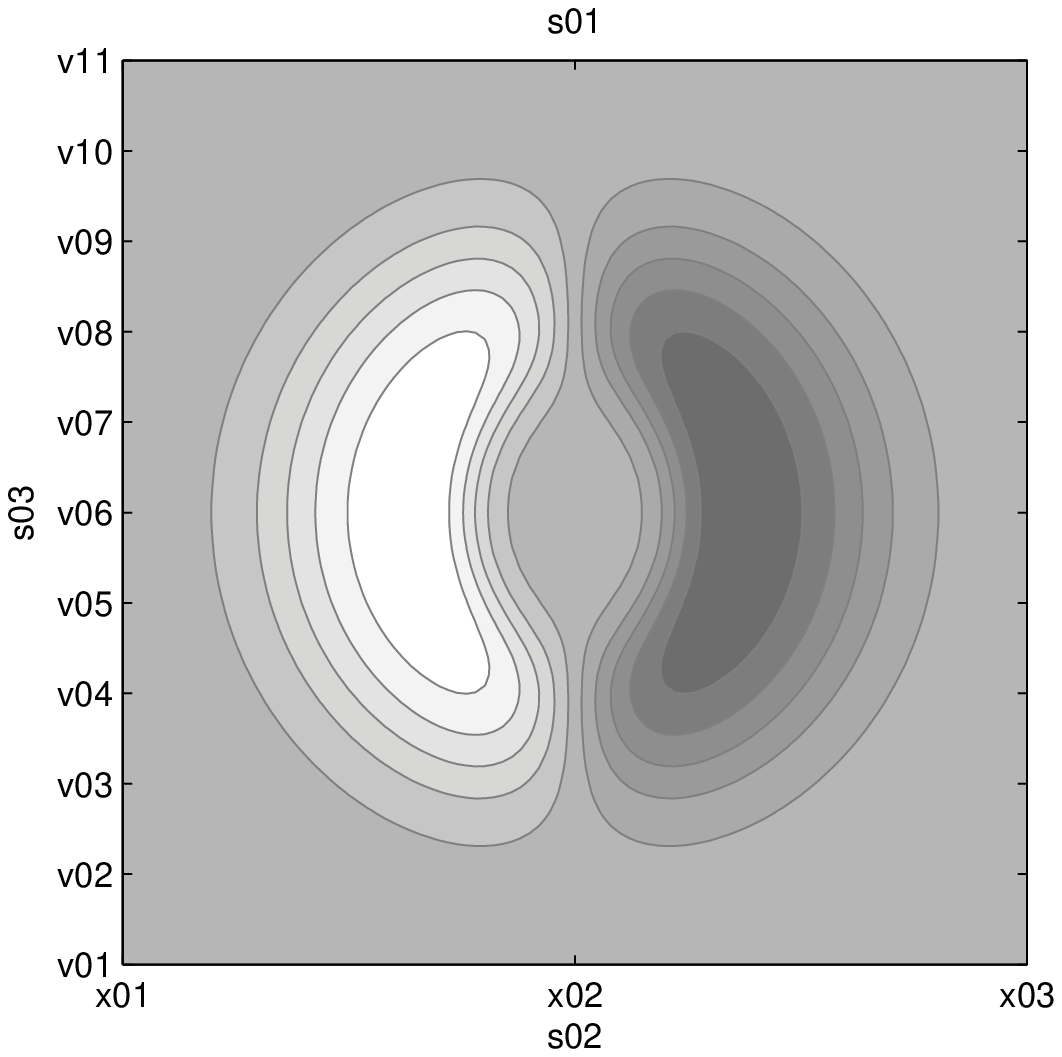}%
\end{psfrags}%
      \hspace{1cm}
\begin{psfrags}%
\psfragscanon%
%
\psfrag{s01}[b][b]{\color[rgb]{0,0,0}\setlength{\tabcolsep}{0pt}\begin{tabular}{c}$\phi_2$\end{tabular}}%
\psfrag{s02}[t][t]{\color[rgb]{0,0,0}\setlength{\tabcolsep}{0pt}\begin{tabular}{c}$x$\end{tabular}}%
\psfrag{s03}[b][b]{\color[rgb]{0,0,0}\setlength{\tabcolsep}{0pt}\begin{tabular}{c}$y$\end{tabular}}%
%
\psfrag{x01}[t][t]{-5}%
\psfrag{x02}[t][t]{0}%
\psfrag{x03}[t][t]{5}%
%
\psfrag{v01}[r][r]{-5}%
\psfrag{v02}[r][r]{-4}%
\psfrag{v03}[r][r]{-3}%
\psfrag{v04}[r][r]{-2}%
\psfrag{v05}[r][r]{-1}%
\psfrag{v06}[r][r]{0}%
\psfrag{v07}[r][r]{1}%
\psfrag{v08}[r][r]{2}%
\psfrag{v09}[r][r]{3}%
\psfrag{v10}[r][r]{4}%
\psfrag{v11}[r][r]{5}%
%
\psfrag{z01}[r][r]{-0.2}%
\psfrag{z02}[r][r]{-0.1}%
\psfrag{z03}[r][r]{0}%
\psfrag{z04}[r][r]{0.1}%
\psfrag{z05}[r][r]{0.2}%
\psfrag{z06}[r][r]{0.3}%
\psfrag{z07}[r][r]{0.4}%
\psfrag{z08}[r][r]{0.5}%
%
\includegraphics[scale=\figurescale]{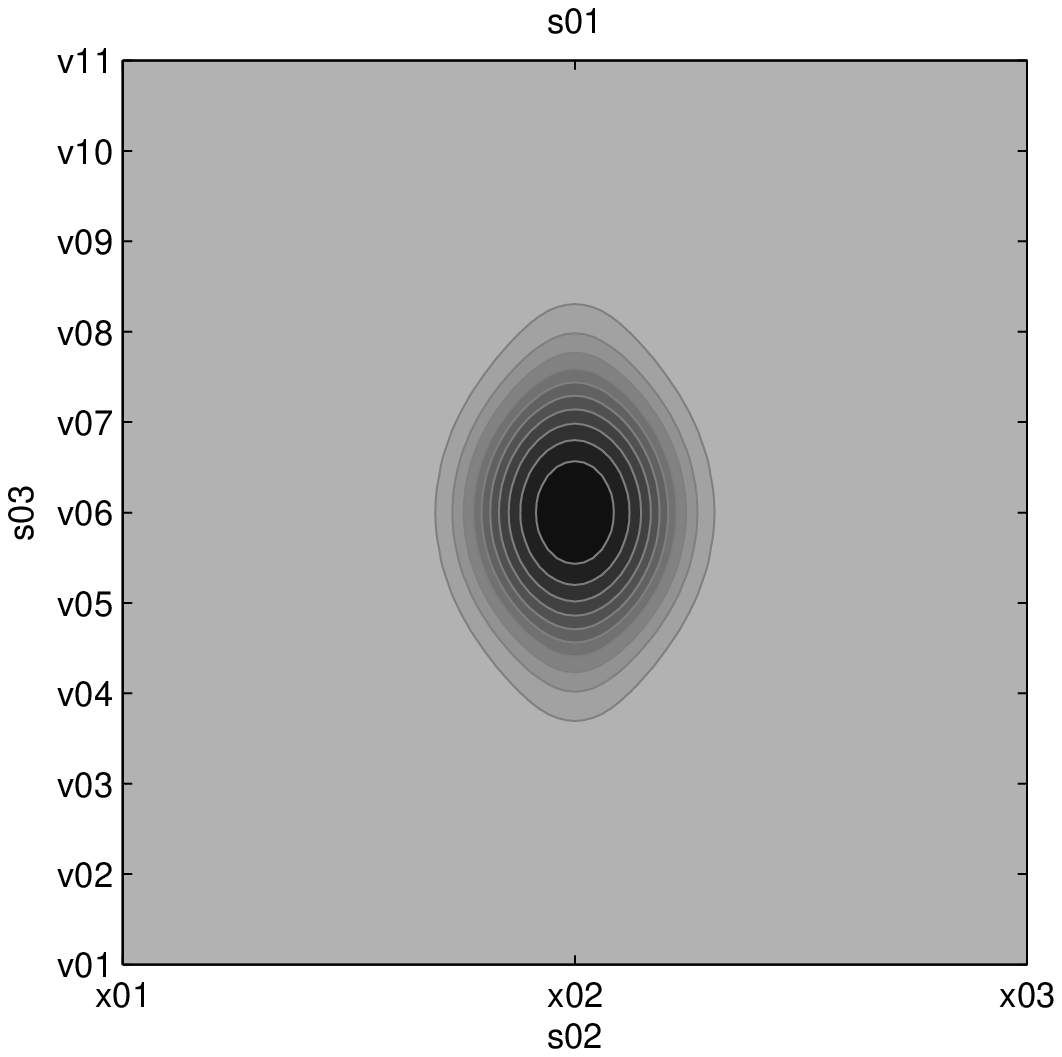}%
\end{psfrags}%
\caption{2D contour plots, in the square $[-5,5]^2$, of an excited state solution of the system.
In the left figures (two nodal regions) we have the $\phi_1$ component corresponding 
to $N_1=N_2=m_1=m_2=1$, $l_1=1$, $l_2=0$ (initial guess), 
$\theta_{11}=50$,  $\theta_{22}=5$, $\theta_{12}=0$ (top) and $\theta_{12}=120$ (bottom). 
In the right figures (no nodal regions) we have 
the corresponding $\phi_2$ component, with $l_1=l_2=0$ (initial guess).
The potentials are centered at the origin.}\label{aaaaaaa}
\end{center}
\end{figure}

\begin{figure}[h!!!]
\begin{center}
\begin{psfrags}%
\psfragscanon%
%
\psfrag{s01}[b][b]{\color[rgb]{0,0,0}\setlength{\tabcolsep}{0pt}\begin{tabular}{c}$\phi_1$\end{tabular}}%
\psfrag{s02}[t][t]{\color[rgb]{0,0,0}\setlength{\tabcolsep}{0pt}\begin{tabular}{c}$x$\end{tabular}}%
\psfrag{s03}[b][b]{\color[rgb]{0,0,0}\setlength{\tabcolsep}{0pt}\begin{tabular}{c}$y$\end{tabular}}%
%
\psfrag{x01}[t][t]{-5}%
\psfrag{x02}[t][t]{0}%
\psfrag{x03}[t][t]{5}%
%
\psfrag{v01}[r][r]{-5}%
\psfrag{v02}[r][r]{-4}%
\psfrag{v03}[r][r]{-3}%
\psfrag{v04}[r][r]{-2}%
\psfrag{v05}[r][r]{-1}%
\psfrag{v06}[r][r]{0}%
\psfrag{v07}[r][r]{1}%
\psfrag{v08}[r][r]{2}%
\psfrag{v09}[r][r]{3}%
\psfrag{v10}[r][r]{4}%
\psfrag{v11}[r][r]{5}%
%
\psfrag{z01}[r][r]{-0.3}%
\psfrag{z02}[r][r]{-0.2}%
\psfrag{z03}[r][r]{-0.1}%
\psfrag{z04}[r][r]{0}%
\psfrag{z05}[r][r]{0.1}%
\psfrag{z06}[r][r]{0.2}%
\psfrag{z07}[r][r]{0.3}%
\psfrag{z08}[r][r]{0.4}%
%
\includegraphics[scale=\figurescale]{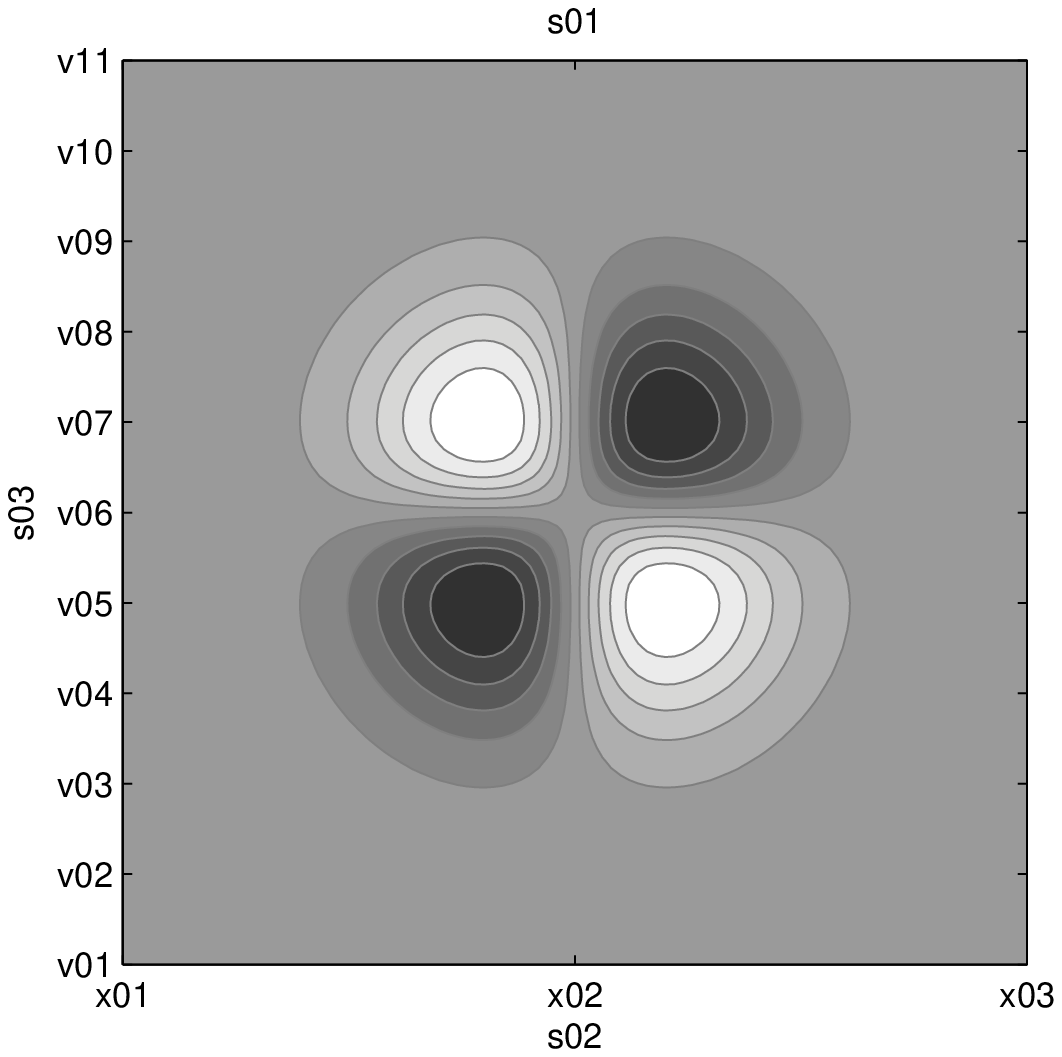}%
\end{psfrags}%
      \hspace{1cm}
\begin{psfrags}%
\psfragscanon%
%
\psfrag{s01}[b][b]{\color[rgb]{0,0,0}\setlength{\tabcolsep}{0pt}\begin{tabular}{c}$\phi_2$\end{tabular}}%
\psfrag{s02}[t][t]{\color[rgb]{0,0,0}\setlength{\tabcolsep}{0pt}\begin{tabular}{c}$x$\end{tabular}}%
\psfrag{s03}[b][b]{\color[rgb]{0,0,0}\setlength{\tabcolsep}{0pt}\begin{tabular}{c}$y$\end{tabular}}%
%
\psfrag{x01}[t][t]{-5}%
\psfrag{x02}[t][t]{0}%
\psfrag{x03}[t][t]{5}%
%
\psfrag{v01}[r][r]{-5}%
\psfrag{v02}[r][r]{-4}%
\psfrag{v03}[r][r]{-3}%
\psfrag{v04}[r][r]{-2}%
\psfrag{v05}[r][r]{-1}%
\psfrag{v06}[r][r]{0}%
\psfrag{v07}[r][r]{1}%
\psfrag{v08}[r][r]{2}%
\psfrag{v09}[r][r]{3}%
\psfrag{v10}[r][r]{4}%
\psfrag{v11}[r][r]{5}%
%
\psfrag{z01}[r][r]{-0.3}%
\psfrag{z02}[r][r]{-0.2}%
\psfrag{z03}[r][r]{-0.1}%
\psfrag{z04}[r][r]{0}%
\psfrag{z05}[r][r]{0.1}%
\psfrag{z06}[r][r]{0.2}%
\psfrag{z07}[r][r]{0.3}%
\psfrag{z08}[r][r]{0.4}%
%
\includegraphics[scale=\figurescale]{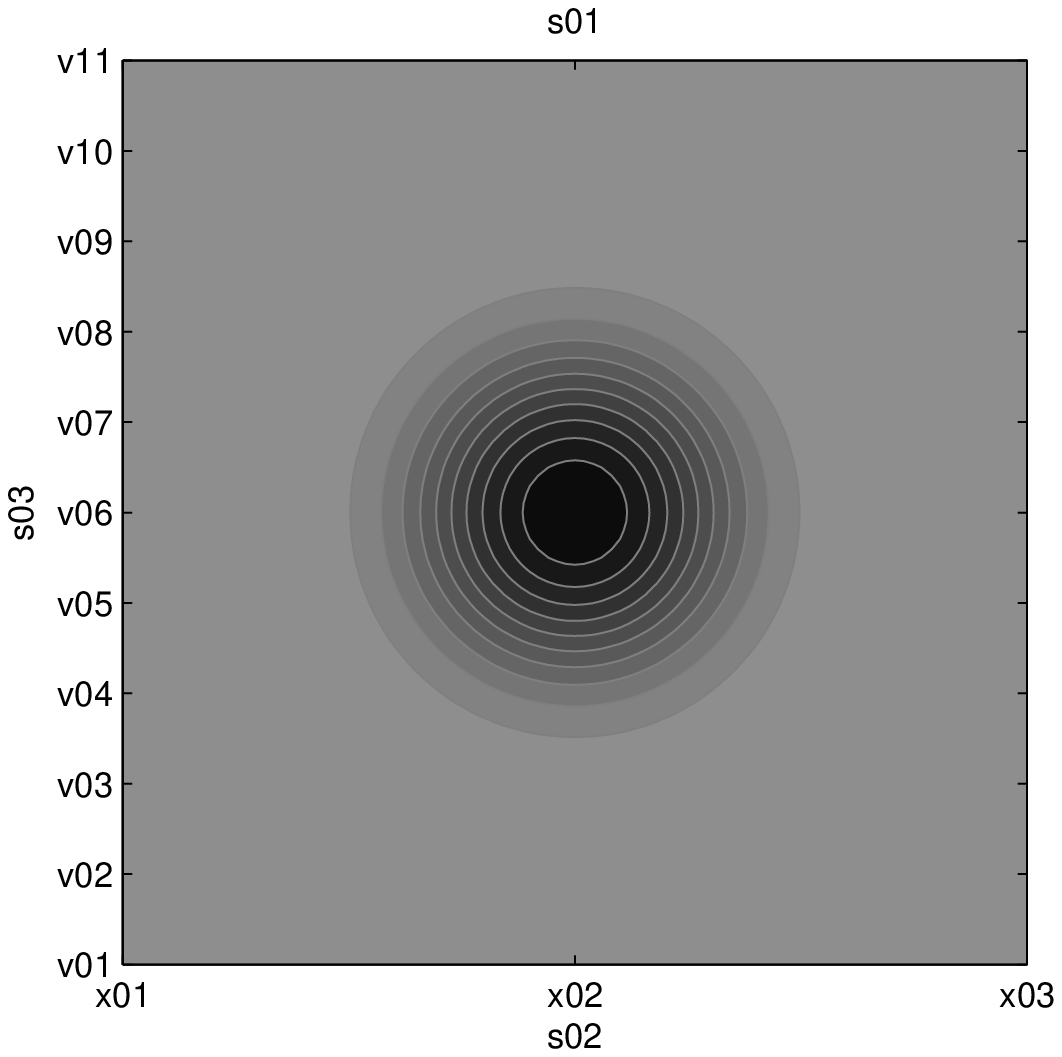}\\[1cm]%
\end{psfrags}%
\end{center}
\begin{center}
\begin{psfrags}%
\psfragscanon%
%
\psfrag{s01}[b][b]{\color[rgb]{0,0,0}\setlength{\tabcolsep}{0pt}\begin{tabular}{c}$\phi_1$\end{tabular}}%
\psfrag{s02}[t][t]{\color[rgb]{0,0,0}\setlength{\tabcolsep}{0pt}\begin{tabular}{c}$x$\end{tabular}}%
\psfrag{s03}[b][b]{\color[rgb]{0,0,0}\setlength{\tabcolsep}{0pt}\begin{tabular}{c}$y$\end{tabular}}%
%
\psfrag{x01}[t][t]{-5}%
\psfrag{x02}[t][t]{0}%
\psfrag{x03}[t][t]{5}%
%
\psfrag{v01}[r][r]{-5}%
\psfrag{v02}[r][r]{-4}%
\psfrag{v03}[r][r]{-3}%
\psfrag{v04}[r][r]{-2}%
\psfrag{v05}[r][r]{-1}%
\psfrag{v06}[r][r]{0}%
\psfrag{v07}[r][r]{1}%
\psfrag{v08}[r][r]{2}%
\psfrag{v09}[r][r]{3}%
\psfrag{v10}[r][r]{4}%
\psfrag{v11}[r][r]{5}%
%
\psfrag{z01}[r][r]{-0.3}%
\psfrag{z02}[r][r]{-0.2}%
\psfrag{z03}[r][r]{-0.1}%
\psfrag{z04}[r][r]{0}%
\psfrag{z05}[r][r]{0.1}%
\psfrag{z06}[r][r]{0.2}%
\psfrag{z07}[r][r]{0.3}%
\psfrag{z08}[r][r]{0.4}%
\psfrag{z09}[r][r]{0.5}%
%
\includegraphics[scale=\figurescale]{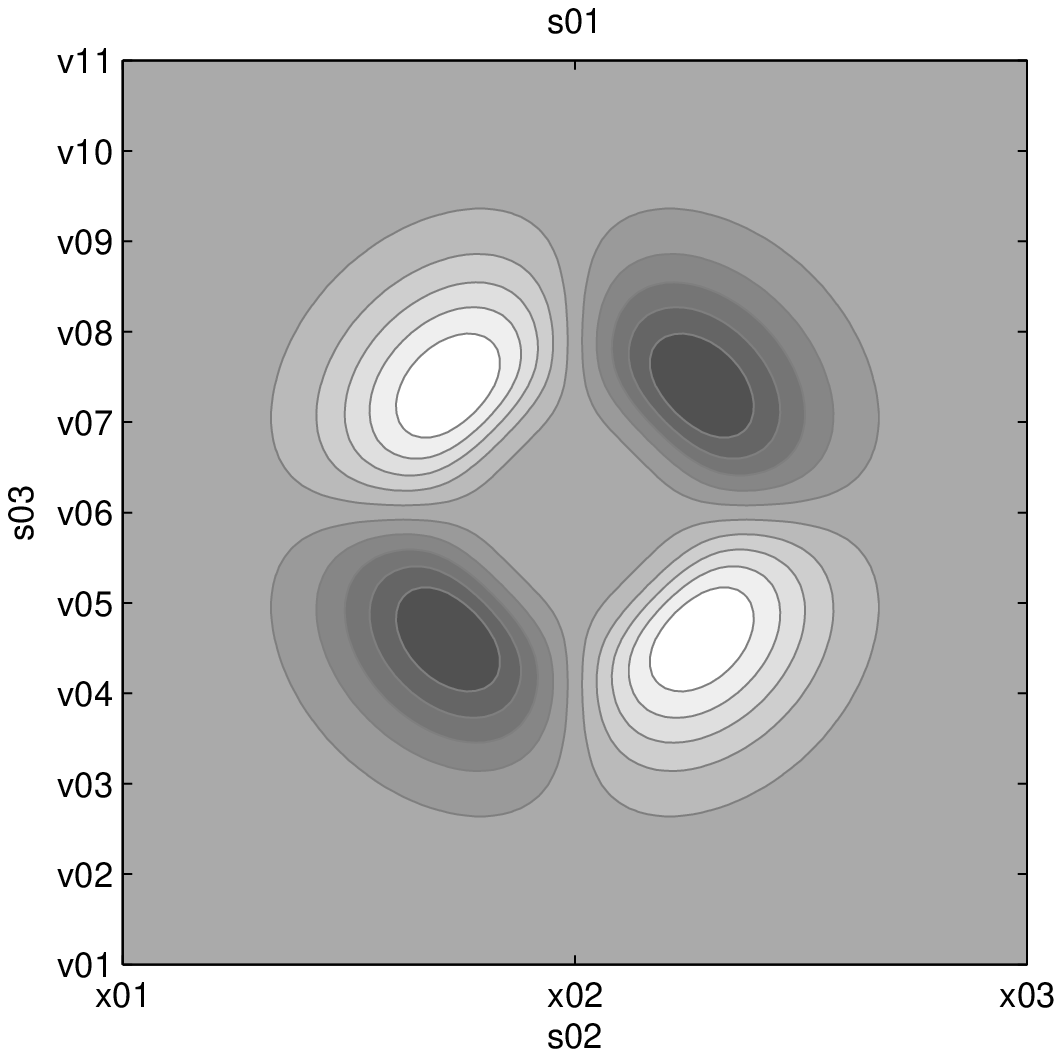}
\end{psfrags}%
      \hspace{1cm}
\begin{psfrags}%
\psfragscanon%
%
\psfrag{s01}[b][b]{\color[rgb]{0,0,0}\setlength{\tabcolsep}{0pt}\begin{tabular}{c}$\phi_2$\end{tabular}}%
\psfrag{s02}[t][t]{\color[rgb]{0,0,0}\setlength{\tabcolsep}{0pt}\begin{tabular}{c}$x$\end{tabular}}%
\psfrag{s03}[b][b]{\color[rgb]{0,0,0}\setlength{\tabcolsep}{0pt}\begin{tabular}{c}$y$\end{tabular}}%
%
\psfrag{x01}[t][t]{-5}%
\psfrag{x02}[t][t]{0}%
\psfrag{x03}[t][t]{5}%
%
\psfrag{v01}[r][r]{-5}%
\psfrag{v02}[r][r]{-4}%
\psfrag{v03}[r][r]{-3}%
\psfrag{v04}[r][r]{-2}%
\psfrag{v05}[r][r]{-1}%
\psfrag{v06}[r][r]{0}%
\psfrag{v07}[r][r]{1}%
\psfrag{v08}[r][r]{2}%
\psfrag{v09}[r][r]{3}%
\psfrag{v10}[r][r]{4}%
\psfrag{v11}[r][r]{5}%
%
\psfrag{z01}[r][r]{-0.3}%
\psfrag{z02}[r][r]{-0.2}%
\psfrag{z03}[r][r]{-0.1}%
\psfrag{z04}[r][r]{0}%
\psfrag{z05}[r][r]{0.1}%
\psfrag{z06}[r][r]{0.2}%
\psfrag{z07}[r][r]{0.3}%
\psfrag{z08}[r][r]{0.4}%
\psfrag{z09}[r][r]{0.5}%
%
\includegraphics[scale=\figurescale]{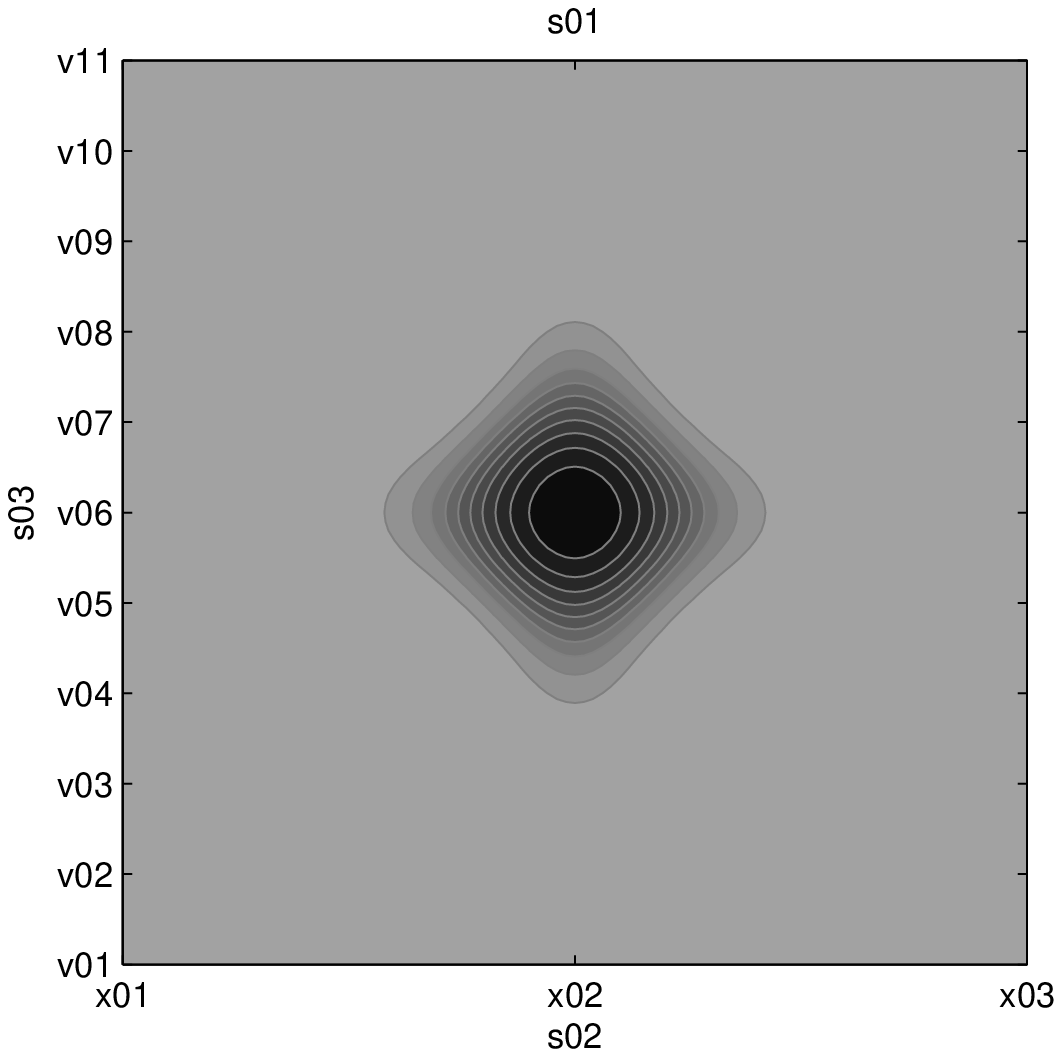}%
\end{psfrags}%
\caption{2D contour plots, in the square $[-5,5]^2$, of an excited state solution of the system.
In the left figures (four nodal regions) we have the $\phi_1$ component corresponding to $N_1=N_2=m_1=m_2=1$, $l_1=l_2=1$ (initial guess), 
$\theta_{11}=10$,  $\theta_{22}=5$, $\theta_{12}=0$ (top) and $\theta_{12}=120$ (bottom). In the right figures (no nodal regions) we have 
the corresponding $\phi_2$ component, with $l_1=l_2=0$ (initial guess).
The potentials are centered at the origin. Visibly, the supports of the $\phi_i$s segregate around the origin. 
}\label{eccitato1}
\end{center}
\end{figure}

As for the ground state solutions, in the strong interaction regime,
also higher energy solutions exhibit a phase separation behaviour. These phenomena have also been confirmed by
some numerical simulations, see e.g.\ the comparison in Figure \ref{aaaaaaa} (top, $\theta_{12}=0$
and bottom, $\theta_{12}=120$) and in Figure \ref{eccitato1} (top, $\theta_{12}=0$ and 
bottom, $\theta_{12}=120$). See also Section~\ref{Numset} for the notations.

Consider the energy functional \eqref{totenergf-Bis} defined on the space ${\mathcal H}$. If we consider
the family of all subsets $A\subset {\mathcal H}\setminus\{(0,0)\}$ which are closed and
symmetric w.r.t.\ the origin, the Krasnoselskii genus of $A\neq\emptyset$, denoted by $\gamma(A)\in\N$, is defined as the smallest
positive integer $n$ such that there exists an odd continuous function $\xi:A\to\R^n\setminus\{0\}$.
We also set $\gamma(\emptyset)=0$. When such an integer $n$ fails to exist, we put $\gamma(A)=\infty$. 
Given a positive integer $m$, we can now introduce the families ${\mathcal E}$, ${\mathcal E}_0$, $\Gamma^m$ and $\Gamma^m_0$ 
of subsets of ${\mathcal H}$, defined as follows:
\begin{align*}
{\mathcal E}&=\big\{A\subset {\mathcal H}\setminus\{(0,0)\}:\,\,\text{$A$ is closed and symmetric w.r.t.\ the origin}\big\}; \\
\noalign{\vskip3pt}
{\mathcal E}_0&=\big\{A\in{\mathcal E}:\,\,\text{if $(\phi_1,\phi_2)\in A$ then $\phi_1\phi_2=0$ a.e.\ in $\R^2$}\big\}; \\
\noalign{\vskip3pt}
\Gamma^m&=\big\{A\in {\mathcal E}:\,\,\,\text{$\gamma(A)\geq m$ and 
$\textstyle\int_{\R^2}\phi_1^2=N_1,\,\,\textstyle\int_{\R^2}\phi_2^2=N_2$ for any $(\phi_1,\phi_2)\in A$}\big\}; \\
\noalign{\vskip3pt}
\Gamma^m_0&=\big\{A\in {\mathcal E}_0:\,\text{$\gamma(A)\geq m$ and 
$\textstyle\int_{\R^2}\phi_1^2=N_1,\,\,\textstyle\int_{\R^2}\phi_2^2=N_2$ for any $(\phi_1,\phi_2)\in A$}\big\}.
\end{align*}
Then the candidates values to detect some critical (higher) levels of $E_\kappa$ are 
$$
c^m_\kappa=\inf_{A\in\Gamma^m}\sup_{(\phi_1,\phi_2)\in A} E_\kappa(\phi_1,\phi_2).
$$
We also introduce the following $\kappa$-independent values (notice that $E_\kappa|_{{\mathcal E}_0}=E_\infty$),
$$
c^m_\infty=\inf_{A\in\Gamma^m_0}\sup_{(\phi_1,\phi_2)\in A} E_\infty(\phi_1,\phi_2).
$$
Since ${\mathcal E}_0\subset {\mathcal E}$, we have $\Gamma^m_0\subset \Gamma^m$ for any $m$
and, then, by the above definitions,
\begin{equation}
\label{controlloeccitati}
c_\kappa^m\leq c_\infty^m,\quad\text{for all $m\in\N$ and $\kappa>0$}.
\end{equation}
As we prove, the levels $c^m_\kappa$ (which satisfy $c_\kappa^m\leq c_\kappa^{m+1}$ as
$\Gamma^{m+1}\subset \Gamma^{m}$ for any $m\in\N$) correspond to critical points of 
$E_\kappa$ on ${\mathcal H}$ constrained to the sphere ${\mathcal S}$,
thus yielding a sequence of nonlinear excited states of the Gross--Pitaevskii system \eqref{systemGPGen}.

\begin{proposition}
\label{finalPropES}
Let $m$ a positive integer. Then, there exists a sequence of solutions $(\phi_1^{\kappa,m},\phi_2^{\kappa,m})$
of system \eqref{systemGPGen} at energy levels  $c^m_\kappa$ such that, in the large competition limit $\kappa\to\infty$,
it converges, weakly in ${\mathcal H}$ and strongly in $L^q(\R^2)$ for any $q\geq 2$ 
to a limit configuration  $(\phi_1^{\infty,m},\phi_2^{\infty,m})\in {\mathcal S}_\infty$.
\end{proposition}

\begin{remark}
Contrary to the case of ground state solutions it seems not possible to show that the limiting
configuration $(\phi_1^{\infty,m},\phi_2^{\infty,m})$ corresponds to the energy level $c^m_\infty$
for the functional $E_\infty$ and satisfies suitable variational inequalities.
\end{remark}

In order to prove Proposition \ref{finalPropES}, we first show that, since $V_i\to\infty$ for $(x_1,x_2)\to\infty$, 
$E_\kappa$ satisfies a technical compactness condition, the Palais--Smale condition.
For the sake of completeness, we shall include a proof of this fact.

\begin{lemma}
\label{lemmapscond}
For any $\kappa>0$ the functional $E_\kappa|_{\mathcal S}$ satisfies the Palais--Smale condition, namely for any sequence
$(\phi^1_n,\phi_n^2)$ in ${\mathcal S}$ such that $E_\kappa(\phi^1_n,\phi_n^2)$ is bounded
and $dE_\kappa|_{\mathcal S}(\phi^1_n,\phi_n^2)\to 0$ as $n\to\infty$ in the dual space ${\mathcal H}^*$
of ${\mathcal H}$ (called Palais--Smale sequence) there exists a strongly convergent subsequence in ${\mathcal H}$.
\end{lemma}
\begin{proof}
Let $\kappa>0$ and let $(\phi^1_n,\phi_n^2)\subset {\mathcal S}$ be a Palais--Smale sequence for $E_\kappa$. 
In particular, 
$$
\sup_{n\geq 1}\|(\phi^1_n,\phi^2_n)\|_{{\mathcal H}}^2\leq \sup_{n\geq 1} E_\kappa(\phi^1_n,\phi^2_n)<\infty
$$
Hence $(\phi^1_n,\phi_n^2)$ is bounded in ${\mathcal H}$ and, up to a subsequence, it converges weakly in 
${\mathcal H}$, and for a.e.\ $(x_1,x_2)$ in $\R^2$, to a function $(\phi^1_\infty,\phi_\infty^2)\in {\mathcal H}$. 
Notice that ${\mathcal H}$ is compactly embedded into $L^2(\R^2)\times L^2(\R^2)$ as, for any 
$i=1,2$, we have
\begin{equation}
\label{embedkeyineq}
\sup_{n\geq 1}\sup_{R>0} R^2\int_{\R^2\setminus B_R(x_{i1},x_{i2})}\!\!\!(\phi_n^i)^2<\infty.
\end{equation}
Then, up to a further subsequence, $\phi^i_n \to \phi^i_\infty$ 
in $L^2(\R^2)$ as $n\to\infty$, which yields $(\phi^1_\infty,\phi_\infty^2)\in {\mathcal S}$. 
Hence, by the Gagliardo--Nirenberg interpolation inequality 
\begin{equation}
\label{gagliardonirenberginequr2}
\|\phi\|_{L^{\frac{2}{1-\alpha}}(\R^2)}\leq c\|\nabla\phi\|_{L^2(\R^2)}^{\alpha} \|\phi\|_{L^2(\R^2)}^{1-\alpha},\qquad
\forall\alpha\in[0,1),\,\,\,\forall \phi\in H^1(\R^2),
\end{equation}
taking, in particular, $\alpha=1/2$ we get ($c>0$ changes from inequality to inequality)
$$
\|\phi^i_n-\phi^i_\infty\|_{L^4(\R^2)}\leq c\|\nabla\phi^i_n-\nabla\phi^i_\infty\|_{L^2(\R^2)}^{1/2} 
\|\phi^i_n-\phi^i_\infty\|_{L^2(\R^2)}^{1/2}\leq c\|\phi^i_n-\phi^i_\infty\|_{L^2(\R^2)}^{1/2}
$$
so that $\phi^i_n$ converges to $\phi^i_\infty$ strongly in $L^4(\R^2)$ as $n\to\infty$ (actually in any $L^q$),  for $i=1,2$. 
Now, by virtue of the condition $dE_\kappa |_{\mathcal S}(\phi^1_n,\phi_n^2)\to 0$ as $n\to\infty$, 
there exists a sequence $(w_n)$ in ${\mathcal H}^*$ with $w_n\to 0$ in ${\mathcal H}^*$
as $n\to\infty$ and two sequences $(\mu^i_n)\subset\R$, $i=1,2$, such that, for all $(\varphi,\eta)\in {\mathcal H}$, 
\begin{align}
\label{gradzero}
\frac{1}{2m_1}\irn\nabla \phi^1_n\cdot\nabla\varphi &+\irn V_1(x_1,x_2)\phi^1_n\varphi
+\theta_{11}\irn |\phi^1_n|^2\phi^1_n\varphi
+\kappa\irn |\phi^2_n|^2\phi^1_n\varphi  \notag \\
\noalign{\vskip6pt}
+\frac{1}{2m_2}\irn\nabla \phi^2_n\cdot\nabla\eta &+\irn V_2(x_1,x_2)\phi^2_n\eta
+\kappa\irn |\phi^1_n|^2\phi^2_n\eta
+\theta_{22}\irn |\phi^2_n|^2\phi^2_n\eta \notag \\
\noalign{\vskip6pt}
& =\mu^1_n\irn\phi^1_n\varphi+\mu^2_n\irn\phi^2_n\eta+\big\langle w_n,(\varphi,\eta) \big\rangle
\end{align}
Observe that, by choosing $\varphi=\phi^1_n$  and $\eta=0$
(resp.\  $\varphi=0$ and $\eta=\phi^2_n$) and recalling that $\irn(\phi^i_n)^2=N_i$, we get a representation formula
for $\mu_n^1$ (resp.\  $\mu_n^2$). It follows that $(\mu_n^i)$ is bounded
in $\R$ so that, up to a subsequence, it converges to some positive number $\mu_\infty^i$. Finally,
choosing any arbitrary $(\varphi,0)\in {\mathcal H}$ and $(0,\eta)\in {\mathcal H}$ as test functions in
the previous identity and taking the limit as $n\to\infty$, it holds
\begin{align}
\label{solo1}
\frac{1}{2m_1}\irn\nabla \phi^1_\infty\cdot\nabla\varphi &+\irn V_1(x_1,x_2)\phi^1_\infty\varphi
+\theta_{11}\irn |\phi^1_\infty|^2\phi^1_\infty\varphi  \\
&+\kappa\irn |\phi^2_\infty|^2\phi^1_\infty\varphi=\mu^1_\infty\irn\phi^1_\infty\varphi,  \notag\\
\noalign{\vskip4pt}
\label{solo2}
\frac{1}{2m_2}\irn\nabla \phi^2_\infty\cdot\nabla\eta &+\irn V_2(x_1,x_2)\phi^2_\infty\eta
+\kappa\irn |\phi^1_\infty|^2\phi^2_\infty\eta \\
&+\theta_{22}\irn |\phi^2_\infty|^2\phi^2_\infty\eta=\mu^2_\infty\irn\phi^2_\infty\eta.    \notag
\end{align}
In particular $(\phi^1_\infty,\phi_\infty^2)\in {\mathcal H}$ is a weak solution of
\begin{equation*}
\begin{cases}
-\frac{1}{2m_1}\Delta \phi^1_\infty+V_1(x_1,x_2)\phi^1_\infty+\theta_{11}|\phi^1_\infty|^2\phi^1_\infty
+\kappa|\phi^2_\infty|^2\phi^1_\infty=\mu^1_\infty\phi^1_\infty, \\
\noalign{\vskip6pt}
-\frac{1}{2m_2}\Delta \phi^2_\infty+V_2(x_1,x_2)\phi^2_\infty
+\kappa|\phi^1_\infty|^2\phi^2_\infty+\theta_{22}|\phi^2_\infty|^2\phi^2_\infty
=\mu^2_\infty\phi^2_\infty.
\end{cases}
\end{equation*}
Now, choosing $\varphi=\phi^1_n$ and $\eta=\phi^2_n$ in \eqref{gradzero},
$\varphi=\phi^1_\infty$ in \eqref{solo1} and $\eta=\phi^2_\infty$ in \eqref{solo2},
taking into account the strong convergence of $\phi_n^i$ to $\phi^i_\infty$ in $L^4(\R^2)$, 
that $(\phi^1_n,\phi_n^2)$ is bounded in ${\mathcal H}$
and $w_n\to 0$ in ${\mathcal H}^*$, by the resulting identities we get
\begin{align*}
&  \lim_{n\to\infty}\|(\phi^1_n,\phi^2_n)\|_{{\mathcal H}}^2=\lim_{n\to\infty}\sum_{i=1}^2\frac{1}{2m_i}\irn |\nabla \phi^i_n|^2+\irn V_i(x_1,x_2)(\phi^i_n)^2 \\
&= \lim_{n\to\infty}\Big[N_1\mu^1_n+N_2\mu^2_n
-\sum_{i=1}^2\theta_{ii}\irn |\phi^i_n|^4-2\kappa\irn |\phi^1_n|^2|\phi^2_n|^2+\big\langle w_n,(\phi^1_n,\phi^2_n)\big\rangle\Big] \\
& =N_1\mu^1_\infty+N_2\mu^2_\infty-\sum_{i=1}^2\theta_{ii}\irn |\phi^i_\infty|^4-2\kappa\irn |\phi^1_\infty|^2|\phi^2_\infty|^2 \\
& =\sum_{i=1}^2\frac{1}{2m_i}\irn |\nabla \phi^i_\infty|^2+\irn V_i(x_1,x_2)(\phi^i_\infty)^2=\|(\phi^1_\infty,\phi^2_\infty)\|_{{\mathcal H}}^2,
\end{align*}
where we used the fact that $\phi_n^1\phi_n^2\to \phi^1_\infty\phi^2_\infty$ in $L^2(\R^2)$, following by
\begin{equation*}
\irn |\phi_n^1\phi_n^2-\phi_\infty^1\phi_\infty^2|^2  
\leq 2\|\phi^1_n\|_{L^4(\R^2)}^2\|\phi_n^2-\phi^2_\infty\|_{L^4(\R^2)}^2
 +2\|\phi^2_\infty\|_{L^4(\R^2)}^2\|\phi_n^1-\phi^1_\infty\|_{L^4(\R^2)}^2.
\end{equation*}
Hence $(\phi^1_n,\phi^2_n)$ converges in ${\mathcal H}$, 
proving the Palais--Smale condition.
\end{proof}

We now recall the following existence result (see e.g.\ \cite[Theorem 5.7]{struwebook}).

Let $(X,\|\cdot\|)$ be a infinite dimensional Banach space and let $Y\subset X\setminus\{0\}$
be a complete symmetric $C^{1,1}$-manifold. Let $f:Y\to\R$ be an even functional of class $C^1$. Assume that $f$ 
satisfies the Palais--Smale condition and is bounded from below on $Y$.
Then $f$ admits at least $N=\sup\{\gamma(K):\, K\subset Y\,\text{compact and symmetric}\}$ critical points.
\vskip3pt
\noindent
We are now ready to prove Proposition \ref{finalPropES}.
\vskip3pt
\noindent
{\em Proof of Proposition \ref{finalPropES}.}
Since $E_\kappa$ is a $C^1$ functional, satisfies the Palais--Smale condition by Lemma \ref{lemmapscond}, is even and bounded
from below (as $E_\kappa\geq 0$), the above mentioned result applies with $Y={\mathcal S}$ yielding (it holds $N=\infty$) a sequence 
of solutions $(\phi_1^{\kappa,m},\phi_2^{\kappa,m})$ in ${\mathcal S}$ to \eqref{systemGPGen} with 
$E_\kappa(\phi_1^{\kappa,m},\phi_2^{\kappa,m})=c_\kappa^m$, $m\geq1$ and $\kappa>0$. 
With reference to \eqref{totenergf-Bis}, by means of  \eqref{controlloeccitati} this implies that
\begin{equation*}
\kappa\irn |\phi_1^{\kappa,m}|^2|\phi_2^{\kappa,m}|^2 \leq E_\infty(\phi_1^{\kappa,m},\phi_2^{\kappa,m})
+\kappa\irn |\phi_1^{\kappa,m}|^2|\phi_2^{\kappa,m}|^2=c_\kappa^m\leq c_\infty^m,
\end{equation*}
for every $\kappa>0$. As a consequence, as $c_\infty^m$ is independent of $\kappa$, for any $m\geq1$,
\begin{equation}
\label{weaksegreg-m}
\lim_{\kappa\to\infty}\irn |\phi_1^{\kappa,m}|^2|\phi_2^{\kappa,m}|^2=0.
\end{equation}
Similarly, as $\|(\phi_1^{\kappa,m},\phi_2^{\kappa,m})\|_{\mathcal H}^2\leq c_\kappa^m\leq c_\infty^m$,
it follows that $(\phi_1^{\kappa,m},\phi_2^{\kappa,m})_{\kappa>0}$ is bounded in ${\mathcal H}$. Hence,
up to a subsequence, $(\phi_1^{\kappa,m},\phi_2^{\kappa,m})_{\kappa>0}$ weakly
converges in ${\mathcal H}$ (and strongly in $L^q(\R^2)$ for any $q\geq 2$, by 
combining formulas \eqref{embedkeyineq}-\eqref{gagliardonirenberginequr2}) to a function 
$(\phi_1^{\infty,m},\phi_2^{\infty,m})\in{\mathcal H}$. In particular,
$$
\irn|\phi_i^{\infty,m}|^2=N_i\quad\text{and}\quad \phi_1^{\infty,m}\phi_2^{\infty,m}=0\quad\text{a.e.\ in $\R^2$},
$$
which proves the assertion.
\qed

\section{Numerical computation of solutions}\label{Numset}
\begin{figure}[h!!!]
\begin{center}
\begin{psfrags}%
\psfragscanon%
%
\psfrag{s01}[b][b]{\color[rgb]{0,0,0}\setlength{\tabcolsep}{0pt}\begin{tabular}{c}${\mathcal H}_{1,0}(x_1,x_2)$\end{tabular}}%
\psfrag{s02}[t][t]{\color[rgb]{0,0,0}\setlength{\tabcolsep}{0pt}\begin{tabular}{c}$x$\end{tabular}}%
\psfrag{s03}[b][b]{\color[rgb]{0,0,0}\setlength{\tabcolsep}{0pt}\begin{tabular}{c}$y$\end{tabular}}%
%
\psfrag{x01}[t][t]{-4}%
\psfrag{x02}[t][t]{-2}%
\psfrag{x03}[t][t]{0}%
\psfrag{x04}[t][t]{2}%
\psfrag{x05}[t][t]{4}%
%
\psfrag{v01}[r][r]{-4}%
\psfrag{v02}[r][r]{-3}%
\psfrag{v03}[r][r]{-2}%
\psfrag{v04}[r][r]{-1}%
\psfrag{v05}[r][r]{0}%
\psfrag{v06}[r][r]{1}%
\psfrag{v07}[r][r]{2}%
\psfrag{v08}[r][r]{3}%
\psfrag{v09}[r][r]{4}%
%
\psfrag{z01}[r][r]{-0.4}%
\psfrag{z02}[r][r]{-0.3}%
\psfrag{z03}[r][r]{-0.2}%
\psfrag{z04}[r][r]{-0.1}%
\psfrag{z05}[r][r]{0}%
\psfrag{z06}[r][r]{0.1}%
\psfrag{z07}[r][r]{0.2}%
\psfrag{z08}[r][r]{0.3}%
\psfrag{z09}[r][r]{0.4}%
%
\includegraphics[scale=\figurescale]{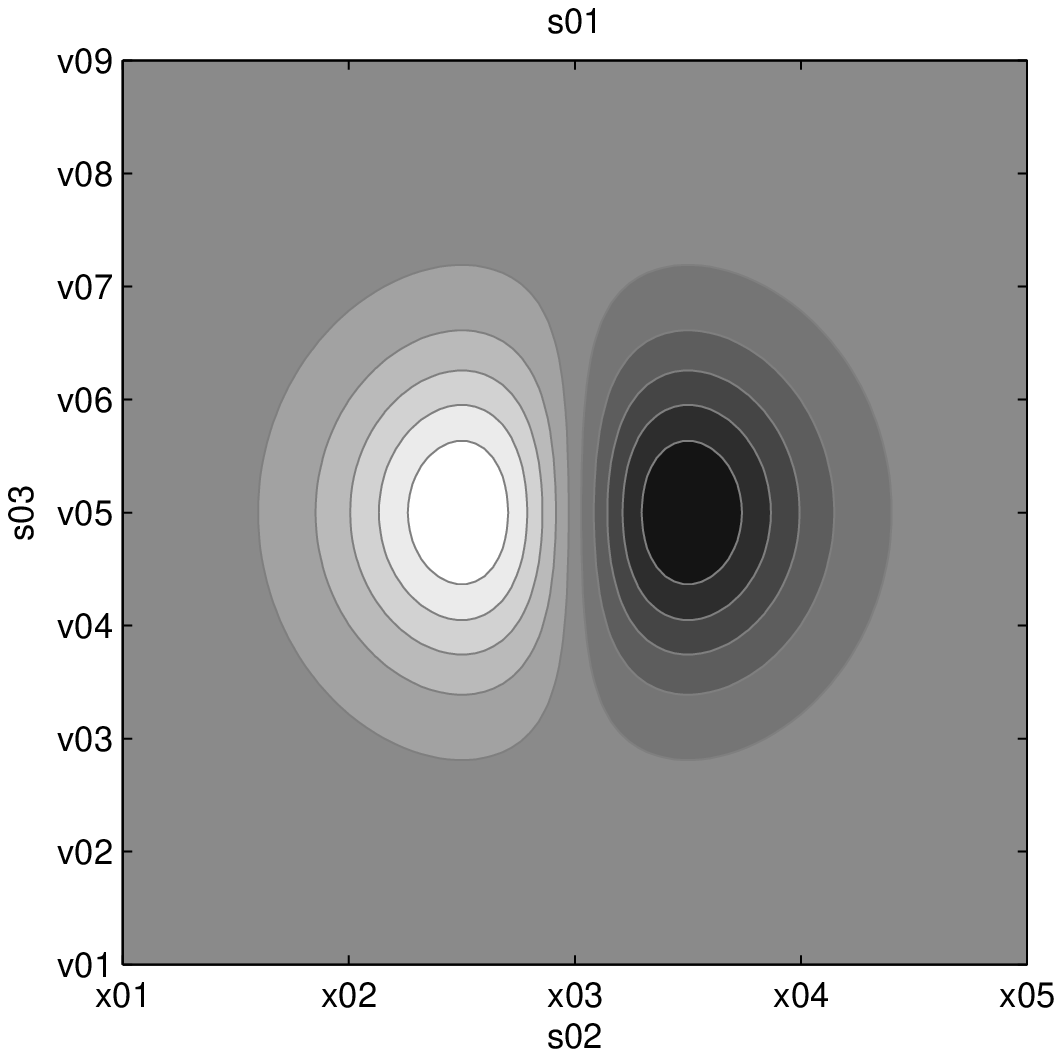}%
\end{psfrags}%
%
      \hspace{1cm}
\begin{psfrags}%
\psfragscanon%
%
\psfrag{s01}[b][b]{\color[rgb]{0,0,0}\setlength{\tabcolsep}{0pt}\begin{tabular}{c}${\mathcal H}_{1,1}(x_1,x_2)$\end{tabular}}%
\psfrag{s02}[t][t]{\color[rgb]{0,0,0}\setlength{\tabcolsep}{0pt}\begin{tabular}{c}$x$\end{tabular}}%
\psfrag{s03}[b][b]{\color[rgb]{0,0,0}\setlength{\tabcolsep}{0pt}\begin{tabular}{c}$y$\end{tabular}}%
%
\psfrag{x01}[t][t]{-4}%
\psfrag{x02}[t][t]{-2}%
\psfrag{x03}[t][t]{0}%
\psfrag{x04}[t][t]{2}%
\psfrag{x05}[t][t]{4}%
%
\psfrag{v01}[r][r]{-4}%
\psfrag{v02}[r][r]{-3}%
\psfrag{v03}[r][r]{-2}%
\psfrag{v04}[r][r]{-1}%
\psfrag{v05}[r][r]{0}%
\psfrag{v06}[r][r]{1}%
\psfrag{v07}[r][r]{2}%
\psfrag{v08}[r][r]{3}%
\psfrag{v09}[r][r]{4}%
%
\psfrag{z01}[r][r]{-0.4}%
\psfrag{z02}[r][r]{-0.3}%
\psfrag{z03}[r][r]{-0.2}%
\psfrag{z04}[r][r]{-0.1}%
\psfrag{z05}[r][r]{0}%
\psfrag{z06}[r][r]{0.1}%
\psfrag{z07}[r][r]{0.2}%
\psfrag{z08}[r][r]{0.3}%
\psfrag{z09}[r][r]{0.4}%
%
\includegraphics[scale=\figurescale]{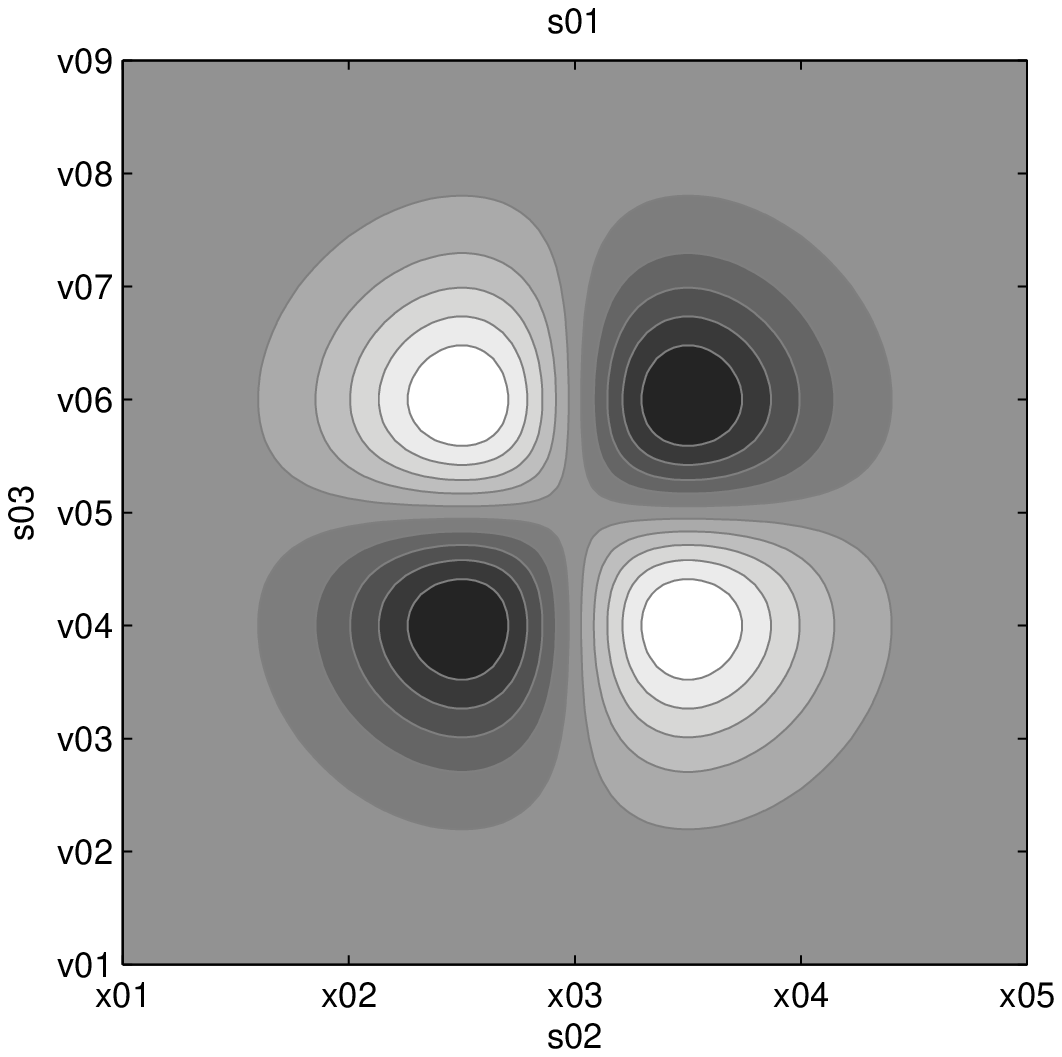}%
\end{psfrags}
\caption{2D contour plots of Hermite functions with
$l_1=1$, $l_2=0$, $\beta_1=\beta_2=1$ (left picture, one nodal region) and 
$l_1=l_2=1$, $\beta_1=\beta_2=1$ (right picture, two nodal regions).}
\label{uncoupled-excited-pattern}
\end{center}
\end{figure}

\begin{figure}[h!!!]
\begin{center}
%
\begin{psfrags}%
\psfragscanon%
%
\psfrag{s01}[b][b]{\color[rgb]{0,0,0}\setlength{\tabcolsep}{0pt}\begin{tabular}{c}${\mathcal H}_{2,1}(x_1,x_2)$\end{tabular}}%
\psfrag{s02}[t][t]{\color[rgb]{0,0,0}\setlength{\tabcolsep}{0pt}\begin{tabular}{c}$x$\end{tabular}}%
\psfrag{s03}[b][b]{\color[rgb]{0,0,0}\setlength{\tabcolsep}{0pt}\begin{tabular}{c}$y$\end{tabular}}%
%
\psfrag{x01}[t][t]{-4}%
\psfrag{x02}[t][t]{-2}%
\psfrag{x03}[t][t]{0}%
\psfrag{x04}[t][t]{2}%
\psfrag{x05}[t][t]{4}%
%
\psfrag{v01}[r][r]{-4}%
\psfrag{v02}[r][r]{-3}%
\psfrag{v03}[r][r]{-2}%
\psfrag{v04}[r][r]{-1}%
\psfrag{v05}[r][r]{0}%
\psfrag{v06}[r][r]{1}%
\psfrag{v07}[r][r]{2}%
\psfrag{v08}[r][r]{3}%
\psfrag{v09}[r][r]{4}%
%
\psfrag{z01}[r][r]{-0.3}%
\psfrag{z02}[r][r]{-0.2}%
\psfrag{z03}[r][r]{-0.1}%
\psfrag{z04}[r][r]{0}%
\psfrag{z05}[r][r]{0.1}%
\psfrag{z06}[r][r]{0.2}%
\psfrag{z07}[r][r]{0.3}%
\psfrag{z08}[r][r]{0.4}%
%
\includegraphics[scale=\figurescale]{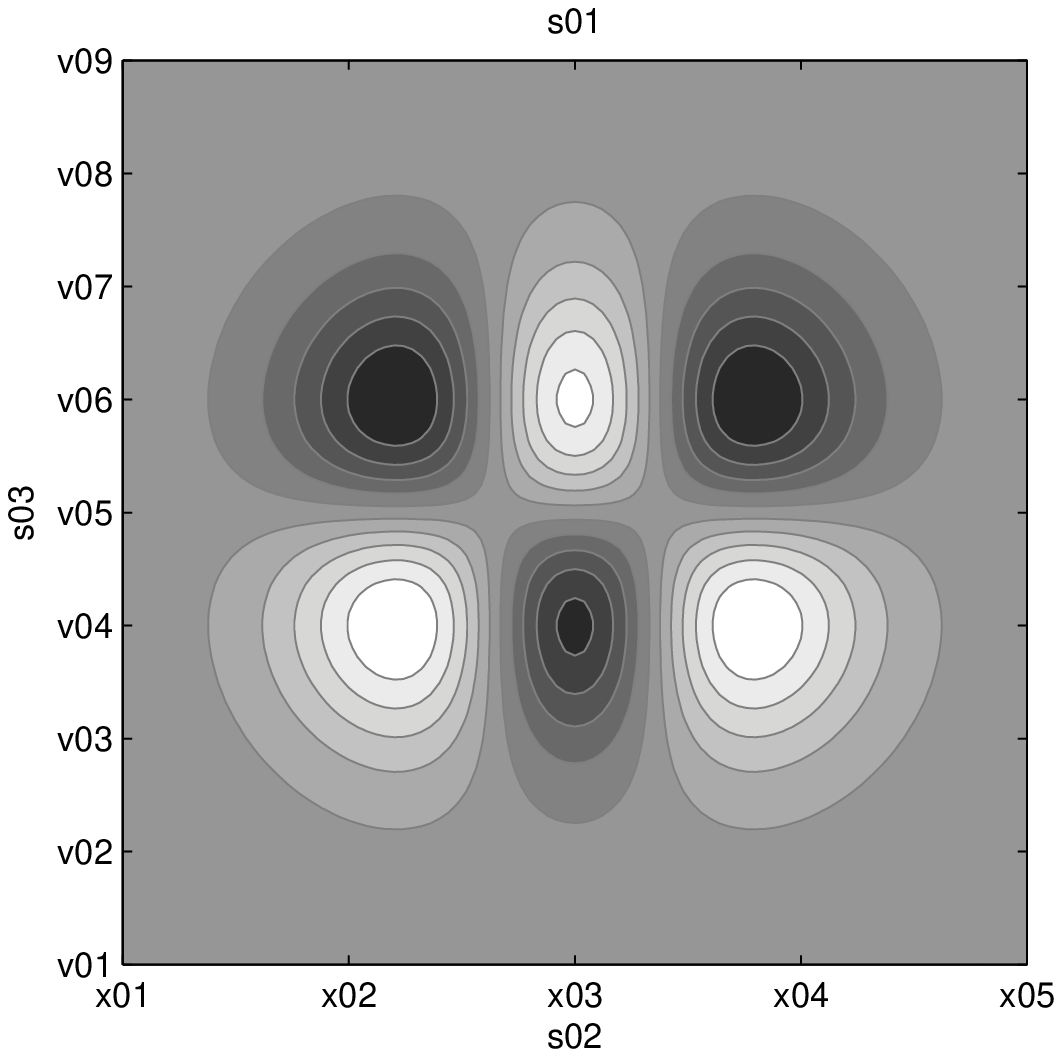}%
\end{psfrags}%
      \hspace{1cm}
\begin{psfrags}%
\psfragscanon%
%
\psfrag{s01}[b][b]{\color[rgb]{0,0,0}\setlength{\tabcolsep}{0pt}\begin{tabular}{c}${\mathcal H}_{2,2}(x_1,x_2)$\end{tabular}}%
\psfrag{s02}[t][t]{\color[rgb]{0,0,0}\setlength{\tabcolsep}{0pt}\begin{tabular}{c}$x$\end{tabular}}%
\psfrag{s03}[b][b]{\color[rgb]{0,0,0}\setlength{\tabcolsep}{0pt}\begin{tabular}{c}$y$\end{tabular}}%
%
\psfrag{x01}[t][t]{-4}%
\psfrag{x02}[t][t]{-2}%
\psfrag{x03}[t][t]{0}%
\psfrag{x04}[t][t]{2}%
\psfrag{x05}[t][t]{4}%
%
\psfrag{v01}[r][r]{-4}%
\psfrag{v02}[r][r]{-3}%
\psfrag{v03}[r][r]{-2}%
\psfrag{v04}[r][r]{-1}%
\psfrag{v05}[r][r]{0}%
\psfrag{v06}[r][r]{1}%
\psfrag{v07}[r][r]{2}%
\psfrag{v08}[r][r]{3}%
\psfrag{v09}[r][r]{4}%
%
\psfrag{z01}[r][r]{-0.3}%
\psfrag{z02}[r][r]{-0.2}%
\psfrag{z03}[r][r]{-0.1}%
\psfrag{z04}[r][r]{0}%
\psfrag{z05}[r][r]{0.1}%
\psfrag{z06}[r][r]{0.2}%
\psfrag{z07}[r][r]{0.3}%
\psfrag{z08}[r][r]{0.4}%
%
\includegraphics[scale=\figurescale]{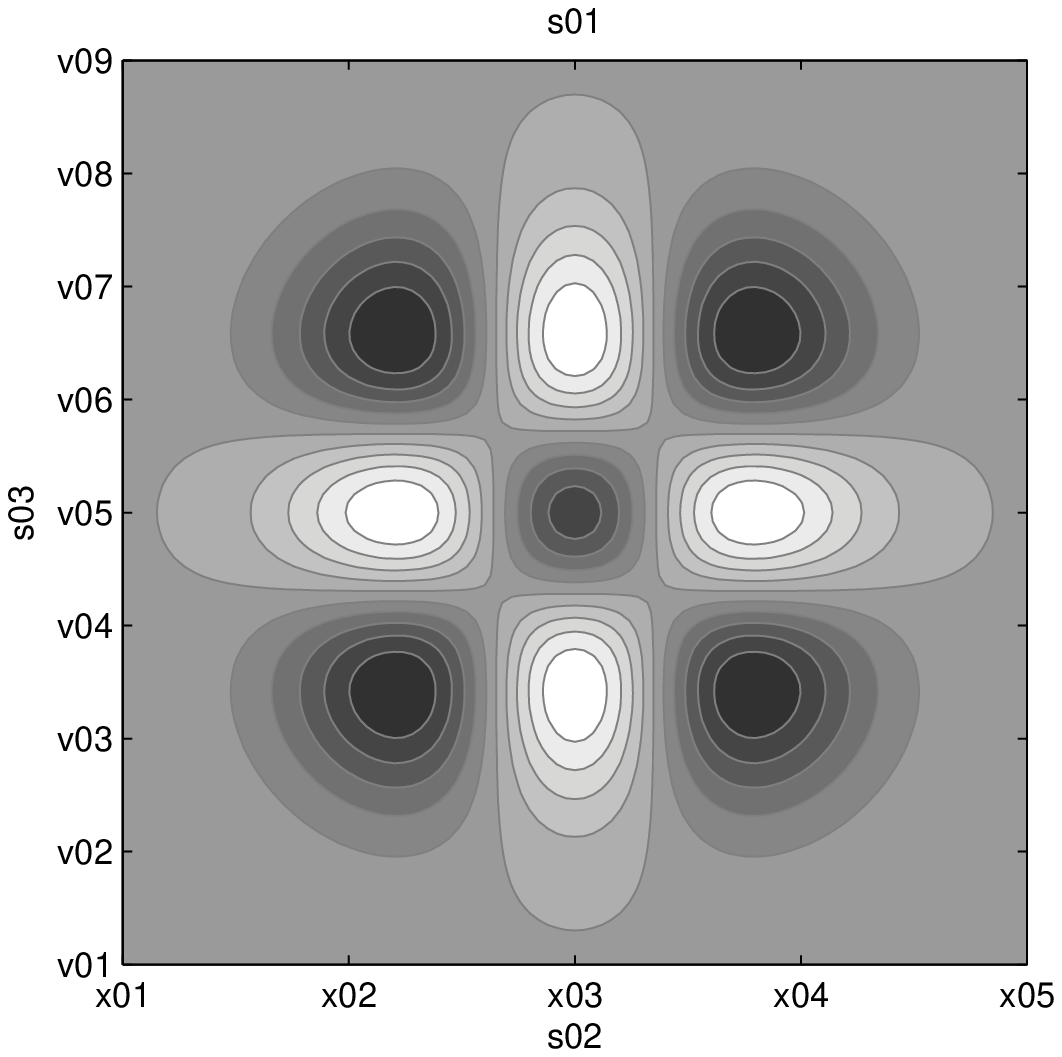}%
\end{psfrags}%
\caption{
2D contour plots of Hermite functions 
with $l_1=2$, $l_2=1$, $\beta_1=\beta_2=1$ (left picture, three nodal regions) and 
$l_1=l_2=2$, $\beta_1=\beta_2=1$ (right picture, four nodal regions).}
\label{uncoupled-excited-patternbis}
\end{center}
\end{figure}

We describe the numerical algorithm used for the computation of the ground 
states for the single one-dimensional Gross--Pitaevskii equation and we 
mention at the end of this section how the same technique can be applied to a 
system of any number of coupled equations in $\R^2$. 
Moreover, without loss of generality, we reduce to the case $\hbar=m=1$. 
The main idea is to directly minimize the energy $E(\phi)$ associated to a 
wave function $\psi(x)=e^{-i\mu t}\phi(x)$, discretized by 
Hermite functions.
As it is known, the Hermite functions $({\mathcal H}^\beta_l)_{l\in\N}$ are
defined by
$$
{\mathcal H}^\beta_l(x)=H^\beta_{l}(x)e^{-\frac12\beta^2x^2},\quad l\in\N,
$$
where $(H_{l}^\beta)_{l\in\N}$ are the {\em Hermite polynomials} 
\cite{boyd}, orthonormal in $L^2$ with respect to the weight $e^{-\beta^2x^2}$.
The Hermite functions are the solutions (ground state, for $l=0$, 
and excited states, if else) to the eigenvalue problem for the linear 
Schr\"odinger equation with standard harmonic potential
$$
\frac{1}{2}\left(-\frac{d^2}{dx^2}+(\beta^2x)^2\right){\mathcal H}_l=
\lambda_l{\mathcal H}_l,\qquad
\lambda_l=\beta^2\left(l+\frac{1}{2}\right).
$$
If we set
$$
\phi=\sum_{l\in\N}\phi_l{\mathcal H}_l,
$$
where
$$
\phi_l=(\phi,{\mathcal H}_l)_{L^2}=\int_\R \phi{\mathcal H}_l,
$$
the energy functional rewrites as
\begin{equation*}
E(\phi)=\sum_{l\in\N}\lambda_l\phi_l^2+
\int_\R \left(V(x)-\frac{(\beta^2x)^2}{2}\right)
\left(\sum_{l\in\N}\phi_l{\mathcal H}_l \right)^2+
\frac12\theta\int_\R\left(\sum_{l\in\N}\phi_l{\mathcal H}_l \right)^4,
\end{equation*}
and the chemical potential turns into
\begin{equation}
\label{frapreMu}
N\mu=E(\phi)+
\frac12\theta\int_\R\left(\sum_{l\in\N}\phi_l{\mathcal H}_l \right)^4   
\end{equation}
By minimizing $E$, under the constraint $\|\phi\|_{L^2}^2=N$, 
we look for local minima of
\begin{align*}
E(\phi;\lambda)&=E(\phi)+\lambda\left(N-\sum_{l\in\N}\phi_l^2\right)
\end{align*}
which solve the system, with $k\in\N$, 
\begin{equation*}
\begin{cases}
{\displaystyle
(\lambda_\kappa-\lambda)\phi_\kappa+
\int_\R \left(V(x)-\frac{(\beta^2x)^2}{2}\right){\mathcal H}_k
\left(\sum_{l\in\N}\phi_l{\mathcal H}_l\right)+
\theta\int_\R {\mathcal H}_k\left(\sum_{l\in\N}\phi_l{\mathcal H}_l \right)^3\!\!\!
=0},\\
{\displaystyle\sum_{l\in\N}\phi_l^2=N}.
\end{cases}
\end{equation*}
We notice that, if $\phi$ is a solution of the above system, then it is 
immediately seen, by multiplying times $\phi_k$, summing up over $k$ and using \eqref{frapreMu},
that the Lagrange multiplier $\lambda$ equals the chemical potential $\mu$.
Next, we truncate to degree $L-1$ and introduce an additional parameter
$\rho=1$ in front of the first integral (its usage will be clear later), to
obtain a corresponding truncated energy functional 
$E_L(\phi;\lambda;\rho)$, whose local minima 
solve the system, with $0\le k\le L-1$,
\begin{equation*}
\begin{cases}
\displaystyle{
(\lambda_\kappa-\lambda)\phi_\kappa+
\rho\int_\R \left(V(x)-\frac{(\beta^2x)^2}{2}\right){\mathcal H}_k
\left(\sum_{l=0}^{L-1}\phi_l{\mathcal H}_l\right)+
\theta\int_\R {\mathcal H}_k\left(\sum_{l=0}^{L-1}\phi_l{\mathcal H}_l \right)^3\!\!\!
=0},\\
{\displaystyle\sum_{l=0}^{L-1}\phi_l^2=N}.
\end{cases}
\end{equation*}
In order to approximate the integrals, we used a Gauss--Hermite quadrature 
formula with $2L-1$ nodes relative to the weight $e^{-2\beta^2x^2}$.
Using the tensor basis of the Hermite functions, i.e.
\begin{equation*}
{\mathcal H}_l(x_1,x_2)=
H_{l_1}^{\beta_1}(x_1)H_{l_2}^{\beta_2}(x_2)
e^{-\frac12(\beta_1^2x_1^2+\beta_2^2x_2^2)}
\end{equation*}
the extension to the two-dimensional case is straightforward. 
In particular, in $\R^2$, ${\mathcal H}_{0,0}(x_1,x_2)$ is the ground eigenstate and ${\mathcal H}_{l_1,l_2}(x_1,x_2)$
with any $l_1\neq 0$ or $l_2\neq 0$ is an excited eigenstate of the
Schr\"odinger equation with standard harmonic potential. 
See Figures \ref{uncoupled-excited-pattern} and \ref{uncoupled-excited-patternbis}
representing  ${\mathcal H}_{1,0}$, ${\mathcal H}_{1,1}$, ${\mathcal H}_{2,1}$ and 
${\mathcal H}_{2,2}$. For small coupling constants excited states solutions of the GPE system
look like these profiles, see e.g.\ Figures \ref{aaaaaaa} and \ref{eccitato1}.
The extension to a system of any number of equations is not difficult, too. In fact, it
is sufficient to consider the total energy of the system as the functional
to be minimized, with a normalization constraint (Lagrange multiplier)
for each wave function. 

The system is solved by a modified Newton method
with backtracking line-search, which guarantees global
convergence to the ground states. We refer to \cite{Baoal,CalThal} and, 
in particular, to \cite{CORT08} for the details. 
Here we just mention that only the diagonal part of the 
Jacobian relative to $\phi_l$ is computed, thus leading to
a dramatic reduction of the computational cost for the solution of each linear 
system.
Moreover, the
initial guess for the Newton iteration is obtained by a continuation 
technique over $\rho$ and $\theta$, starting from the ground state of
the Schr\"odinger equation with the standard harmonic potential,
which corresponds to $\rho=\theta=0$. The convergence to the excited
states is not guaranteed, although we observed numerical convergence for the
examples reported in the previous section.
Of course the case of very large values of the coefficients $\theta_{ij}$ 
can be treated within the framework of the Thomas--Fermi approximation.

\bigskip
\section*{Acknowledgments}
The authors wish to thank Prof.\ Mark Ablowitz for providing
some useful bibliographic references.

\bigskip
\bibliographystyle{amsplain}

\bigskip
\medskip

\end{document}